\documentclass[preprint]{elsarticle}

\journal{Journal of Computational Physics}

\def\d{\delta} 
 
\def\e{{\epsilon}}

\def\norm#1{\|#1\|} 
\def\wb{{\bar w}}

\usepackage{algorithm}
\usepackage{algorithmic}
\usepackage{amssymb}
\usepackage{amsmath}
\usepackage{amsthm}
\usepackage{caption}
\usepackage{comment}
\usepackage{float}
\usepackage[latin1]{inputenc}
\usepackage[pdfborder={0 0 0}, colorlinks=true, linkcolor=blue, citecolor=red]{hyperref}\hypersetup{pdfborder=0 0 0}
\usepackage{mathrsfs}
\usepackage{multicol}
\usepackage{multirow}
\usepackage{subfigure}
\usepackage{wrapfig}
\usepackage{ulem}

\usepackage{soul,xargs}
\usepackage[pdftex,dvipsnames]{xcolor}

\usepackage[colorinlistoftodos,prependcaption,textsize=tiny]{todonotes}
\newcommandx{\question}[2][1=]{\todo[linecolor=red,backgroundcolor=red!20,bordercolor=red,#1]{#2}}
\newcommandx{\change}[2][1=]{\todo[linecolor=blue,backgroundcolor=blue!25,bordercolor=blue,#1]{#2}}
\newcommandx{\info}[2][1=]{\todo[linecolor=OliveGreen,backgroundcolor=OliveGreen!25,bordercolor=OliveGreen,#1]{#2}}
\newcommandx{\improve}[2][1=]{\todo[linecolor=Plum,backgroundcolor=Plum!25,bordercolor=Plum,#1]{#2}}
\newcommandx{\thiswillnotshow}[2][1=]{\todo[disable,#1]{#2}}
\newcommandx{\answer}[2][1=]{\todo[linecolor=blue,backgroundcolor=White!25,bordercolor=Plum,#1]{#2}}

\newtheorem{remark}{Remark}
\graphicspath{{./plots/}}

\begin{document}

\newcommand{\TODO}[1]{ \fbox{\parbox{3in}{\bf TODO: #1}}}

\newcommand{\grbf}[1] {\mbox{\boldmath${#1}$\unboldmath}}
\newcommand{\gbf}[1] {\mathbf{#1}}

\newcommand{\beq} {\begin{equation}}
\newcommand{\eeq} {\end{equation}}
\newcommand{\bdm} {\begin{displaymath}}
\newcommand{\edm} {\end{displaymath}}
\newcommand{\bit}{\begin{itemize}}
\newcommand{\eit}{\end{itemize}}
\newcommand{\bde}{\begin{description}}
\newcommand{\ede}{\end{description}}
\newcommand{\bce}{\begin{center}}
\newcommand{\ece}{\end{center}}
\newcommand{\ben} {\begin{enumerate}}
\newcommand{\een} {\end{enumerate}}
\newcommand{\bea} {\begin{eqnarray}}
\newcommand{\eea} {\end{eqnarray}}
\newcommand{\barr} {\begin{array}}
\newcommand{\earr} {\end{array}}
\newcommand{\bean} {\begin{eqnarray*}}
\newcommand{\eean} {\end{eqnarray*}}
\newcommand{\edoc} {\end{document}}
\newcommand{\On}{\hspace{0.2cm} \mbox{on} \hspace{0.2cm}}
\newcommand{\In}{\hspace{0.2cm} \mbox{in} \hspace{0.2cm}}
\newcommand{\Minimize}{\hspace{0.2cm} \mbox{minimize} \hspace{0.2cm}}
\newcommand{\Maximize}{\hspace{0.2cm} \mbox{maximize} \hspace{0.2cm}}
\newcommand{\SubjectTo}{\hspace{0.2cm} \mbox{subject} \hspace{0.15cm}
            \mbox{to} \hspace{0.2cm}}

\def\etal{{\it et al.~}}

\newcommand{\point}[1] {\ensuremath{\boldsymbol{#1}}}
\newcommand{\vvec}[1] {\ensuremath{\boldsymbol{#1}}}
\renewcommand{\vec}[1] {\ensuremath{\boldsymbol{#1}}}
\newcommand{\cvec}[1] {\ensuremath{\boldsymbol{{#1}}}}
\newcommand{\ten}[1] {\ensuremath{\boldsymbol{#1}}}
\newcommand{\tenfour}[1] {\ensuremath{\boldsymbol{\mathsf{#1}}}}
\newcommand{\comp}[1]{\begin{pmatrix}{ #1 }\end{pmatrix}}
\newcommand{\vv}{\ensuremath{\vec{v}}}
\newcommand{\vr}{\left(\ensuremath{\vec{r}}\right)}
\newcommand{\att}{\left(t\right)}
\newcommand{\tvv}{\ensuremath{\tilde{\vec{v}}}}
\newcommand{\ww}{\ensuremath{\vec{w}}}
\newcommand{\GG}{\ensuremath{\ten{G}}}
\newcommand{\FF}{\ensuremath{\boldsymbol{\mathfrak{F}}}}
\newcommand{\tEE}{\ensuremath{\tilde{\vec{E}}}}
\newcommand{\nn}{\ensuremath{\vec{n}}}
\renewcommand{\SS}{\ensuremath{\ten{S}}}
\newcommand{\tSS}{\ensuremath{\tilde{\ten{S}}}}
\newcommand\tK{{\bf K}}
\newcommand{\cp}{\ensuremath{{c_p}}}
\newcommand{\cs}{\ensuremath{{c_s}}}
\newcommand{\tcs}{\ensuremath{{\tilde{c}_s}}}
\newcommand{\Vp}{\ensuremath{{V^e_p}}}
\newcommand{\Vs}{\ensuremath{{V^e_s}}}
\newcommand{\Ss}{\ensuremath{{S^e_s}}}
\newcommand{\nxnx}[1] {\ensuremath{\nn\times\left(\nn\times{#1}\right)}}

\newcommand{\dd}[2] {\ensuremath{\frac{\partial {#1}}{\partial {#2}}}}
\newcommand{\DD}[2] {\ensuremath{\frac{d {#1}}{d {#2}}}}
\newcommand{\Grad} {\ensuremath{\nabla}}
\newcommand{\Gradv}[1]{\ensuremath{\nabla_{#1}}}
\newcommand{\Div} {\ensuremath{\nabla\cdot}}
\newcommand{\Divv}[1] {\ensuremath{\nabla_{#1}\cdot}}
\newcommand{\Curl} {\ensuremath{\nabla\times}}
\newcommand{\Cten}{\ensuremath{\tenfour{C}}}

\newcommand{\eval}[2][\right]{\relax
  \ifx#1\right\relax \left.\fi#2#1\rvert}

\providecommand{\abs}[1]{\left\lvert#1\right\rvert}
\providecommand{\norm}[1]{\left\lVert#1\right\rVert}

\newcommand{\Nel} {\ensuremath{{N_T}}}
\newcommand{\Ned} {\ensuremath{{N_\text{ed}}}}
\newcommand{\De} {\ensuremath{{\mathsf{D}^e}}}
\newcommand{\Dep} {\ensuremath{{\mathsf{D}^{e'}}}}
\newcommand{\Dhat} {\ensuremath{\hat{\mathsf{D}}}}
\newcommand{\Dehat} {\ensuremath{{\hat{\mathsf{D}}^e}}}
\newcommand{\half} {\ensuremath{\frac{1}{2}}}
\newcommand{\IN} {\ensuremath{{\mathcal{I}_N}}}
\newcommand{\INe} {\ensuremath{{\mathcal{I}_{N_e}}}}

\newcommand{\mc}[1]{\mathcal{#1}}
\newcommand{\mcC}[1]{\mathcal{#1}}
\newcommand{\mb}[1]{\mathbf{#1}}
\newcommand{\mbb}[1]{\mathbb{#1}}
\newcommand{\mi}[1]{\mathit{#1}}
\newcommand{\mf}[1]{\mathfrak{#1}}
\newcommand{\Oi}[1]{\int_{\mc{D}} #1 \, d\vec{x}}
\newcommand{\TOi}[1]{\int_{0}^T \int_{\mc{D}} #1 \, d\vec{x}dt}
\newcommand{\Ti}[1]{\int_{0}^T #1 \, dt}
\newcommand{\TDei}[1]{\int_{0}^T\int_{\De} #1 \, d\vec{x}dt}
\newcommand{\TDeid}[1]{\int_{0}^T\int_{\De,N} #1 \, d\vec{x}dt}
\newcommand{\Dei}[1]{\int_{\De} #1 \, d\vec{x}}
\newcommand{\DeI}[1]{\int_{D} #1 \, d\vec{x}}
\newcommand{\DeiM}[1]{\int_{\Dhat} #1 \, d\vec{r}}
\newcommand{\DeMd}[1]{\int_{\De,N_e} #1 \, d\vec{x}}
\newcommand{\DeiMd}[1]{\int_{\Dhat,N_e} #1 \, d\vec{r}}
\newcommand{\DeMdN}[1]{\int_{D,N} #1 \, d\vec{x}}
\newcommand{\DeiMdN}[1]{\int_{\Dhat,N} #1 \, d\vec{r}}
\newcommand{\DeiMdNC}[2]{\int_{\Dhat,N_{#2}} #1 \, d\vec{r}}
\newcommand{\pDei}[1]{\int_{\partial \De} #1 \, d\vec{x}}
\newcommand{\pDeiM}[1]{\int_{\partial \Dhat} #1 \, d\vec{r}}
\newcommand{\pDeiMd}[1]{\int_{\partial \Dhat,N_e} #1 \, d\vec{r}}
\newcommand{\pDeiMdN}[1]{\int_{\partial \Dhat,N} #1 \, d\vec{r}}
\newcommand{\pDeiMdNC}[2]{\int_{\partial \Dhat,N_{#2}} #1 \, d\vec{r}}
\newcommand{\pMortar}[2]{\int_{\mc{M}_{#2}} #1 \, d\vec{r}}
\newcommand{\pMortard}[2]{\int_{\mc{M}^{#2},N_{#2}} #1 \, d\vec{r}}
\newcommand{\pMortardd}[3]{\int_{\mc{M}^{#2},N_{#3}} #1 \, d\vec{r}}
\newcommand{\pMortardh}[3]{\int_{\mc{M}_{#2},N_{#3}} #1 \, d\vec{r}}
\newcommand{\pMortardhn}[4]{\int_{{#2}_{#3},N_{#4}} #1 \, d\vec{r}}
\newcommand{\TpDei}[1]{\int_{0}^T\int_{\partial \De} #1 \, d\vec{x}}
\newcommand{\TpDeid}[1]{\int_{0}^T\int_{\partial \De,N} #1 \, d\vec{x}}
\newcommand{\Deid}[1]{\int_{\De,N_e} #1 \, d\vec{x}}
\newcommand{\DeidN}[1]{\int_{\De,N} #1 \, d\vec{x}}
\newcommand{\pDeidN}[1]{\int_{\partial \De,N} #1 \, d\vec{x}}
\newcommand{\pDeid}[1]{\int_{\partial \De,N_e} #1 \, d\vec{x}}
\newcommand{\nor}[1]{\left\| #1 \right\|}
\newcommand{\snor}[1]{\left| #1 \right|}
\newcommand{\LRp}[1]{\left( #1 \right)}
\newcommand{\LRs}[1]{\left[ #1 \right]}
\newcommand{\LRa}[1]{\left< #1 \right>}
\newcommand{\LRc}[1]{\left\{ #1 \right\}}
\newcommand{\pp}[2]{\frac{\partial #1}{\partial #2}}
\newcommand{\ppcp}[1]{\frac{\partial #1}{\partial \vec{c_p}}}
\newcommand{\ppcpj}[1]{\frac{\partial #1}{\partial c_p}}
\newcommand{\ppcpmf}[1]{\frac{\partial #1}{\partial \vec{\mf{c_p}}}}
\newcommand{\ppm}[1]{\frac{\partial #1}{\partial \vec{m}}}
\newcommand{\ppho}[3]{\frac{\partial^{#3} #1}{{\partial #2}^{#3}}}
\newcommand{\ppmi}[1]{\frac{\partial #1}{\partial m_i}}
\newcommand{\ppmj}[1]{\frac{\partial #1}{\partial m_j}}
\newcommand{\ppcpm}[1]{\frac{\partial #1}{\partial \vec{c_p}^-}}
\newcommand{\ppcpp}[1]{\frac{\partial #1}{\partial \vec{c_p}^+}}
\newcommand{\ppcs}[1]{\frac{\partial #1}{\partial \vec{c_s}}}
\newcommand{\ppcsm}[1]{\frac{\partial #1}{\partial \vec{c_s}^-}}
\newcommand{\ppcsp}[1]{\frac{\partial #1}{\partial \vec{c_s}^+}}
\newcommand{\td}[2]{\frac{d #1}{d #2}}
\newcommand{\Rvert}[1]{\left. #1 \right|}
\newcommand{\bs}[1]{\boldsymbol{#1}}
\newcommand{\ipDed}[3]{\LRp{ #1, #2}_{\De,N}^{#3}}

\newcommand{\bigO}{\mathcal{O}}

\newcommand{\iOm}[1]{\int_{\Omega} #1 \, d\Omega}
\newcommand{\iGb}[1]{\int_{\Gamma_{b}} #1 \, ds}
\newcommand{\iGs}[1]{\int_{\Gamma_{s}} #1 \, ds}
\newcommand{\iGsy}[1]{\int_{\Gamma_{s}} #1 \, ds(\mb{y})}
\newcommand{\RE}{\mbox{{I}}\hspace{-0.5ex}\mbox{{R}}}
\newcommand{\iC}[1]{\int_{0}^{2\pi} #1 \, d\theta}
\newcommand{\iGr}[1]{\int_{\Gamma_{\infty}} #1 \, ds}

\newcommand{\jump}[1] {\ensuremath{[\![{#1}]\!]}}
\newcommand{\diff}[1] {\ensuremath{\LRs{#1}}}
\newcommand{\average}[1] {\ensuremath{\LRc{\!\!\LRc{#1}\!\!}}}

\newcounter{tablecell}

\newcommand*{\numbercell}{%
  \stepcounter{tablecell}%
  \thetablecell. 
}

\newtheorem{propo}{Proposition}
\newtheorem{coro}{Corollary}
\newtheorem{theo}{Theorem}
\newtheorem{lem}{Lemma}

\newtheorem{defi}{definition}

\newtheorem{rema}{remark}

\newcommand{\figlab}[1]{\label{fig:#1}}
\newcommand{\eqnlab}[1]{\label{eq:#1}}
\newcommand{\theolab}[1]{\label{theo:#1}}
\newcommand{\corolab}[1]{\label{coro:#1}}
\newcommand{\propolab}[1]{\label{propo:#1}}
\newcommand{\lemlab}[1]{\label{lem:#1}}
\newcommand{\defilab}[1]{\label{defi:#1}}
\newcommand{\remalab}[1]{\label{rema:#1}}
\newcommand{\tablab}[1]{\label{tab:#1}}

\newcommand{\figref}[1]{\ref{fig:#1}}
\newcommand{\theoref}[1]{\ref{theo:#1}}
\newcommand{\defiref}[1]{\ref{defi:#1}}
\newcommand{\remaref}[1]{\ref{rema:#1}}
\newcommand{\cororef}[1]{\ref{coro:#1}}
\newcommand{\proporef}[1]{\ref{propo:#1}}
\newcommand{\lemref}[1]{\ref{lem:#1}}
\newcommand{\eqnref}[1]{\eqref{eq:#1}}
\newcommand{\alglab}[1]{\label{alg:#1}}
\newcommand{\algref}[1]{\ref{alg:#1}}
\newcommand{\seclab}[1]{\label{sect:#1}}
\newcommand{\secref}[1]{\ref{sect:#1}}
\newcommand{\tabref}[1]{\ref{tab:#1}}

\renewcommand{\u}{u}
\newcommand{\uk}{\u^k}
\newcommand{\upk}{\u^{+k}}
\newcommand{\ukp}{\u^{k+1}}
\newcommand{\umkp}{\u^{-k+1}}
\newcommand{\umk}{\u^{-k}}
\renewcommand{\v}{v}
\newcommand{\vk}{{\v}^k}
\newcommand{\vkp}{{\v}^{k+1}}
\newcommand{\un}{\u_h}
\renewcommand{\e}{e}
\newcommand{\w}{w}
\renewcommand{\wb}{\mb{\w}}
\newcommand{\J}{J}
\newcommand{\Jb}{\mb{J}}
\newcommand{\q}{q}
\renewcommand{\r}{r}
\newcommand{\qh}{\hat{\q}}
\newcommand{\qs}{\q^*}
\newcommand{\qht}{\tilde{\qh}}
\newcommand{\qbar}{\overline{\q}}
\newcommand{\qb}{{\mb{\q}}}
\newcommand{\sig}{{\sigma}}
\newcommand{\thetah}{{\hat{\theta}}}
\newcommand{\thetas}{{\theta^*}}
\newcommand{\thetahn}{{\thetah_h}}
\newcommand{\sign}[1]{\text{ sgn}\!\LRp{#1}}
\newcommand{\sigh}{\hat{\sig}}
\newcommand{\sigs}{\sig^*}
\newcommand{\sigk}{\sig^{k}}
\newcommand{\sigpk}{\sig^{+k}}
\newcommand{\sigkp}{\sig^{k+1}}
\newcommand{\sigmkp}{\sig^{-k+1}}
\newcommand{\sigmk}{\sig^{-k}}
\newcommand{\sigb}{{\bs{\sig}}}
\newcommand{\sigbk}{\sigb^k}
\newcommand{\sigbkp}{\sigb^{k+1}}
\newcommand{\sigbn}{\sigb_h}
\newcommand{\sigbs}{\sigb^*}
\newcommand{\taub}{{\bs{\tau}}}
\newcommand{\taustar}{\tau_*}
\newcommand{\kappab}{{\bs{\kappa}}}
\newcommand{\gamb}{{\bs{\gamma}}}
\newcommand{\ut}{\tilde{\u}}
\newcommand{\sigbt}{\tilde{\sigb}}
\newcommand{\vt}{\tilde{\v}}
\newcommand{\taubt}{\tilde{\taub}}
\newcommand{\ub}{\mb{\u}}
\newcommand{\ubl}{\mb{\u}_{\lambda}}
\newcommand{\ubf}{\mb{\u}_{f}}
\newcommand{\ubp}{\mb{\u}^+}
\newcommand{\ubm}{\mb{\u}^-}
\newcommand{\vb}{\mb{\v}}
\newcommand{\um}{\u^-}
\newcommand{\up}{\u^+}
\newcommand{\ubk}{\ub^k}
\newcommand{\ubkp}{\ub^{k+1}}

\newcommand{\m}{{m}}

\newcommand{\ud}{{\dot{\u}}}
\newcommand{\vd}{{\dot{\v}}}
\newcommand{\sigbd}{{\dot{\sigb}}}
\newcommand{\taubd}{\dot{\taub}}

\renewcommand{\d}{{d}}
\newcommand{\R}{{\mathbb{R}}}
\newcommand{\kap}{\kappa}
\newcommand{\F}{\mb{F}}
\newcommand{\f}{f}
\newcommand{\g}{g}
\newcommand{\fb}{\mb{\f}}
\newcommand{\gb}{\mb{\g}}
\newcommand{\gD}{g_D}
\newcommand{\pOmega}{\partial \Omega}
\newcommand{\K}{T}
\newcommand{\Kp}{\K^+}
\newcommand{\Km}{\K^-}
\newcommand{\pK}{{\partial \K}}
\newcommand{\pKin}{{\pK^{\text{in}}}}
\newcommand{\pKout}{{\pK^{\text{out}}}}
\newcommand{\Kj}{{\K_j}}
\newcommand{\Gh}{{\mc{E}_h}}
\newcommand{\Gho}{{\mc{E}_h^o}}
\newcommand{\Ghb}{{\mc{E}_h^{\partial}}}
\newcommand{\Poly}{{\mc{P}}}
\newcommand{\D}{{D}}
\newcommand{\p}{{p}}
\newcommand{\pres}{{q}}
\newcommand{\presl}{{q}_{\lambda}}
\newcommand{\presf}{{q}_{f}}
\newcommand{\Ts}{{\mb{T}}}
\newcommand{\Vb}{{\mb{V}}}
\newcommand{\V}{{V}}
\newcommand{\Lam}{{\Lambda}}
\newcommand{\Lamb}{\bs{\Lambda}}
\newcommand{\mub}{\bs{\mu}}
\newcommand{\lambdab}{\bs{\lambda}}
\newcommand{\Thn}{\mathcal{T}_h}
\newcommand{\Th}{{\Ts_h}}
\newcommand{\Vh}{{\V_h}}
\newcommand{\Vbh}{{\Vb_h}}
\newcommand{\Lamh}{{\Lam_h}}
\newcommand{\Lambh}{{\Lamb_h}}
\newcommand{\Lamho}{{\overline{\Lam}_h}}
\newcommand{\Lambht}{{\Lamb_h^t}}
\newcommand{\Ltwo}{{L^2\LRp{\Omega}}}
\newcommand{\Ltwoe}{{L^2\LRp{\e}}}
\newcommand{\Ltw}{{L^2\LRp{\K}}}
\newcommand{\VhK}{{\V_h\LRp{\K}}}
\newcommand{\ThK}{{\Ts_h\LRp{\K}}}
\newcommand{\VbhK}{{\Vb_h\LRp{\K}}}
\newcommand{\Lamhe}{{\Lam_h\LRp{\e}}}
\newcommand{\Lambhe}{{\Lamb_h\LRp{\e}}}
\renewcommand{\L}{\mc{L}}
\renewcommand{\D}{\mc{D}}
\newcommand{\T}{\mc{T}}
\newcommand{\Tt}{\tilde{\T}}
\newcommand{\A}{\mb{A}}
\newcommand{\B}{\mb{B}}
\renewcommand{\D}{\mb{D}}
\newcommand{\Am}{\A^-}
\newcommand{\Ap}{\A^+}
\newcommand{\Aa}{\snor{\A}}
\newcommand{\Rb}{\mb{R}}
\newcommand{\C}{\mb{C}}
\renewcommand{\S}{\mb{S}}
\newcommand{\Sm}{\S^{-}}
\newcommand{\Sp}{\S^{+}}
\newcommand{\Spm}{\S^{\pm}}

\newcommand{\Lte}{{L^2\LRp{\Gh}}}
\newcommand{\n}{\mb{n}}
\newcommand{\nm}{\n^-}
\newcommand{\np}{\n^+}
\newcommand{\uh}{{\hat{\u}}}
\newcommand{\uhk}{{\uh}^k}
\newcommand{\uhkp}{{\uh}^{k+1}}
\newcommand{\us}{{\u^*}}
\newcommand{\uhn}{{\uh_h}}
\newcommand{\ubh}{{\hat{\ub}}}
\newcommand{\ubs}{{\ub^*}}
\newcommand{\ubht}{{\tilde{\ubh}}}
\newcommand{\uht}{{\tilde{\uh}}}
\newcommand{\uhtn}{{\uht_h}}
\newcommand{\uhd}{{\dot{\uh}}}
\newcommand{\ubhk}{\ubh^k}
\newcommand{\ubhkp}{\ubh^{k+1}}
\newcommand{\sigbh}{{\hat{\sigb}}}
\newcommand{\Fh}{{\hat{\F}}}
\newcommand{\Fs}{\F^*}
\renewcommand{\c}{{c}}

\newcommand{\zb}{\mb{z}}
\renewcommand{\sb}{\mb{s}}
\newcommand{\rb}{\mb{r}}
\newcommand{\betab}{\bs{\beta}}
\newcommand{\veps}{{\varepsilon}}
\newcommand{\vepsb}{\bs{\veps}}
\newcommand{\Hb}{\mb{H}}
\newcommand{\Hbt}{\Hb^t}
\newcommand\E{\mathcal{E}}
\newcommand\Eh{\mathcal{E}_h}
\newcommand\Ehi{\mathcal{E}_{h,i}}
\newcommand\Ehiint{\mathcal{E}_{h,i}^\mathrm{int}}
\newcommand\EhG{\mathcal{E}_{h}^{\Gamma}}
\newcommand\EhiG{\mathcal{E}_{h,i}^{\Gamma}}
\newcommand\EhijG{\mathcal{E}_{h,i,j}^{\Gamma}}
\newcommand\EhjiG{\mathcal{E}_{h,j,i}^{\Gamma}}
\newcommand{\Eb}{\mb{E}}
\newcommand{\Ebt}{\Eb^t}
\newcommand{\hb}{\mb{h}}
\newcommand{\eb}{\mb{e}}
\newcommand{\Hbh}{\hat{\mb{H}}}
\newcommand{\Ebh}{\hat{\mb{E}}}
\newcommand{\Hbs}{{\mb{H}^*}}
\newcommand{\Ebs}{{\mb{E}^*}}
\newcommand{\Ebht}{\Ebh^t}
\newcommand{\Hbht}{\Hbh^t}
\newcommand{\Lb}{\mb{L}}
\newcommand{\Gb}{\mb{G}}
\newcommand{\Lbh}{\hat{\Lb}}
\newcommand{\Lbs}{\Lb^*}
\newcommand{\Ib}{\mb{I}}
\newcommand{\I}{I}
\newcommand{\Q}{Q}
\newcommand{\Sigb}{\bs{\Sigma}}
\newcommand{\Piu}{\Pi_\u}
\newcommand{\Pisigb}{\Pi_\sigb}
\newcommand{\Z}{Z}
\newcommand{\Y}{Y}
\newcommand{\rhou}{\rho\u}
\newcommand{\rhov}{\rho\v}
\renewcommand{\P}{P}
\newcommand{\h}{H}

\newcommand{\mygreen}[1]{{\color[rgb]{0,.65,0} #1}}
\newcommand{\mywhite}[1]{{\color[rgb]{1.0,1.0,1.0} #1}}
\newcommand{\myred}[1]{{\color[rgb]{0.65,0.0,0.0} #1}}
\newcommand{\myblue}[1]{{\color[rgb]{0.0,0.0,1.0} #1}}
\newcommand{\mygray}[1]{{\color[rgb]{.15,.35,.60} #1}}
\newcommand{\mythemec}[1]{{\color[RGB]{51,108,121} #1}}
\newcommand{\mydarkred}[1]{{\color[rgb]{.6,.1,.1} #1}}
\newcommand{\mydarkblue}[1]{{\color[rgb]{.1,.1,.9} #1}}
\newcommand{\mygreenback}[1]{{\color[rgb]{.75,1,.75} #1}}
\newcommand{\myorange}[1]{{\color[rgb]{.76,.39,.13} #1}}
\newcommand{\mygrass}[1]{{\color[rgb]{.19,.64,.13} #1}}
\newcommand{\mysierp}[1]{{\color[RGB]{209,28,209} #1}}

\newcommand{\tanbui}[2]{\textcolor{blue}{#1} \textcolor{red}{#2}}

\newcommand{\vel}{\bs{\vartheta}}
\newcommand{\vels}{\vel^*}
\newcommand{\vele}{\vel^e}
\newcommand{\velh}{\hat{\vel}}
\newcommand{\velk}{\vel^k}
\newcommand{\velkp}{\vel^{k+1}}
\newcommand{\vhk}{\hat\v^{k}}

\newcommand{\phis}{\phi^*}
\newcommand{\phie}{\phi^e}
\newcommand{\phih}{\hat{\phi}}
\newcommand{\phihk}{\phih^k}
\newcommand{\phik}{\phi^k}
\newcommand{\phikp}{\phi^{k+1}}

\begin{frontmatter}

\title{A Multilevel Approach for Trace System in HDG Discretizations\tnoteref{t1}}
   \tnotetext[t1]{The research of   T. Bui-Thanh and S.  Muralikrishnan was partially supported by DOE grant DE-SC0018147 and NSF Grant NSF-DMS1620352. The research of J. Shadid was funded by the Department of Energy Office of Science, Advanced Scientific Computing Research (ASCR) Applied Math Program. We are grateful for the support. An earlier form of this work was first presented at the FEM Rodeo, February 2015, Southern Methodist University. The views expressed in the article do not necessarily represent the views of the U.S. Department of Energy or the United States Government. Sandia National
    Laboratories is a multimission laboratory managed and operated by National Technology and Engineering Solutions of Sandia, LLC., a wholly owned subsidiary of Honeywell International, Inc., for the U.S. Department of Energy's National Nuclear Security Administration under contract DE-NA-0003525.}

    \author[srk]{Sriramkrishnan Muralikrishnan}
    \author[srk,Tan]{Tan Bui-Thanh}
    \author[John1,John2]{John N. Shadid}
    \address[srk]{Department of Aerospace Engineering and Engineering Mechanics,\\ The University of Texas at Austin, TX 78712, USA.} 
    \address[Tan]{The Institute for Computational Engineering
    \& Sciences,\\ The University of Texas at Austin, Austin, TX 78712,
    USA.}
    \address[John1]{Computational Mathematics Department, Sandia National Laboratories, P.O. Box 5800, MS 1321,
    Albuquerque, NM 87185.}
    \address[John2]{Department of Mathematics and Statistics, University of New Mexico, Albuquerque, NM 87131.}
\begin{abstract}
    We propose a multilevel approach for trace systems resulting from
    hybridized discontinuous Galerkin (HDG) methods. The key is to
    blend ideas from nested dissection, domain decomposition, and
    high-order characteristic of HDG discretizations. Specifically, we first
    create a coarse solver by eliminating and/or limiting the front
    growth in nested dissection. This is accomplished by projecting
    the trace data into a sequence of same or high-order polynomials
    on a set of increasingly $h-$coarser edges/faces. We then combine the
    coarse solver with a block-Jacobi fine scale solver to form a 
    two-level solver/preconditioner. Numerical experiments indicate that the performance of the
    resulting two-level solver/preconditioner depends only on the
    smoothness of the solution and  is independent of the nature of
    the PDE under consideration. While the proposed algorithms are developed within
    the HDG framework, they are applicable to other hybrid(ized) high-order finite
    element methods. Moreover, we  show that our multilevel algorithms can be
    interpreted as a multigrid method with specific intergrid transfer
    and smoothing operators. With several numerical
    examples from Poisson, pure transport, and convection-diffusion
    equations we demonstrate the robustness and scalability of the
    algorithms. 
\end{abstract}


\end{frontmatter}


\section{Introduction}
Hybridized discontinuous Galerkin (HDG) methods introduced
a decade ago
\cite{CockburnGopalakrishnanLazarov:2009:UHO} have now been developed
for a wide range of PDEs including, but not limited to, Poisson-type
equation \cite{CockburnGopalakrishnanLazarov:2009:UHO,
  CockburnGopalakrishnanSayas10}, Stokes equation
\cite{CockburnGopalakrishnan09, NguyenPeraireCockburn10,
  rhebergen2017analysis}, Euler and Navier-Stokes equations
\cite{NguyenPeraireCockburn11, MoroNguyenPeraire11,
  rhebergen2018hybridizable}, wave equations
\cite{NguyenPeraireCockburn11b, NguyenPeraireCockburn11a,
  GriesmaierMonk11, lee2017analysis}. In \cite{Bui-Thanh15, Bui-Thanh15a,
  bui2016construction}, an upwind HDG framework was proposed that
provides a unified and a systematic construction of HDG methods for a
large class of PDEs.

Roughly speaking, HDG methods combine the advantages of hybrid(ized)
methods and discontinuous Galerkin (DG) discretizations. In particular,
its inherent characteristics from DG include: i) arbitrary high-order
with compact stencil; ii) ability to handle complex geometries; iii)
local conservation; and iv) upwinding for hyperbolic systems. On the
other hand, it also possesses the advantages of hybrid(ized) methods,
namely, i) having smaller and sparser linear system for steady state
problems or time-dependent problems with implicit time integrators;
ii) {$h\p$-adaptivity-ready using  the trace space;}
iii) facilitating
multinumerics with different hybrid(ized) methods in different parts
of the domain; and iv) when applicable, providing superconvergence by
local post-processing
\cite{AB85,CockburnGopalakrishnanLazarov:2009:UHO}. Thus for complex
multiphysics applications with disparate spatial and temporal scales (e.g. magnetohydrodynamics and atmospheric flows), high-order HDG spatial
discretization together with high-order implicit time integrator is a
strong candidate for large scale simulations in modern extreme-scale computing
architectures owing to its high computation-to-communication ratio.

The main challenge facing hybridized methods is, however, the
construction of scalable solvers/preconditioners for the resulting trace systems.
Over the past 30 years, a
tremendous amount of research has been devoted to the convergence of
multigrid methods for such linear systems, both as iterative methods
and as preconditioners for Krylov subspace methods.  Optimal
convergence with respect to the number of unknowns is usually obtained
under mild elliptic regularity assumptions
\cite{bramble1993multigrid,Bramble94uniformconvergence,bpx1991}.
Multigrid algorithms have been developed for mortar domain
decomposition methods \cite{339353,YotovMultigrid}.  Several multigrid
algorithms have been proposed for hybridized mixed finite element
methods \cite{Bren_MG_MFE_92, Chen_equiv_96}, whose optimal convergence has
already been established \cite{Braess:1990,
  Brenner:nonconfmg}. Multigrid algorithms based on restricting the
trace (skeletal) space to linear continuous finite element space has been
proposed for hybridized mixed methods
\cite{Gopalakrishnan09aconvergent}, hybridized discontinuous Galerkin
methods \cite{cockburn2014multigrid} and weak Galerkin methods
\cite{chen2015auxiliary}. 


Iterative solvers/preconditioners for solving HDG trace systems are,
however, scarced.  Recent efforts on multigrid methods have been
presented for elliptic PDEs
\cite{cockburn2014multigrid,chen2014robust,kronbichler2018performance,wildey2018unified}. Attempts
using domain decomposition type solvers/preconditioners have been
proposed for elliptic PDEs
\cite{GanderHajian:2015:ASM,gander2018analysis}, Maxwell's equations
\cite{li2014hybridizable,he2016optimized}, and hyperbolic systems
\cite{iHDG,iHDGII,diosady2011domain}. Recently, an approximate block factorization preconditioner
for Stokes equations have been developed in
\cite{rhebergen2018preconditioning}.  Thus, there is a critical need
for developing robust, scalable solvers/preconditioners for high-order
HDG methods to tackle high-fidelity large-scale
simulations of multiphysics/multiscale systems. As a step towards to
achieve this goal, we propose a multilevel approach for both solving
and preconditioning the trace system of HDG discretizations. As will
be demonstrated, unlike existing approaches our proposed algorithms are reasonably robust and
scalable beyond elliptic PDEs.

Now let us briefly discuss the main idea behind our approach.  The
goal is to {\em advance the nested dissection}
\cite{george1973nested}\textemdash a fill-in reducing direct solver
strategy\textemdash to create a {\em scalable and robust solver
utilizing the high-order and variational structure of HDG
methods}. This is achieved by  projecting the skeletal data
at different levels to either same or high-order polynomial on a set
of increasingly $h-$coarser edges/faces. Exploiting the concept of two-level domain decomposition
methods we make use of our approach as a coarse solver together with a fine scale solver
(e.g. block-Jacobi) to create a solver/preconditioner for solving the
trace system iteratively. Thanks to its root in direct solver
strategy, {\em the behavior of our approach seems to depend only on the
solution smoothness, but otherwise is independent of the nature
of the underlying PDE}. Indeed,
the numerical experiments
show that the algorithms are robust even for transport equation with
discontinuous solution and elliptic equations with highly
heterogeneous and discontinuous permeability. For convection-diffusion
equations our multilevel preconditioning algorithms are scalable and reasonably robust
for not only  diffusion-dominated but also moderately
convection-dominated regimes.
We show that {\em the two-level approach can
also be interpreted as a multigrid algorithm with specific intergrid
transfer and smoothing operators}. Our complexity
estimates show that the cost of the multilevel algorithms is somewhat
in between the cost of nested dissection and standard multigrid
solvers.

This paper is organized as follows. Section \secref{model_problem}
introduces the model problem, notations, and an upwind HDG method
considered in this paper. In Section \secref{multilevel}, we first
recall the nested dissection approach and then explain how it can be advanced
 using HDG variational structure and the two-level domain
decomposition approach. We also show that our approach can be
interpreted as a multigrid method, and estimate the complexity of the
proposed multilevel solvers. Section \secref{numerical} presents
several numerical examples to study the robustness and scalability of
the proposed algorithm for Poisson, transport and convection-diffusion
equations as the mesh and the solution order are refined. Finally,
Section \secref{conclusion} summarizes our findings and discusses
future research directions.

\section{Model Problem, notations, and an upwind HDG method}
\seclab{model_problem}
We consider the following model problem 

\begin{subequations}
  \label{model_problem}
\begin{align}
  -\Div\LRp{{\bf K} \Grad \u} + \Div\LRp{\betab \u} &= f, \quad &\text{in} \quad \Omega, \\
    \u &=g_{D},  \quad  &\text{on} \quad \pOmega.
\end{align}
\end{subequations}
where $\Omega$ is an open, bounded, and connected subset of $\mathbb{R}^d$, with $d\in \LRc{2,3}$\footnote{Note that the treatment for 1D problems is trivial and hence omitted.}. Here, $\tK$ is a
symmetric, bounded, and uniformly positive definite tensor. 
Let $\Thn$ be a conforming partition of $\Omega$ into
$\Nel$ non-overlapping elements $\Kj$, $j = 1, \hdots, \Nel$, with
Lipschitz boundaries such that $\Thn := \cup_{j=1}^\Nel \Kj$ and
$\overline{\Omega} = \overline{\Thn}$. The mesh size $h$ is defined as $h
:= \max_{j\in \LRc{1,\hdots,\Nel}}\text{diam}\LRp{\Kj}$. We denote the
skeleton of the mesh by $\Gh := \cup_{j=1}^\Nel \pK_j$:
the set of all (uniquely defined) interfaces $\e$ between elements. We conventionally identify $\nm$ as the outward 
normal vector on the boundary $\pK$ of element $\K$ (also denoted as $\Km$) and $\np = -\nm$ as the outward normal vector of the boundary of a neighboring element (also denoted as $\Kp$). Furthermore, we use $\n$ to denote either $\nm$ or $\np$ in an expression that is valid for both cases, and this convention is also used for other quantities (restricted) on
a face $\e \in \Gh$.

For simplicity, we define $\LRp{\cdot,\cdot}_\K$ as the
$L^2$-inner product on a domain $\K \subset \R^\d$ and
$\LRa{\cdot,\cdot}_\K$ as the $L^2$-inner product on a domain $\K$ if
$\K \subset \R^{\d-1}$. We shall use $\nor{\cdot}_{\K} :=
\nor{\cdot}_{\Ltw}$ as the induced norm for both cases.
Boldface lowercase letters are conventionally used for vector-valued functions and in that
case the inner product is defined as $\LRp{\ub,\vb}_\K :=
\sum_{i=1}^m\LRp{\ub_i,\vb_i}_\K$, and similarly $\LRa{\ub,\vb}_\K :=
\sum_{i = 1}^m\LRa{\ub_i,\vb_i}_\K$, where $\m$ is the number of
components ($\ub_i, i=1,\hdots,\m$) of $\ub$.  Moreover, we define
$\LRp{\ub,\vb}_{\Thn} := \sum_{\K\in \Thn}\LRp{\ub,\vb}_\K$ and
$\LRa{\ub,\vb}_\Gh := \sum_{\e\in \Gh}\LRa{\ub,\vb}_\e$ whose
induced norms are clear, and hence their definitions are
omitted. We  employ boldface uppercase letters, e.g. $\tK$, to
denote matrices and tensors. 
We denote by $\mc{Q}^\p\LRp{\K}$ the space of tensor product polynomials of degree at
most $\p$ in each dimension on a domain $\K$. We use the terms ``skeletal unknowns" and ``trace unknowns" interchangeably and they both refer to the unknowns on the mesh skeleton.

 First, we cast equation \eqref{model_problem} into the following first-order form:
 \begin{subequations}
   \label{model_first_order}
     \begin{align}\label{mixed_eq1}
         \sigb&=-{\bf K}\Grad \u \quad &\text{ in }  \Omega,\\\label{mixed_eq2}
         \Div \sigb + \Div\LRp{\betab \u} &= f \quad &\text{ in } \Omega,\\\label{mixed_boundary}
    \u &= g_{D} \quad &\text {on } \pOmega.
\end{align}
\end{subequations}

The upwind hybridized DG method \cite{Bui-Thanh15,Bui-Thanh15a} for the
discretization of equation \eqref{model_first_order} reads: seek $\LRp{\u, \sigb,\lambda}$ such that
\begin{subequations}
  \label{HDG}
\begin{align}\label{HDG_local1}
    \LRp{{\bf K}^{-1}\sigb,\vb}_\K -\LRp{\u,\Div \vb}_\K + \LRa{\lambda,\vb \cdot \n}_\pK &= 0, \\\label{HDG_local2}
    -\LRp{\sigb,\Grad \w}_\K -\LRp{\betab\u,\Grad \w}_\K + \LRa{\LRp{\widehat{\sigb} + \widehat{\betab\u}} \cdot \n,\w}_\pK &= \LRp{f,w}_\K, \\
    \LRa{\jump{\LRp{\widehat{\sigb}+\widehat{\betab\u}} \cdot \n},\mu}_\e&=0,
    \label{conservation}
\end{align}
\end{subequations}
where the upwind HDG numerical flux \cite{Bui-Thanh15,Bui-Thanh15a} $\LRp{\widehat{\sigb}+\widehat{\betab\u}}$ is given by
\begin{equation}
    \LRp{\widehat{\sigb}+\widehat{\betab\u}} \cdot \n := \sigb \cdot \n + \betab\cdot\n\u + \frac{1}{2}\LRp{\sqrt{|\betab\cdot\n|^2+4}-\betab\cdot\n}(\u - \lambda).
    \label{uHDG_flux}
\end{equation}

For simplicity, we have suppressed the explicit statement that equations
\eqref{HDG_local1}, \eqref{HDG_local2} and \eqref{conservation} must
hold for all test functions $\vb \in \VbhK$, $\w \in W_{h}\LRp{\K}$,
and $\mu \in M_{h}\LRp{\e}$, respectively (this is implicitly
understood throughout the paper), where $\Vbh$, $W_h$ and $M_h$ are
defined as
\begin{subequations}
\begin{align}
\Vbh\LRp{\Thn} &= \LRc{\vb \in \LRs{L^2\LRp{\Thn}}^d:
    \eval{\vb}_{\K} \in \LRs{\mc{Q}^\p\LRp{\K}}^d, \forall \K \in \Thn}, \\
 W_h\LRp{\Thn} &= \LRc{\w \in L^2\LRp{\Thn}:
   \eval{\w}_{\K} \in \mc{Q}^\p\LRp{\K}, \forall \K \in \Thn}, \\
 \eqnlab{traceSpace}
M_h\LRp{\Gh} &= \LRc{\lambda \in \Lte:
    \eval{\lambda}_{\e} \in \mc{Q}^\p\LRp{\e}, \forall \e \in \Gh},
\end{align}
\end{subequations}
and similar spaces  $\VbhK$, $W_h\LRp{\K}$ and $M_h\LRp{\e}$ on $\K$ and $\e$ can be defined by replacing $\Thn$ with
$\K$ and $\Gh$ with $\e$, respectively.

\section{A Multilevel solver for the HDG trace system}
\seclab{multilevel}
The HDG solution process involves the following steps: 
\begin{enumerate}
    \item Express the local volume unknowns $\u$ and $\sigb$,
      element-by-element, as a function of the skeletal unknowns
      $\lambda$ using \eqref{HDG_local1}, \eqref{HDG_local2}. The
      well-posedness of this step can be found in \cite{Bui-Thanh15,Bui-Thanh15a} (and the references therein).
    \item Use the conservation condition \eqref{conservation} to
construct a global linear system involving only the skeletal
unknowns and solve it. Similarly, this step can be rigorously justified as in \cite{Bui-Thanh15,Bui-Thanh15a} (and the references therein).
    \item Recover the local volume
unknowns in an element-by-element fashion completely independent of each other using \eqref{HDG_local1}, \eqref{HDG_local2}. 
\end{enumerate}

        The main advantage of this Schur complement approach is that,
        for high-order, the global trace system is much smaller and
        sparser compared to the linear system for the volume unknowns
        \cite{CockburnGopalakrishnanLazarov:2009:UHO,bui2016construction}.
        Since Steps 1. and 3. are embarassingly parallel, the main
        bottle neck for HDG in large scale simulations is the solution
        of the global trace system (Step 2.).

\subsection{A brief review on solvers/preconditioners for HDG trace system}
In this section we briefly discuss existing works on solvers for HDG
methods and our contributions. The first geometric
multigrid solver for HDG methods was introduced in
\cite{cockburn2014multigrid}. The main idea was to transfer the residual
from the skeletal space to linear continuous Galerkin FEM space, and
then carry out the standard multigrid algorithm. A
similar concept with few modifications was pursued in
\cite{chen2014robust} for the simulation of high frequency Helmholtz
equation discretized by HDG. In \cite{kronbichler2018performance} the
authors studied a version of the multigrid algorithm proposed in
\cite{cockburn2014multigrid} along with multigrid algorithms for
continuous Galerkin and interior penalty discontinuous Galerkin for
standard elliptic equation. They concluded that both continuous and
interior penalty discontinuous Galerkin algorithms with multigrid
outperforms HDG with multigrid in terms of time to solution by a
significant margin. One level Schwarz type domain decomposition
algorithms in the context of HDG have been studied for elliptic equation
\cite{GanderHajian:2015:ASM,gander2018analysis}, hyperbolic systems
\cite{iHDG,iHDGII} and Maxwell's equations
\cite{li2014hybridizable,he2016optimized}. A balancing domain
decomposition by constraints algorithm for HDG was introduced in
\cite{diosady2011domain} and studied for Euler and Navier-Stokes
equations. A unified
geometric multigrid algorithm based on Dirichlet-to-Neumann maps was
developed in \cite{wildey2018unified} for hybridized methods
including HDG. An approximate block factorization based preconditioner for 
HDG discretization of Stokes system was presented in \cite{rhebergen2018preconditioning}.

The objective of our current work is to develop a robust multilevel
solver and preconditioner for HDG discretizations for a wide variety of PDEs. 
The ultimate goal is to significantly reduce factorization and memory costs 
compared to a direct solver. Unlike cost reduction strategies for direct solvers
in \cite{ martinsson2009fast,gillman2014direct,ho2016hierarchical}
which utilizes the elliptic nature of PDEs, here we exploit the
high-order and variational structure of HDG methods. As a result, our
method is applicable to not only elliptic but also parabolic, hyperbolic, and
mixed-type PDEs. For ease of the exposition and implementation, we will
focus only on structured grids. While extending the algorithm to
block-structured or nested grids is fairly straightforward,
applying it to a completely unstructured grid is a non-trivial task
and hence left for future work.

\subsection{Nested dissection}
As nested dissection idea is utilized in the proposed multilevel
algorithm, we briefly review its concept (more details can be found in
\cite{george1973nested,george1978automatic,lipton1979generalized,liu1992multifrontal,duff1983multifrontal,davis2006direct}). Nested
dissection (ND) is a fill-in reducing ordering strategy introduced in
1973 \cite{george1973nested} for efficient solution of linear
systems. 
Consider a $\p=2$ solution on an $8\times8$  quadrilateral HDG skeletal mesh in
Figure \figref{level1_nd} (the boundary nodes are eliminated for
clarity). In the ND algorithm, we identify a set of separators which
divide the mesh into independent subdomains. For example, the black
nodes in Figure \figref{level1_nd} divide the mesh into four
independent subdomains each of which can be recursively divided into
four subdomains and so on. We then order the nodes such that the red
ones are ordered first, followed by blue and then the black ones. This
will enable a recursive Schur complement approach in a multilevel
fashion and the nodes remaining in the system after elimination at
each level is shown in Figure \figref{ND_levels} (for three levels).
    
\begin{figure}[h!b!t!]
    \subfigure[Level 1]{
    \includegraphics[trim=8cm 7cm 10cm 7cm,clip=true,width=0.48\textwidth]{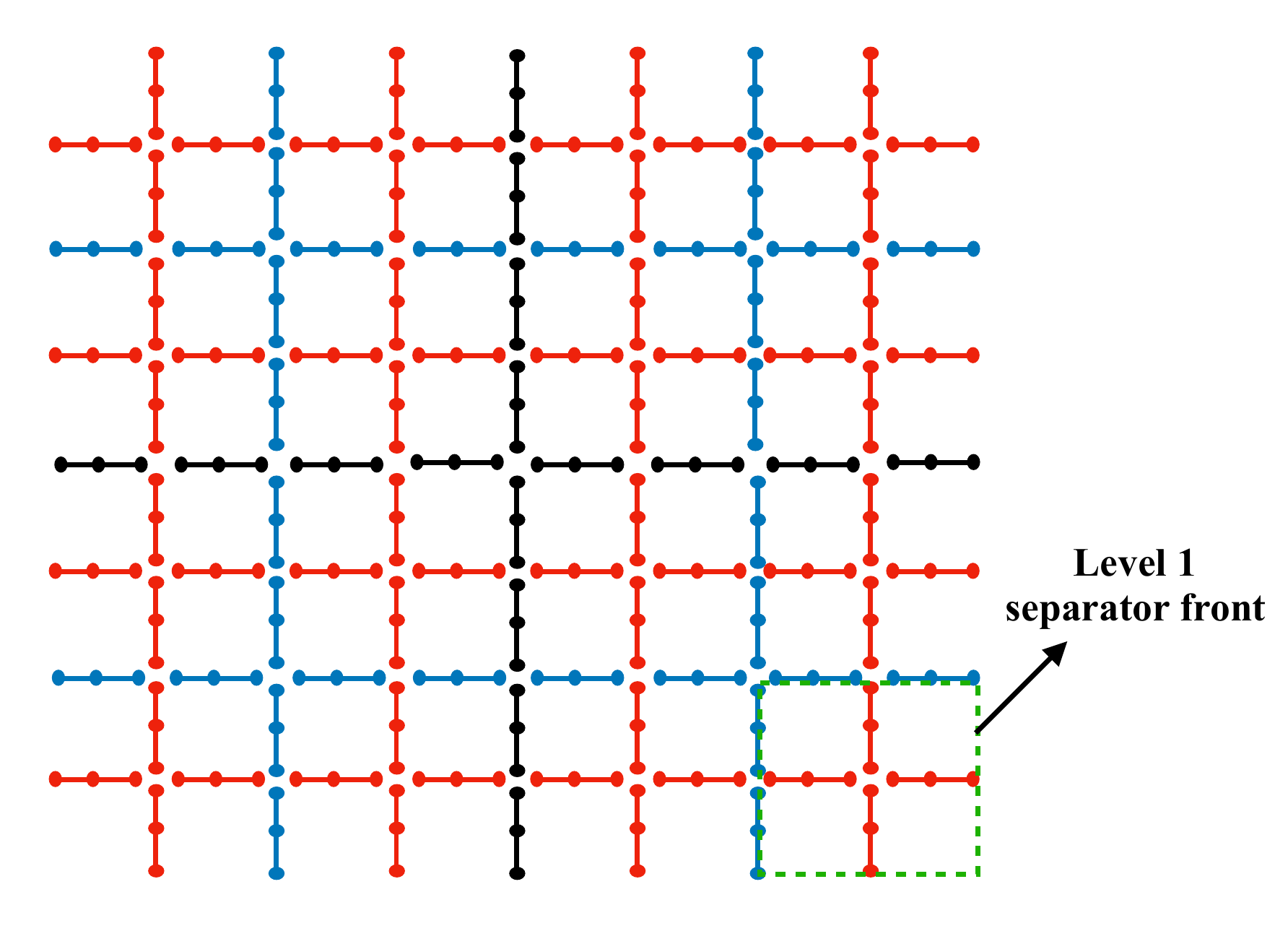}
    \figlab{level1_nd}
  }
    \subfigure[Level 2]{
    \includegraphics[trim=8cm 7cm 10cm 7cm,clip=true,width=0.48\textwidth]{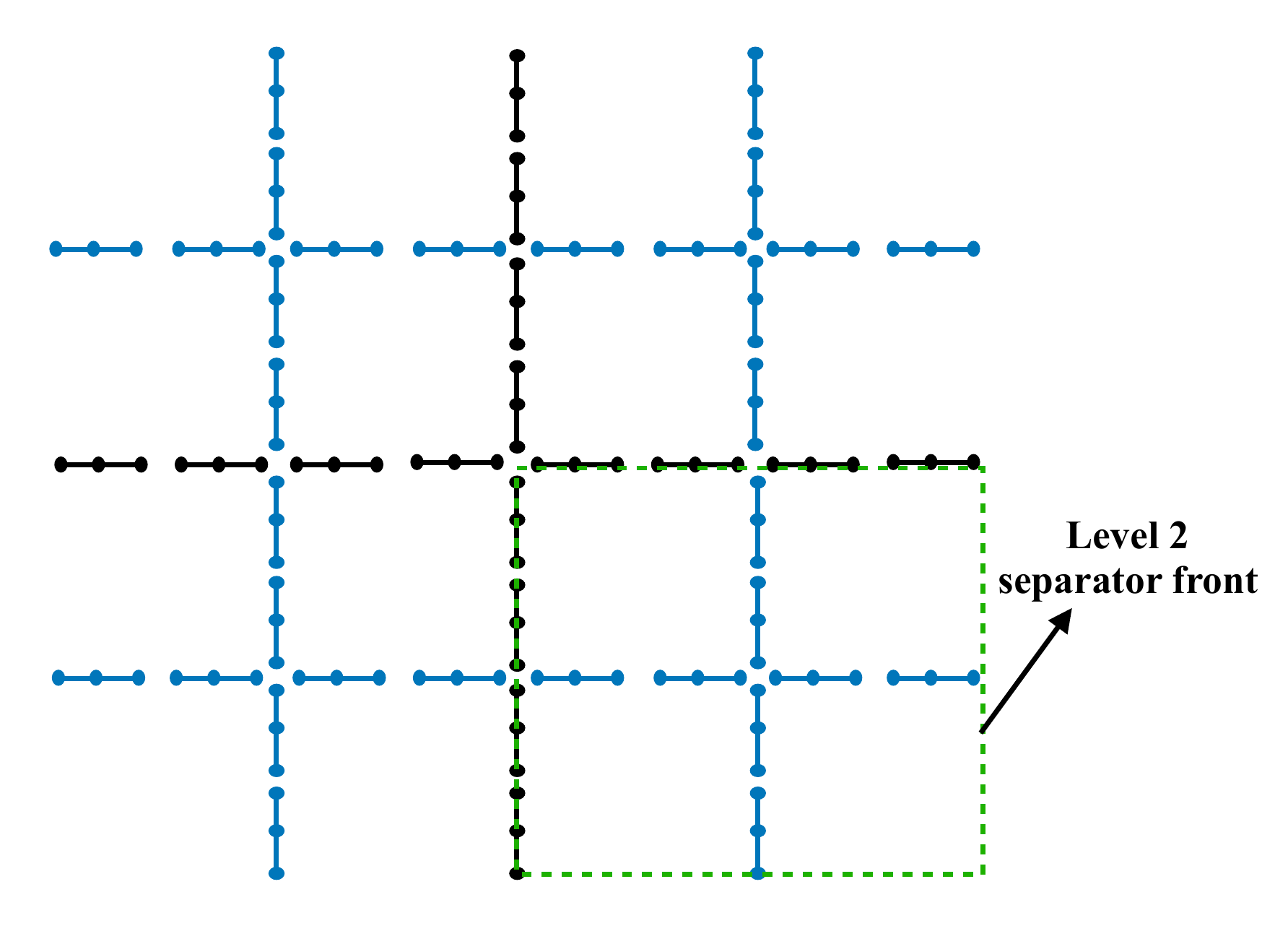}
    \figlab{level2_nd}
  }
    \begin{center}
    \subfigure[Level 3]{
    \includegraphics[trim=8cm 7cm 10cm 7cm,clip=true,width=0.48\textwidth]{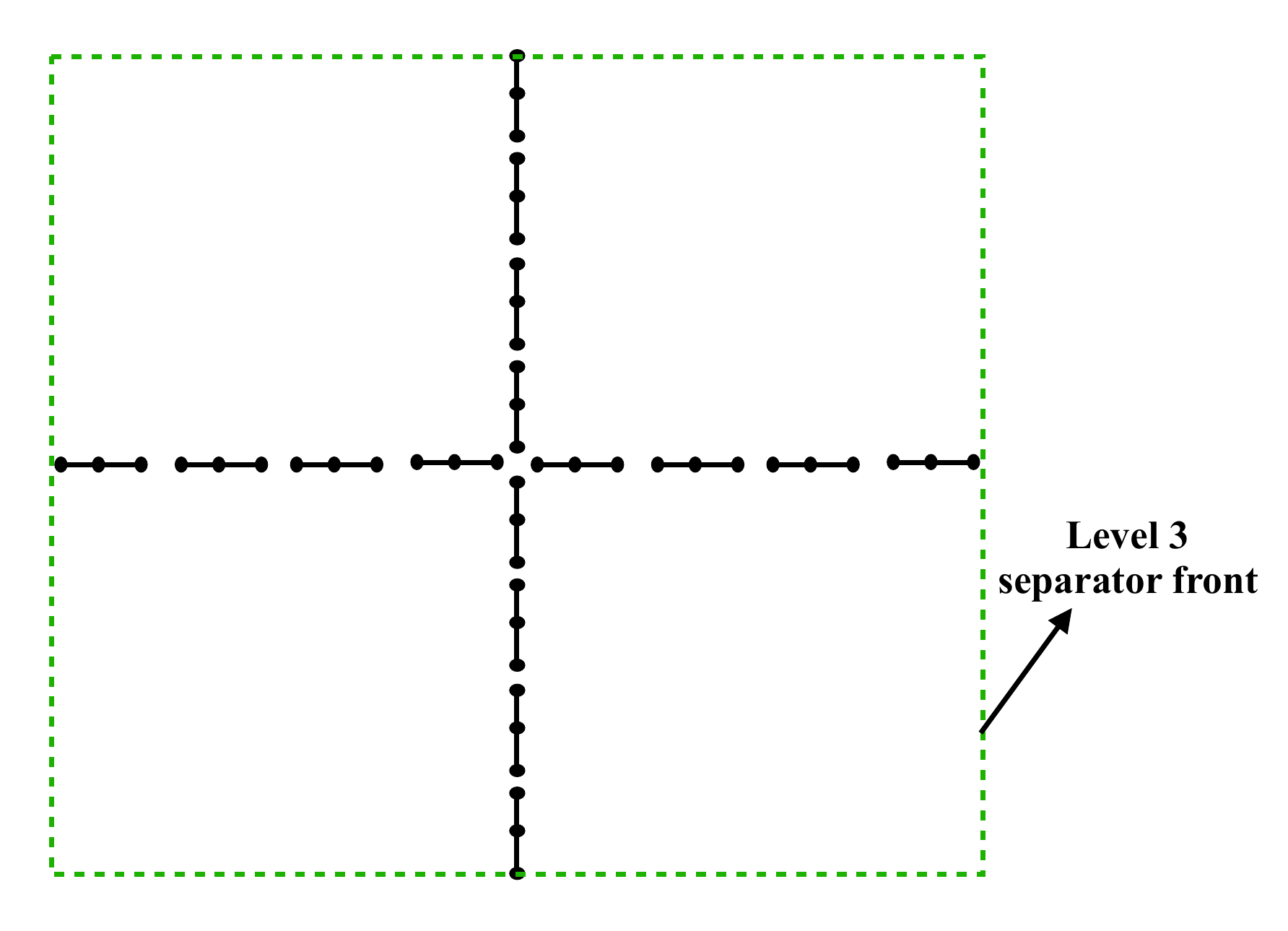}
    \figlab{level3_nd}
  }
    \end{center}
    \caption{An example of three levels in the nested dissection (ND) algorithm. The red crosses correspond to level 1 separator fronts and there
    are 16 fronts in Figure \figref{level1_nd}, each having 4 edges. The blue crosses correspond to level 2 separators and in Figure \figref{level2_nd} there
    are four level 2 fronts, each having 8 edges. The black cross correspond to level 3 separator and in Figure \figref{level3_nd} there is one level 3 front
    with 16 edges. The circles on each edge represent the nodes and there are three nodes in each edge corresponding to a solution order of $\p=2$.}
  \figlab{ND_levels}
\end{figure}
    
    There are several advantages to this algorithm. First, it can be
    shown that the serial complexity of factorization for an $N \times
    N$ matrix arising from 2D problems with this procedure is
    $\mc{O}(N^{3/2})$, and the memory requirement is $\mc{O}(NlogN)$
    \cite{george1973nested}. Whereas with a naive lexicographic
    ordering, it is $\mc{O}(N^2)$ for factorization and
    $\mc{O}(N^{3/2})$ for memory \cite{george1973nested}. Moreover, in
    2D it is optimal in the sense that the lower bound of the cost for factorization using any ordering algorithm
    is  $\mc{O}(N^{3/2})$
    \cite{george1973nested}.  Second, all the leaf calculations at
    any level are independent of each other and hence are amenable
    to parallel implementation.  However, in 3D the cost is
    $\mc{O}(N^2)$ for factorization and $\mc{O}(N^{4/3})$ for
    memory \cite{george1973nested}, and we no longer have the
    tremendous savings as in 2D \cite{eijkhout2014introduction}. This
    can be intuitively understood using Figure \figref{ND_levels}. The separator
    fronts (the crosses at each level in Figure \figref{ND_levels})
    grow in size as the level increases. For example, the black crosses have {more
      edges}
    and  nodes than the blue ones, which in turn has more edges
    and nodes than the red ones. 
    On the  last level, the  size is $\mc{O}(N^{1/2})$ for 2D and the cost of a
    dense matrix factorization for the separator front matrix corresponding to the last level is
    $\mc{O}(N^{3/2})$. In 3D the size of the separator at last level is 
    $\mc{O}(N^{2/3})$ and hence the factorization cost becomes $\mc{O}(N^2)$.  Thus
    in order to cut down the factorization and storage cost of the ND algorithm
    we need to reduce the front growth.

    There have been many efforts in this direction over the past
    decade
    \cite{martinsson2009fast,xia2009superfast,schmitz2012fast,gillman2014direct,ho2016hierarchical}. The
    basic idea in these approaches is to exploit the low rank
    structure of the off-diagonal blocks, a characteristic of
    elliptic PDEs, to compress the fronts. In this way one can obtain a
    solver which is $\mc{O}(N)$ or $\mc{O}(NlogN)$ in both 2D and 3D
    \cite{gillman2014direct,ho2016hierarchical}. Unfortunately, since
    the compression capability is a direct consequence of the ellipticity, it is
      not trivially applicable for convection-dominated or pure hyperbolic PDEs.
      {\em Our goal here is to construct a multilevel
      algorithm that is independent of the nature of PDE and at the
      same time more efficient than ND}. At the heart of our approach is the exploitation of the
    high-order properties of HDG and the underlying variational
    structure.

    \subsection{Direct multilevel solvers}
    \seclab{directSolvers}
    In our multilevel algorithm, we start with the ND ordering of the
    original fine mesh (red, blue, and black edges) as in Figure
    \figref{level1_nd}. Here, by edges we mean the original elemental
    edges (faces) on the fine mesh. Let us denote the fine mesh as
    level 0.
    In Figure \figref{level_0_ML}, 
    all red crosses have 4 edges, blue crosses have 8 edges and black cross has 16 edges. On
    these edges are the trace spaces \eqnref{traceSpace}, and thus
    going from level $k$ to level $\LRp{k+1}$ the separator front grows by a
    factor of two. We propose to reduce the front growth by lumping the
    edges so that each cross at any level has only four (longer) edges
    as on level 1 separator fronts. We accomplish this goal by projecting the traces on
    original fine mesh skeletal edges into a single trace
    space on a single longer edge (obtained by lumping the edges). Below are the details on how we lump edges and how we construct the projection operators.

    The lumping procedure is straightforward. For example, longer edges on each blue cross in Figure \figref{level_1_create} are obtained 
    by lumping the corresponding two blue (shorter) edges. Similarly, longer edges on each black cross in Figure \figref{level_1_create} are obtained 
    by lumping the corresponding four black (shorter) edges. The resulting skeleton mesh with the same number of edges on the separator fronts in
      all levels forms level 1 in our multilevel algorithm.

    
    
    Next, we project the traces spaces on shorter edges into a single
    trace space on the corresponding lumped edge.  The three obvious
    choices for the solution order of the single trace spaces: (1)
    lower than, (2) same as, or (3) higher than the solution order on
    the shorter edges.  Low-order option is not sensible as we have
    already coarsened in $h$. In particular, additional coarsening in
    $\p$, i.e. option (1), makes the solver even further away from
    being ``direct". Moreover, since  we already invert matrices of size
    $\mc{O}((\p+1)^2)$ for separators in level 1,  low-order
    choice will not help in reducing the cost. 
For option (2), we obtain separators which are
    increasingly coarsened in $h$, and clearly when we apply the ND
    algorithm this approach do not yield a direct solution nor
    $h$-convergence. However, it can be used in the construction of an iterative solver/preconditioner as we shall discuss
     in section \secref{twoLevel}.

    
    Option (3), i.e. high-order projection, is more interesting as we
    now explain. Due to exponential convergence in $p$ for smooth
      solution \cite{CockburnGopalakrishnanSayas10}, we
    can compensate for the coarseness in $h$ by means of refinement in
    $\p$. As a result, for sufficiently smooth solution, projecting to
    high-order trace spaces can provide accurate approximations to
    the original ND solution while almost avoiding the front growth.
    In our algorithm we enrich the order in the following fashion: if
    $\p$ is the solution order on the original fine mesh, then level 1
    separators have solution of order $\p$,
$\p+1$ for separators on
    level 2, $\p+2$ for separators on level 3 and so on. For practical
    purposes we also (arbitrarily) limit the growth to order $10$ to
    actually stop the front growth after 10 orders.  Specifically, for
    a generic level $k$ we take the solution order on the separator
    fronts as $\p_k=\min\LRc{\p+(k-1),10}$. We would like to point out that
    this enriching strategy is natural, but by no means
    optimal. Optimality requires balancing accuracy and computational
    complexity, which in turns requires rigorous error analysis of the
    enrichment. This is, however, beyond the scope of this paper and thus left  for future research.

    To the end of the paper, we denote option (2) as multilevel (ML)
    and option (3) as enriched multilevel (EML) to
    differentiate between them. In Figures \figref{ML_levels} and
    \figref{EML_levels} are different levels corresponding to the ML
    and EML approaches for solution order of $\p=2$ on the original
    fine mesh. Level 0 of both cases corresponds to Figure
    \figref{level1_nd}. Note that the number of circles on each edge
    is equal to the solution order plus one. For example, the  solution
    order on each edge of Figure \figref{ML_levels3} is $2$, while the enriched solution order
    is $4$ for each edge in Figure \figref{EML_levels3}.

\begin{figure}[h!b!t!]
    \subfigure[Level 0]{
    \includegraphics[trim=8cm 7cm 10cm 7cm,clip=true,width=0.48\textwidth]{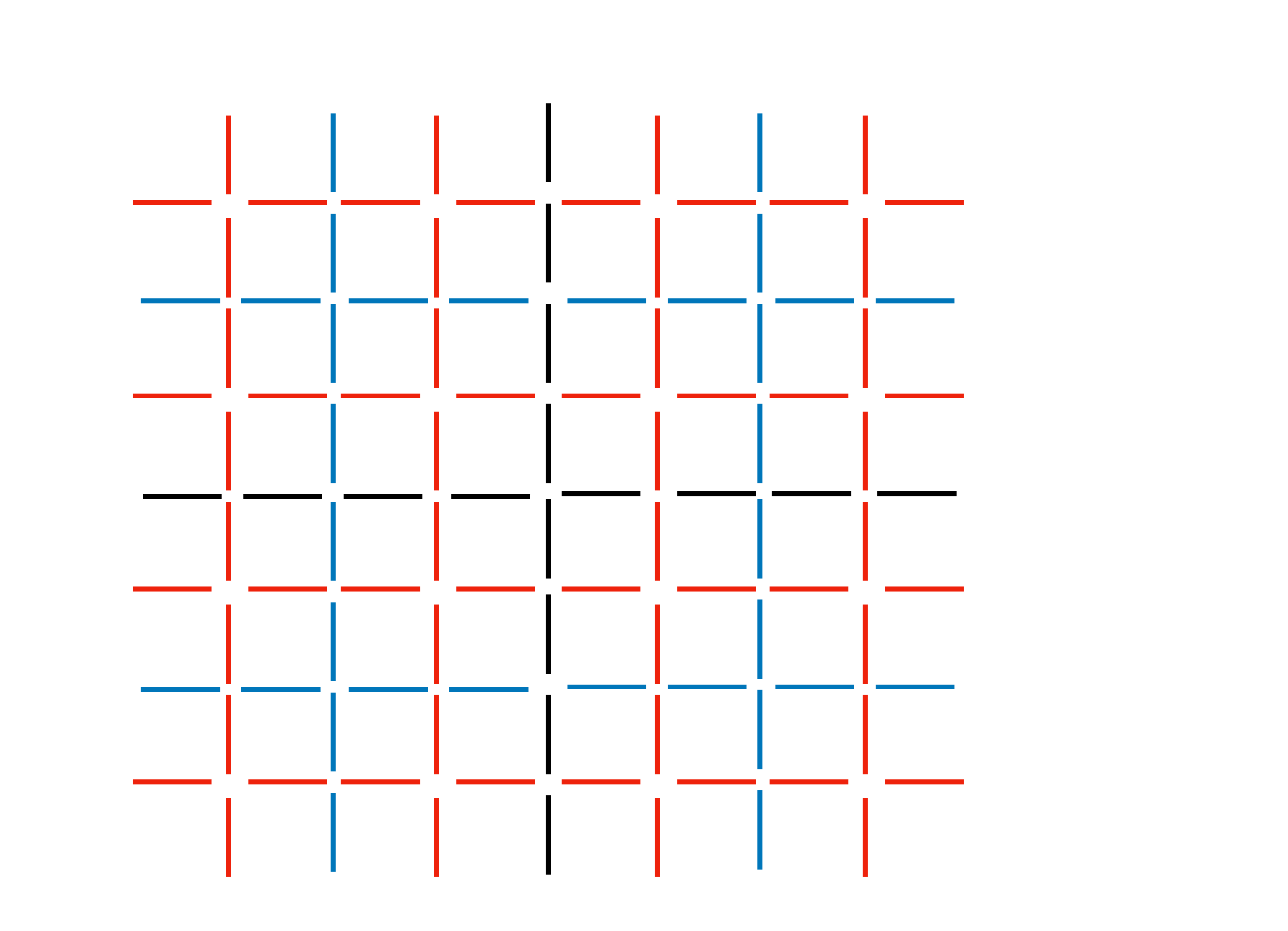}
    \figlab{level_0_ML}
  }
    \subfigure[Level 1]{
    \includegraphics[trim=8cm 7cm 10cm 7cm,clip=true,width=0.48\textwidth]{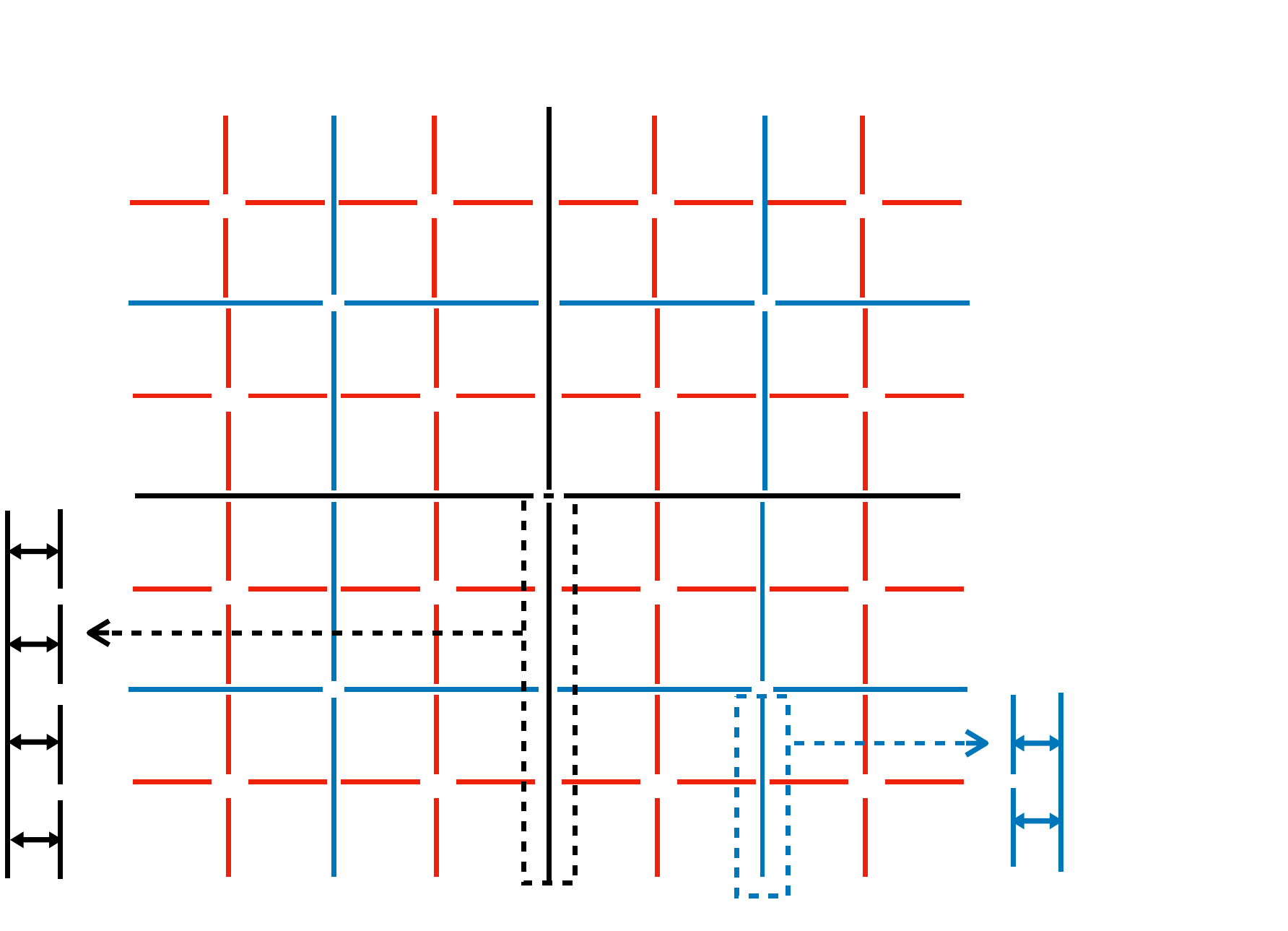}
    \figlab{level_1_create}
  }
    \caption{Creation of level 1 from level 0 in the multilevel
      algorithm: every two short blue edges in Figure \figref{level_0_ML}
      are projected on to the corresponding single long blue edge in
      Figure \figref{level_1_create}. Similarly, every four short black
      edges in Figure \figref{level_0_ML} are projected on to the corresponding single long black edge in Figure \figref{level_1_create}. In level 1,
      all the separator fronts have the same number of edges (of
      different lengths), which is 4. The nodes on each edge (circles
      in Figure \figref{ND_levels}) are not shown in this figure.}
    \figlab{Projection}
\end{figure}

\begin{figure}[h!b!t!]
    \subfigure[Level 1]{
    \includegraphics[trim=8cm 7cm 10cm 7cm,clip=true,width=0.48\textwidth]{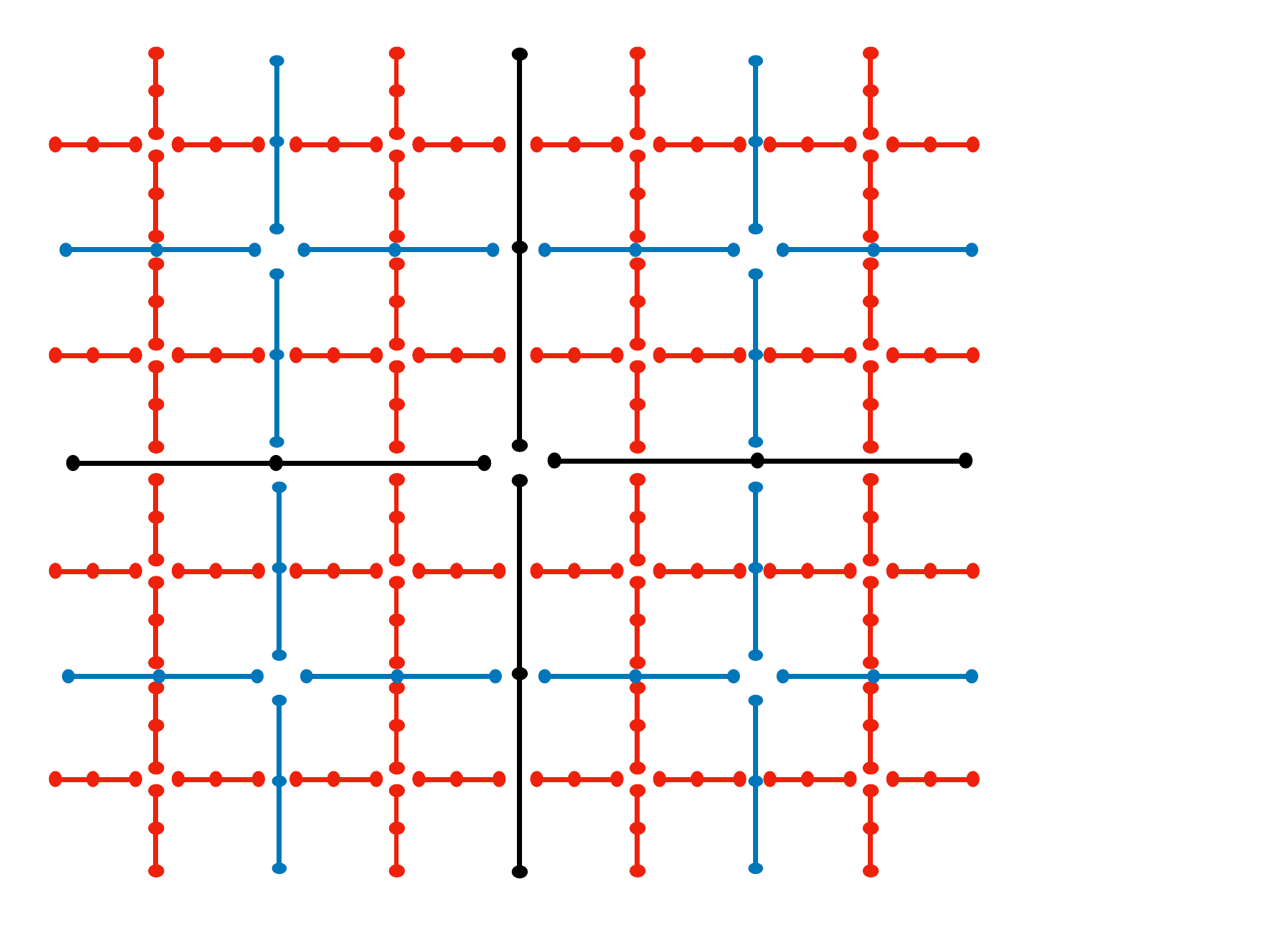}
  }
    \subfigure[Level 2]{
    \includegraphics[trim=8cm 7cm 10cm 7cm,clip=true,width=0.48\textwidth]{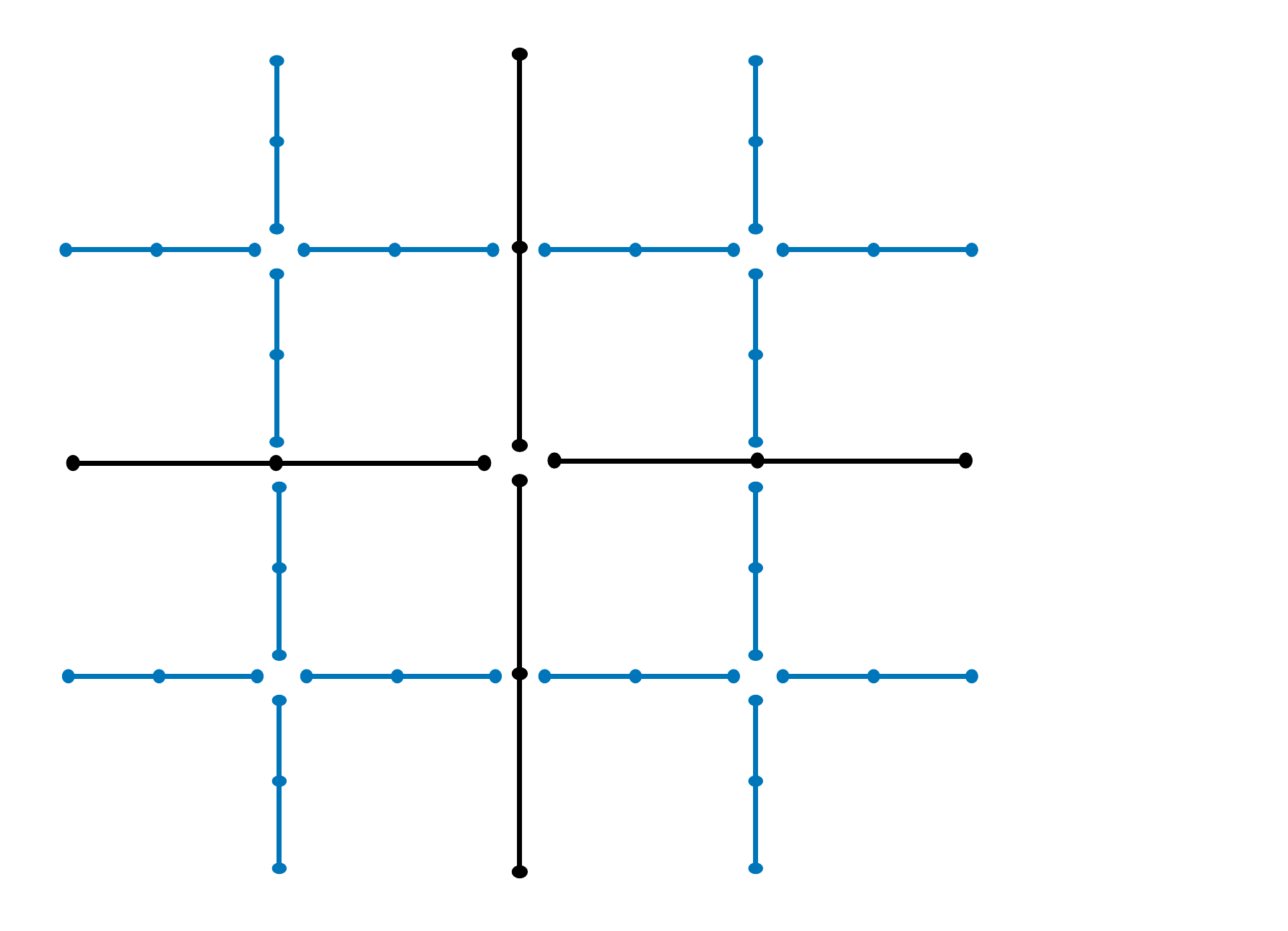}
  }
    \begin{center}
    \subfigure[Level 3]{
    \includegraphics[trim=8cm 7cm 10cm 7cm,clip=true,width=0.48\textwidth]{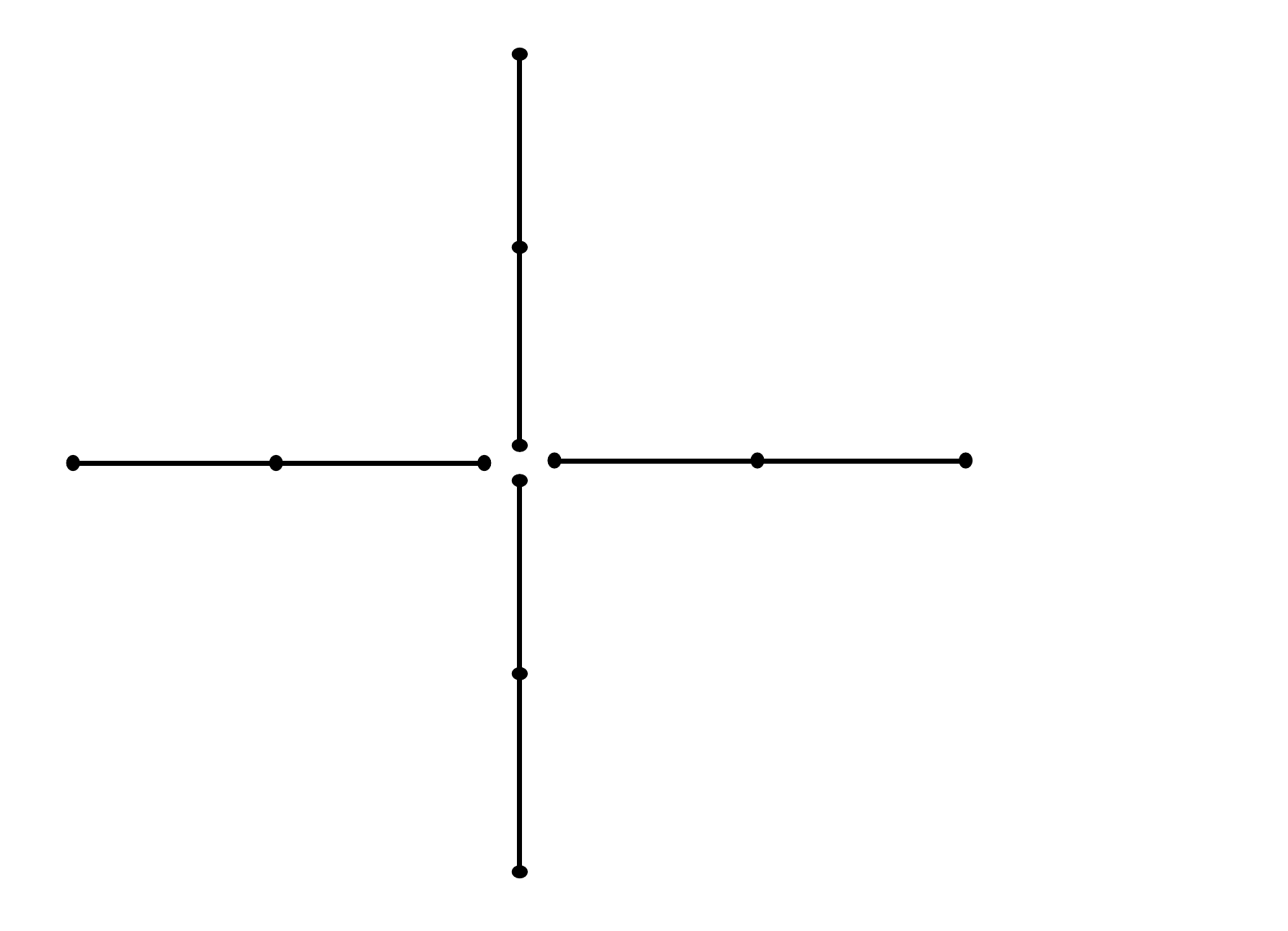}
      \figlab{ML_levels3}
  }
    \end{center}
    \caption{An example of different levels in the multilevel (ML) algorithm. Compared to Figure \figref{ND_levels} for ND, here the 
    separator fronts at all levels have the same number of edges, i.e., 4. Similar to Figure \figref{ND_levels}
    the three circles on each edge represent the nodes corresponding 
    to a solution order of $\p=2$.
      }
  \figlab{ML_levels}
\end{figure}

\begin{figure}[h!b!t!]
    \subfigure[Level 1]{
    \includegraphics[trim=8cm 7cm 10cm 7cm,clip=true,width=0.48\textwidth]{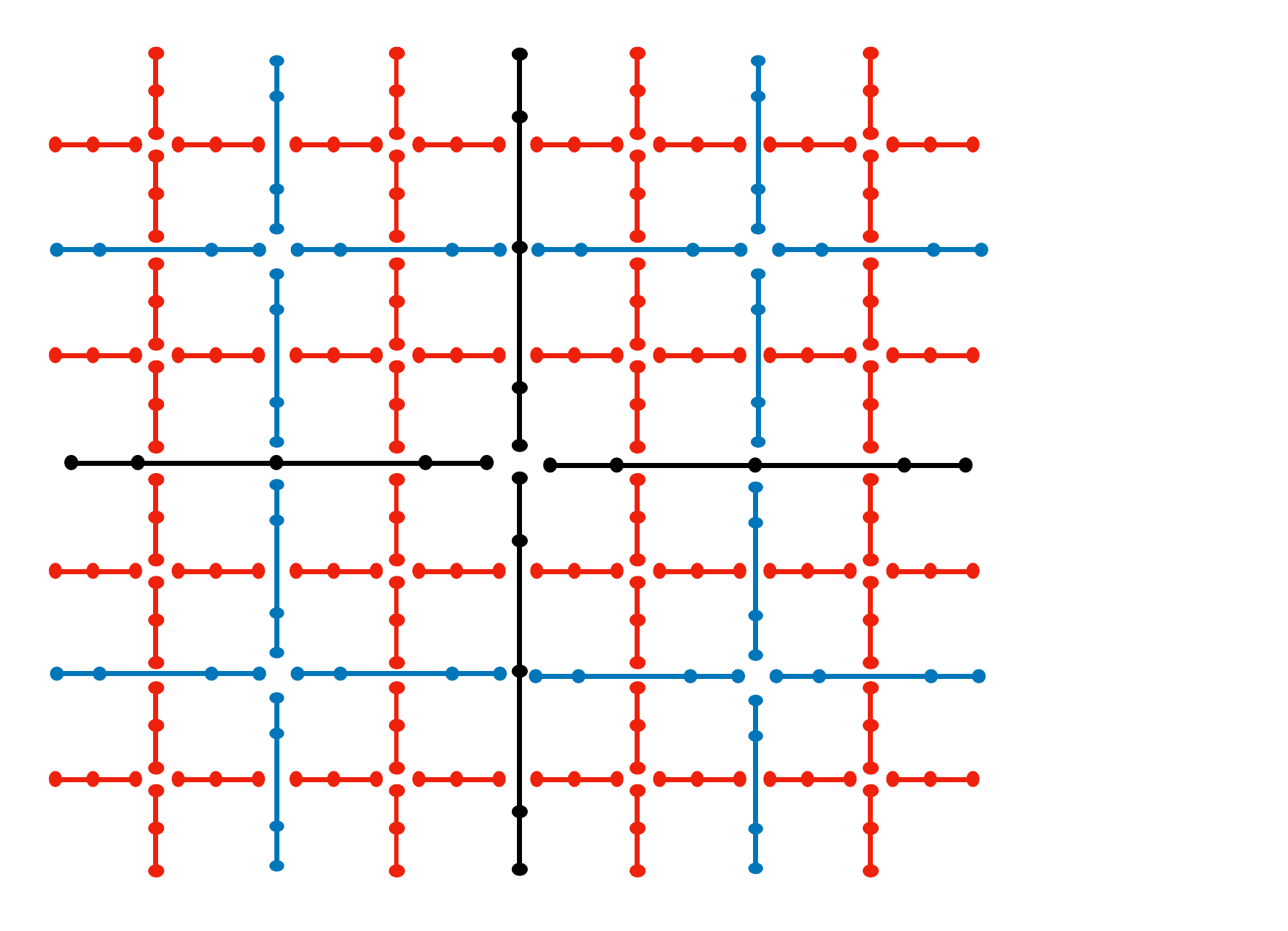}
  }
    \subfigure[Level 2]{
    \includegraphics[trim=8cm 7cm 10cm 7cm,clip=true,width=0.48\textwidth]{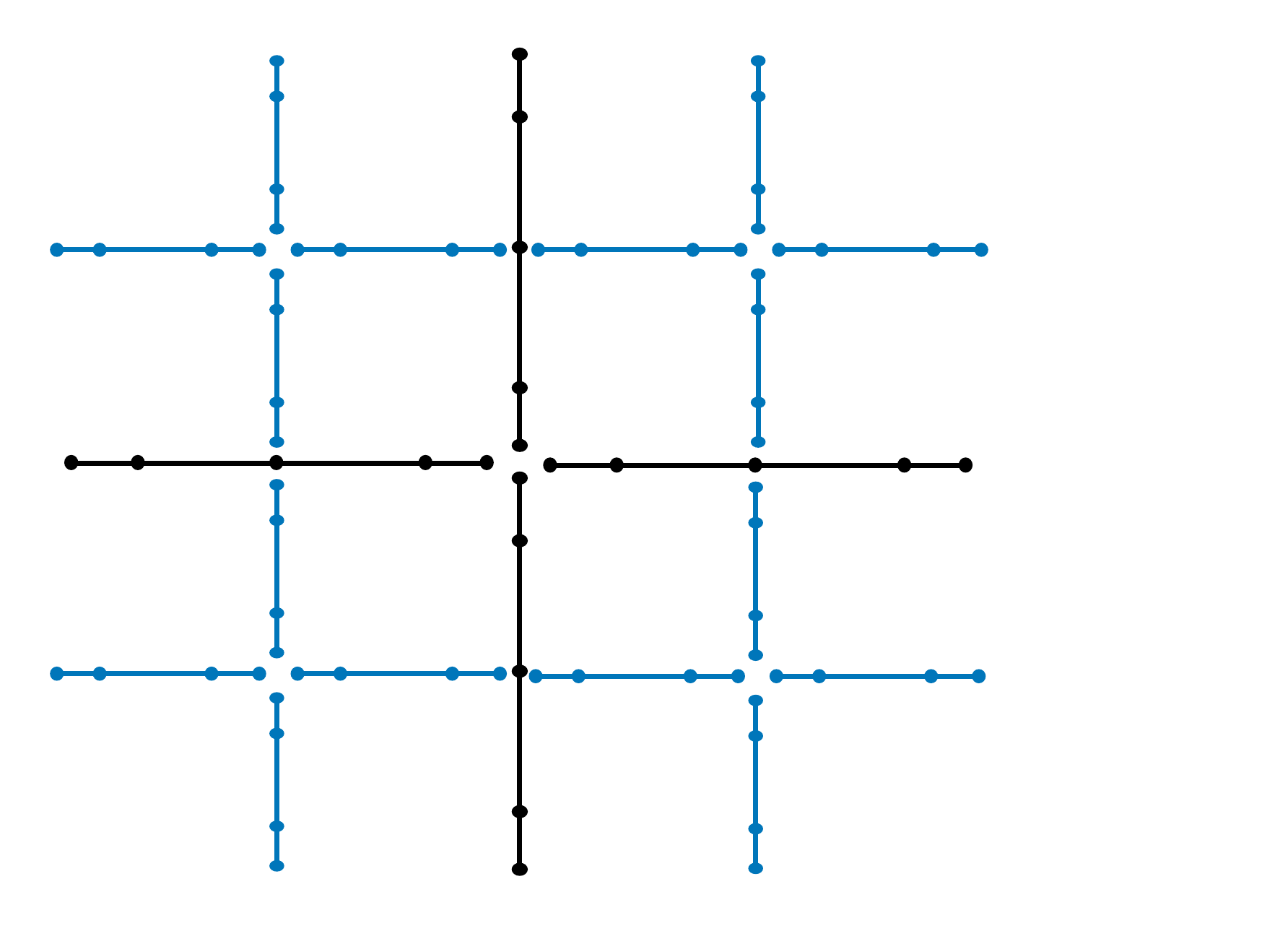}
  }
    \begin{center}
    \subfigure[Level 3]{
    \includegraphics[trim=8cm 7cm 10cm 7cm,clip=true,width=0.48\textwidth]{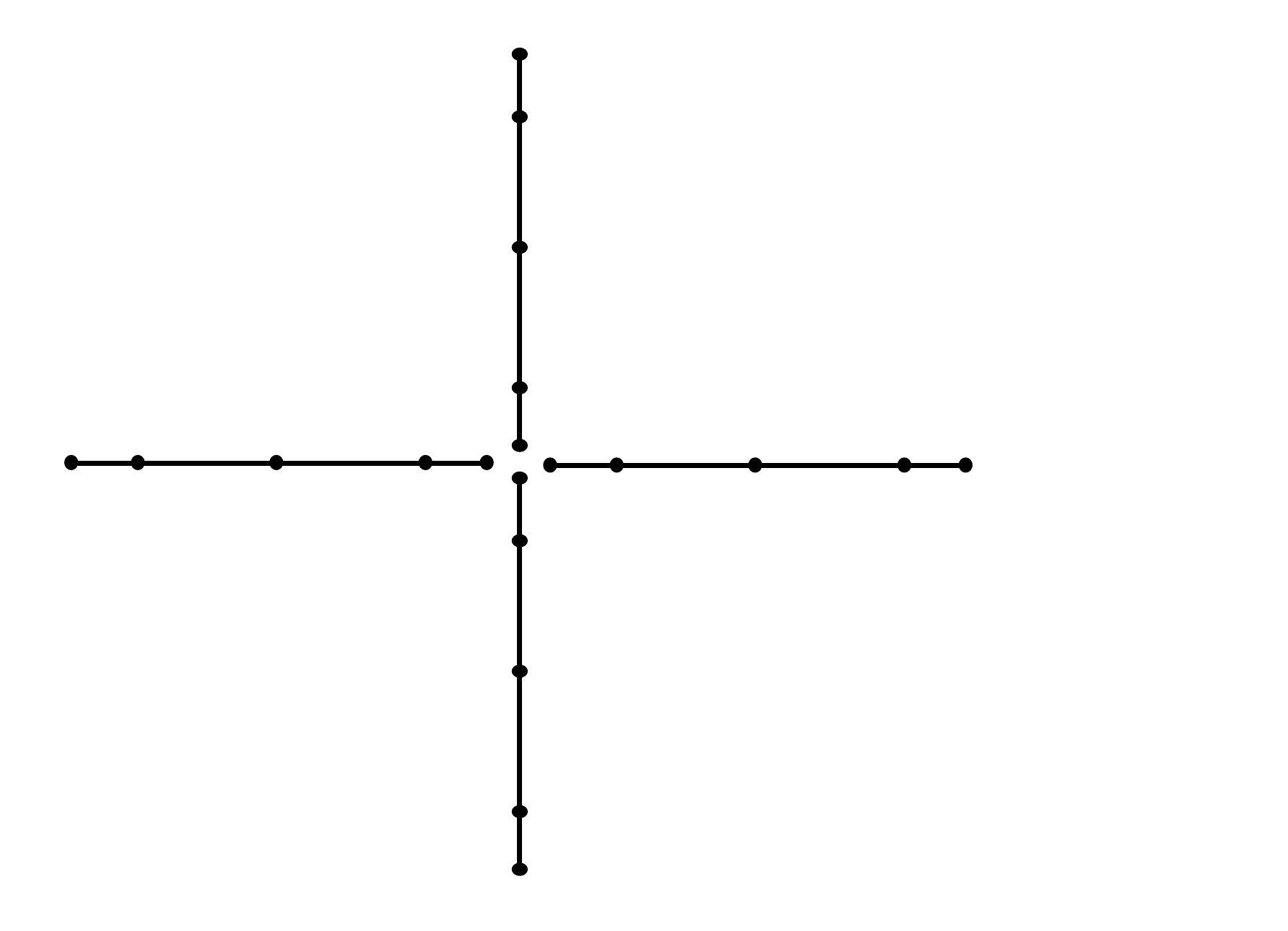}
        \figlab{EML_levels3}
  }
    \end{center}
    \caption{An example of different levels in the enriched multilevel (EML) algorithm.  The 
    number of edges in the separator fronts at all levels is 4. Due to  polynomial enrichment, we have 3 nodes per edge, corresponding to
     $\p=2$, on the red crosses (level 1 separator fronts); 4 nodes per edge, corresponding to $\p=3$, on the blue crosses (level 2 separator fronts); and 
    5 nodes per edge,  corresponding to $\p=4$, on the black cross (level 3 separator front). }
  \figlab{EML_levels}
\end{figure}

\subsection{Combining multilevel approaches with domain decomposition methods}
\seclab{twoLevel}
As discussed in Section \secref{directSolvers}, both ML and EML
strategies are approximations of direct solver. A
natural idea is to use them as ``coarse" scale solvers in a two-level
domain decomposition method
\cite{quarteroni1999domain,smith2004domain,toselli2006domain}. In
particular, either ML or EML approach can be used to capture the smooth
components and to provide global coupling for the algorithm, whereas a
fine scale solver can capture the high-frequency localized error and the small length scale details and sharp
discontinuities. This combined approach can be employed in an iterative manner
as a two-level solver in the domain decomposition
methods.

We select block-Jacobi as our fine scale solver, where each block
corresponds to an edge in the original fine mesh in Figure
\figref{level1_nd}. The reason for this choice is that block-Jacobi is
straightforward to parallelize, insensitive to ordering direction for problems with convection and also reasonably robust with respect to
problem parameters. This is also supported by  our previous work on geometric multigrid
methods \cite{wildey2018unified} for elliptic PDEs, where we compared few different
smoothers and found block-Jacobi to be a competitive choice.  We combine the fine and coarse scale
solvers in a multiplicative way as this is typically more effective than
additive two-level solvers especially for nonsymmetric problems
\cite{toselli2006domain}.

We would like to point out that due to the approximate direct solver characteristic
of our coarse scale solvers, regardless of the nature of the underlying
PDEs, our two-level approaches are applicable. This will be verified in various numerical results in Section
\secref{numerical}. Next we layout the algorithm
for our iterative multilevel approach.

\subsection{Iterative multilevel solvers/preconditioners}

In Figure \figref{multilevel_solver_precond} we show schematic
diagrams of the two-level approaches described in Section
\secref{twoLevel} combining block-Jacobi fine-scale solver and ML or
EML coarse-scale solvers (coarse in the $h-$sense). Algorithm
\ref{al:multilevel_iterative} describes in details every step of these
iterative multilevel solvers and how to implement them.  In
particular, there we present the algorithm for linear problems or
linear systems arising from  Newton linearization or Picard
linearization for nonlinear problems. Note that we perform the
factorizations of both coarse- and fine-scale matrices before the start
of the iteration process so that during each iteration only back
solves are needed.
To precondition the GMRES algorithm we simply use one v-cycle of these iterative
approaches.

\begin{figure}[h!b!t!]
    \subfigure[Two-level solver]{
    \includegraphics[trim=1cm 16cm 1cm 3cm,clip=true,width=0.48\textwidth]{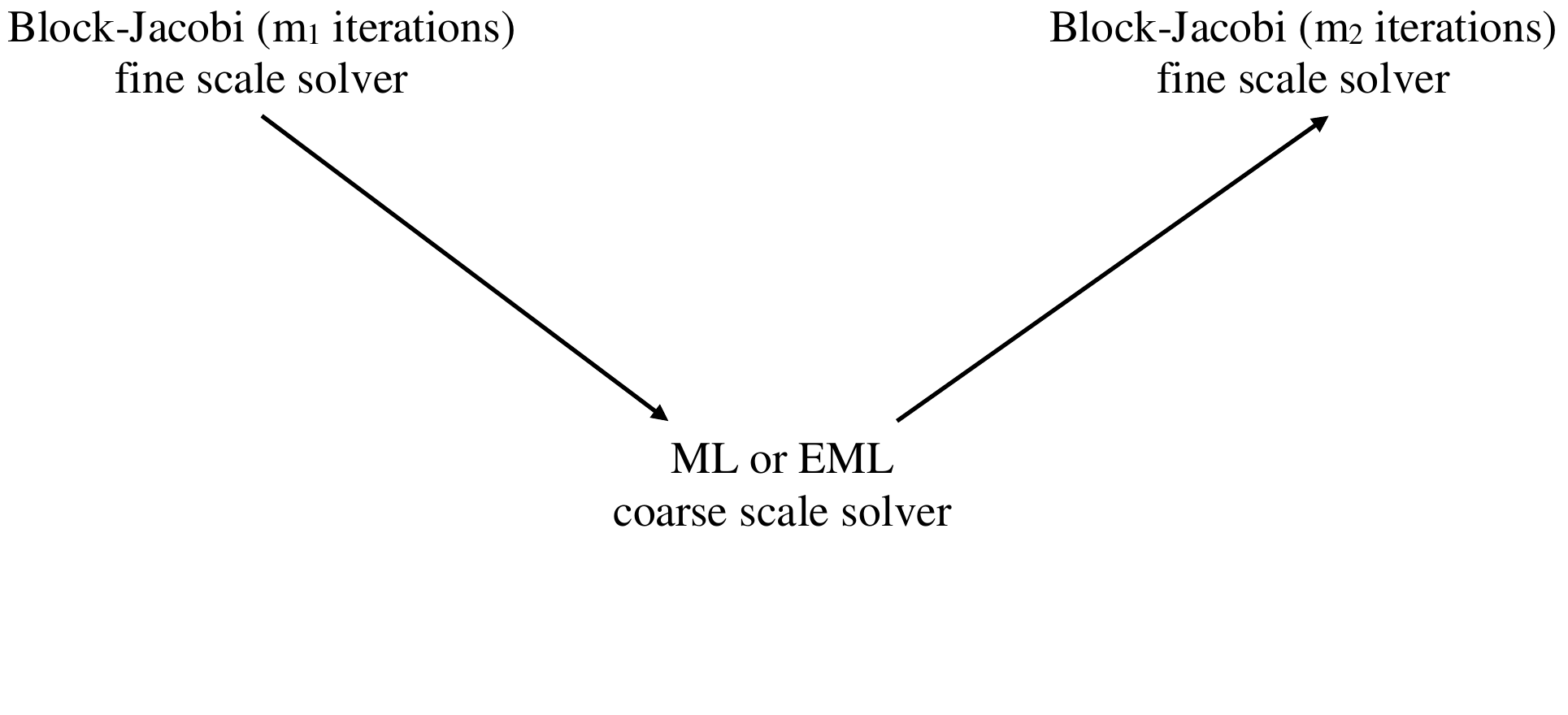}
  }
    \subfigure[Two-level preconditioner]{
    \includegraphics[trim=1cm 16cm 1cm 2.5cm,clip=true,width=0.48\textwidth]{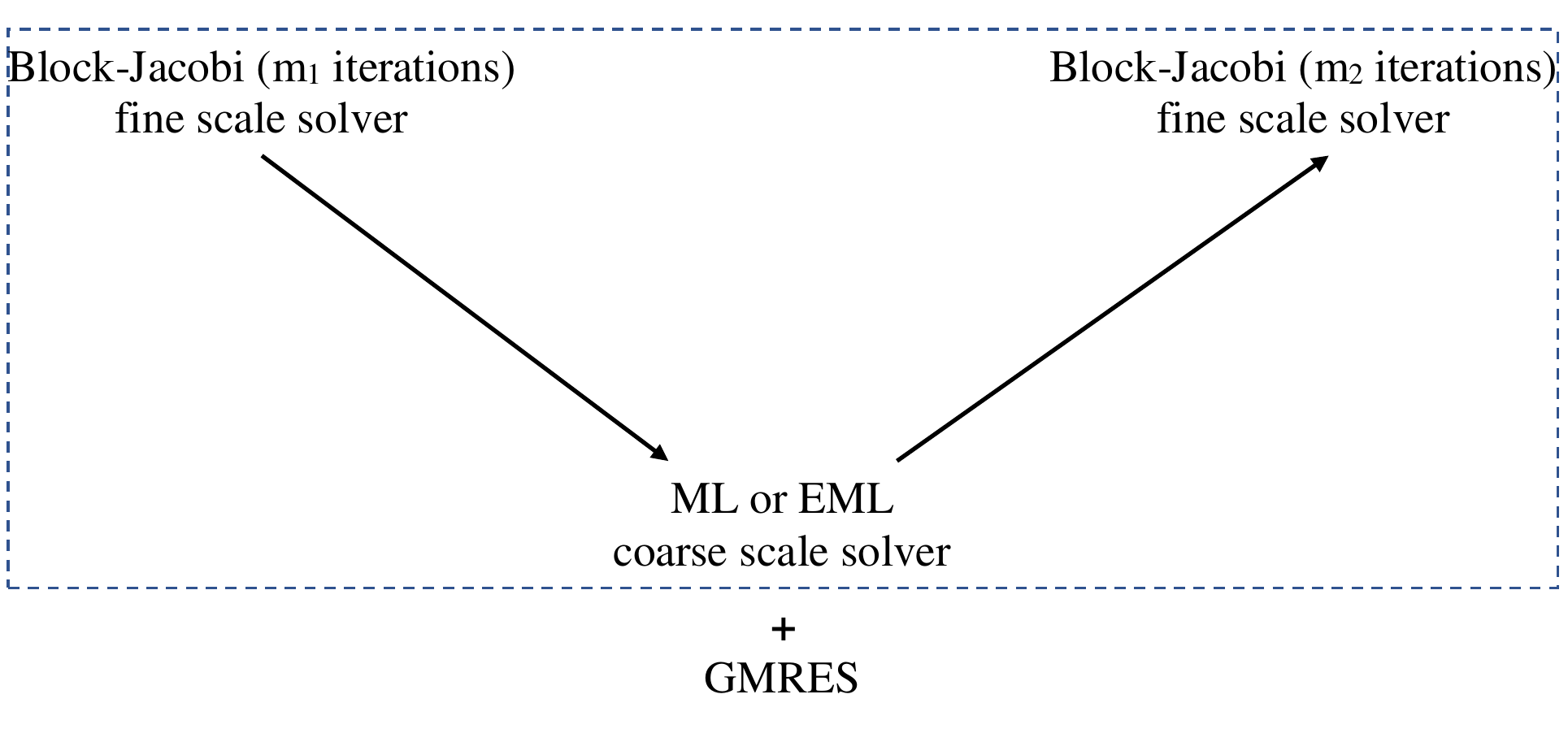}
  }
   \caption{Two-level solvers and preconditioners combining block-Jacobi and ML or EML solvers.}
  \figlab{multilevel_solver_precond}
\end{figure}

\begin{algorithm}
  \begin{algorithmic}[1]
    \STATE Order the unknowns (or permute the matrix) in the nested dissection manner.
    \STATE Construct a set of $L^2$ projections by visiting the edges of the original
      fine mesh skeleton. 
      \STATE Create the level 1 matrices for ML or EML as $A_1=I_0^*AI_0$, where $I_0$ is the projection matrix from level 1 to level 0 and $I_0^*$ is its $L^2$ adjoint.
      \STATE Compute factorizations of level 1 matrices of ML or EML, and 
      the block-Jacobi matrices corresponding to level 0.
      \STATE Compute the initial guess using the coarse scale solver (either ML or EML).
    \WHILE{not converged} 
      \STATE Perform $m_1$ iterations of the block-Jacobi method.
      \STATE Compute the residual.
      \STATE Perform coarse grid correction using either ML or EML.
      \STATE Compute the residual.
      \STATE Perform $m_2$ iterations of the block-Jacobi method.
      \STATE Check convergence. If yes, {\bf exit}, otherwise {\bf set} $i = i+ 1$, and {\bf continue}.
    \ENDWHILE
  \end{algorithmic}
  \caption{An iterative multilevel approach.}
  \label{al:multilevel_iterative}
\end{algorithm}

\subsection{Relationship between iterative multilevel approach and multigrid method}
\seclab{multigrid_interpretation}
In this section we will show that the iterative multilevel approach presented in Algorithm \ref{al:multilevel_iterative} can be
viewed as a multigrid approach with specific prolongation, restriction,
and smoothing operators. To that end, let us consider a sequence of
interface grids $\E_{0} = \Eh,\E_{1},\dots,\E_{N},$
where each $\E_k$ contains the set of edges which remain at level
$k$. Here, $\E_{0}$ is the fine interface grid and $\E_{N}$ is the
coarsest one. Each partition $\E_{k}$ is in turn associated with a
skeletal (trace) space $M_{k}$. We decompose $\E_k$ as $\E_k =
\E_{k,I}\oplus \E_{k,B}$, where $\E_{k,I}$ is the set of interior
edges, corresponding to separator fronts at level $k$, and $\E_{k,B}$
is the set of remaining (boundary) edges. To illustrate this
decomposition, let us consider Figures \figref{ML_levels} and
\figref{EML_levels}. Red edges, blue, and black lumped edges are $\E_{1,I}$, 
$\E_{2,I}$, and $\E_{3,I}$, respectively, and $\E_{k,B} =
\E_{k}\setminus\E_{k,I}$ for $k=1,2,3$. We also decompose the trace
space $M_k$ on $\E_k$ into two parts $M_{k,I}$ and $ M_{k,B}$
corresponding to $\E_{k,I}$ and $\E_{k,B}$,
respectively. Specifically, we require $M_k = M_{k,I}\oplus M_{k,B}$
such that each $\lambda_k \in M_k$ can be uniquely expressed as
$\lambda_k = \lambda_{k,I} + \lambda_{k,B}$, where
\[ \lambda_{k,I} = 
\begin{cases}
 \lambda_k,& \text{ on } M_{k,I},\\
 0,& \text{ on } M_{k,B},
\end{cases}
\quad \text{ and } \quad
\lambda_{k,B} = 
\begin{cases}
 0,& \text{ on } M_{k,I},\\
 \lambda_k,& \text{ on } M_{k,B}.
\end{cases}
\]
The spaces $M_k$ for ML algorithm is given by
\[ M_k = \LRc{\eta\in \mc{Q}^{\p}(e), \forall e \in \E_k} \quad \text{for} \quad k=0,1,2,\dots,N,\]
whereas for EML algorithm it is given by
\[ M_k = 
\begin{cases}
\text{for} \quad  k=0\\
    \LRc{\eta\in \mc{Q}^{\p}(e), \forall e \in \E_k}\\
\text{for} \quad  k=1,2,\dots,N\\
    \LRc{\eta\in \mc{Q}^{\min\LRp{\p+(k-1),10}}(e), \forall e \in \E_{k,I}} \\
    \LRc{\eta\in \mc{Q}^{\min\LRp{{\p+k},10}}(e), \forall e \in \LRc{\E_{k,B}\subset\E_{k+1,I}}}\\
    \LRc{\eta\in \mc{Q}^{\min\LRp{{\p+k+1},10}}(e), \forall e \in \LRc{\E_{k,B}\subset\E_{k+2,I}}}\\
    \vdots\\ 
    \LRc{\eta\in \mc{Q}^{\min\LRp{{\p+k+N-2},10}}(e), \forall e \in \LRc{\E_{k,B}\subset\E_{k+N-1,I}}}.
\end{cases}
\]
If the trace system at level 0 is given by
\begin{equation}
    \label{trace_linear_system}
    A \lambda = g.
\end{equation}
Given the decomposition $M_k = M_{k,I}\oplus M_{k,B}$, the trace system \eqref{trace_linear_system} at the $k$th level
can be written as 
\begin{equation}
  \label{trace_partition}
  A_k \lambda_k = g_k \Leftrightarrow
        \LRs{
     \begin{array}{cc}
       A_{k,II} & A_{k,IB} \\
         A_{k,BI} & A_{k,BB}
     \end{array}
      }
      \LRs{
      \begin{array}{c}
          \lambda_{k,I} \\
          \lambda_{k,B}
      \end{array}
      }
      =
      \LRs{
      \begin{array}{c}
          g_{k,I} \\
          g_{k,B}
      \end{array}
      }.
\end{equation}
We next
 specify the prolongation, restriction and smoothing operators. Since, 
all of the operators except the ones between level 0 and level 1, correspond to ideal operators
in \cite{trottenberg2000multigrid} we explain them briefly here.

\subsection{Prolongation operator}
We define the prolongation operator $I_{k-1}:M_k \to M_{k-1}$ for our
iterative algorithm as
    \begin{equation}
        \I_{k-1}:=
          \begin{cases}
              \Pi_0 &\quad \text{for} \quad k=1,\\
          \LRs{
              \begin{array}{c}
               -A_{k,II}^{-1} A_{k,IB} \\
                  \mc{I}_{BB}
                \end{array}
            } &\quad \text{for} \quad k=2,\dots,N.
          \end{cases}
        \label{Ik}     
    \end{equation}
Here, we denote by $\Pi_0$ the $L^2$ projection from $M_1 \to M_0$ and $\mc{I}_{BB}$
 the identity operator on the boundary. Clearly, apart from $k=1$,
the prolongation operator is nothing but the ideal prolongation in 
algebraic multigrid methods \cite{trottenberg2000multigrid} as well as
in the Schur complement multigrid methods \cite{reusken1994multigrid,wagner1997schur,dendy1982black,de1990matrix}.

\subsection{Restriction operator}

We define the restriction operator $Q_k:M_{k-1} \to M_k$ for our iterative algorithm
as 
    \begin{equation}
          Q_k:=
          \begin{cases}
              \Pi_0^* &\quad \text{for} \quad k=1,\\
          \LRs{
              \begin{array}{c}
                  -A_{k,BI} A_{k,II}^{-1} \quad
                  \mc{I}_{BB}
                \end{array}
            } &\quad \text{for} \quad k=2,\dots,N.
          \end{cases}
        \label{Qk}     
    \end{equation}
Here,  $\Pi_0^*$ is the $L^2$ adjoint of $\Pi_0$. Similar to 
prolongation, apart from $k=1$, the restriction operator  is the ideal
restriction operator \cite{trottenberg2000multigrid}. Given the restriction
and prolongation operators, the Galerkin coarse grid operator is constructed
as
\begin{equation}
    \label{coarse}
 A_{k} := \Q_{k}A_{k-1}\I_{k-1}.
\end{equation}
\subsection{Smoothing}

Recall from the Algorithm \ref{al:multilevel_iterative} and Figure
\figref{multilevel_solver_precond} that  we have both
pre- and post-smoothing steps using block-Jacobi at level 0.
Either ML or EML algorithm implicitly provides additional smoothing.
Indeed, let us consider two generic
levels $k$ and $k+1$. At level $k$, given the decomposition \eqref{trace_partition}
we can write the inverse of $A_k$ as \cite{vassilevski2008multilevel,reusken1994multigrid}
\begin{equation}
    \label{block_inverse}
    A_k^{-1} = 
         \LRs{
         \begin{array}{cc}
             A_{k,II}^{-1} & 0 \\
             0 & 0
         \end{array}
         }
         + I_{k}A_{k+1}^{-1}Q_{k+1},   
\end{equation}
where $A_{k+1}=Q_{k+1}A_k I_{k}$ is the Galerkin coarse grid matrix
(it is also the Schur complement of $A_{k,II}$ in
\eqref{trace_partition} \cite{trottenberg2000multigrid}). As can be
seen, the second term on the right hand side of \eqref{block_inverse} is the coarse grid
correction while the first term is the additive smoothing applied only
to the interior nodes. Another way \cite{reusken1994multigrid} to look at this 
 is the following. If the coarse grid correction is
$z_{k+1} = A_{k+1}^{-1}Q_{k+1}[g_{k,I} \quad g_{k,B}]^T$ then the
block-Jacobi smoothing applied only on the interior nodes with initial guess
as $I_k z_{k+1}$ is given by
\begin{equation}
    \label{block_Jacobi}
    \lambda_k = I_k z_{k+1} + 
         \LRs{
         \begin{array}{cc}
             A_{k,II}^{-1} & 0 \\
             0 & 0
         \end{array}
         }
         \LRp{\LRs{
      \begin{array}{c}
          g_{k,I} \\
          g_{k,B}
      \end{array}
      } - 
      A_k I_k z_{k+1}
      }.
\end{equation}
 From the definition of prolongation operator for $k>1$ in \eqref{Ik}
 we see that $A_k I_k z_{k+1} = [0 \quad \circ]^T$, where ``$\circ$" is
 a term that will subsequently be multiplied by $0$, and is thus not relevant for our
 discussion. As a result,
   $\lambda_k$ obtained from \eqref{block_Jacobi} is the same as
 $A_k^{-1}$ 
 acting on $[g_{k,I} \quad
     g_{k,B}]^T$. In other words, the implicit smoothing is
   equivalent to block-Jacobi smoothing on the interior nodes with the
   initial guess as $I_k z_{k+1}$. In AMG literature
   \cite{trottenberg2000multigrid} this is called F-smoothing, where F
   stands for fine nodes. To summarize, the smoothing operator at
   different levels is given by
    \begin{equation}
          G_k:=
          \begin{cases}
              G_0 &\quad \text{for} \quad k=0,\\
        \LRs{
         \begin{array}{cc}
             A_{k,II}^{-1} & 0 \\
             0 & 0
         \end{array}
         }
              &\quad \text{for} \quad k=1,\dots,N,
          \end{cases}
        \label{Gk}     
    \end{equation}
where $G_0$ is the block-Jacobi smoothing operator at level 0 with each
block corresponding to an edge in the original fine mesh. If we denote by $m_{1,k}$ and $m_{2,k}$ the
number of pre- and post-smoothing steps at level $k$ we have
    \begin{equation}
        m_{1,k}:=
          \begin{cases}
              m_1 &\quad \text{for} \quad k=0,\\
              0  &\quad \text{for} \quad k=1,\dots,N,
          \end{cases}
          \quad
        m_{2,k}:=
          \begin{cases}
              m_2 &\quad \text{for} \quad k=0,\\
              1  &\quad \text{for} \quad k=1,\dots,N.
          \end{cases}
       \label{mk}     
    \end{equation}
 Note that instead of post-smoothing inside the ML or EML solver we 
 could also consider pre-smoothing (i.e. $m_{1,k}=1, m_{2,k}=0 \quad \forall k=1,\dots,N$) and the result remains the same
 \cite{trottenberg2000multigrid}. Now with these specific set of
 operators the iterative multilevel algorithm \ref{al:multilevel_iterative}
 is equivalent to the following multigrid v-cycle 
$$\lambda^{i+1} = \lambda^i + B_0(r_0), \quad i = 0,\hdots$$
where $r_0 = g_0 - A_0\lambda^0$. Here, the action of $B_0$ on a function/vector is defined recursively in the multigrid algorithm \ref{al:multigrid} and the initial guess $\lambda^0$ is computed from either  ML or EML solver. In algorithm \ref{al:multigrid}, $k=0,1,\cdots,N-1$, and $G_{k,m_{1,k}}$, $G_{k,m_{2,k}}$ represent the smoother $G_{k}$ with $m_{1,k}$ and $m_{2,k}$ smoothing steps respectively.
\begin{algorithm}
  \begin{algorithmic}[1]
 \STATE {\it Initialization}: \\
	$e^{\{0\}} = 0,$ 
 \STATE {\it Presmoothing}: \\ 
 $e^{\{1\}} = e^{\{0\}} + G_{k,m_{1,k}}\left(r_k-A_ke^{\{0\}}\right),$
 \STATE {\it Coarse Grid Correction}: \\ 
	$e^{\{2\}} = e^{\{1\}} + I_k B_{k+1} \LRp{Q_{k+1}\left(r_k-A_ke^{\{1\}}\right)},$
 \STATE {\it Postsmoothing}: \\ 
      $B_k\LRp{r_k} = e^{\{3\}} = e^{\{2\}} + G_{k,m_{2,k}}\left(r_k-A_ke^{\{2\}}\right).$
\end{algorithmic}
  \caption{Iterative multilevel approach as a v-cycle multigrid algorithm}
  \label{al:multigrid}
\end{algorithm}
At the coarsest level $M_N$, we set $B_N = A_N^{-1}$ and the inversion
is computed using direct solver.
The above multigrid approach trivially satisfies the following relation \cite{trottenberg2000multigrid}
\begin{equation}
    \eqnlab{energyConser}
    \LRa{A_k\I_k\lambda,\I_k\lambda}_{\E_k}=\LRa{A_{k+1}\lambda,\lambda}_{\E_{k+1}} \quad \forall \lambda \in M_{k+1}, 
\end{equation}
where $\LRa{.,.}_{\E_k}$, $\LRa{.,.}_{\E_{k+1}}$ represents the L$^2$ inner product on $\E_k$ and $\E_{k+1}$ respectively.
This
is a sufficient condition for the stability of intergrid transfer operators in a multigrid algorithm \cite{trottenberg2000multigrid}.
The trivialness is due to the fact that: 1) our prolongation and restriction operators are ideal ones
except for $k=1$, for which they are $L^2$ adjoint of each other; and 2)
the coarse grid matrices are constructed by Galerkin projection 
\eqref{coarse}.

\subsection{Complexity estimations}
\seclab{complexity}

In this section we estimate the serial complexity for ND, ML and EML
algorithms in both 2D and 3D. For simplicity, we consider standard
square or cubical domain discretized by $N_{T}=n^d$ quad/hex elements,
where $d$ is the dimension. The total number of levels is
$N=log_2(n)$. Let $\p_k$ be the solution order on separator fronts at
level $k$ and we denote by $q_k = (\p_k+1)^{d-1}$ the number of nodes
on an edge/face. For simplicity, we consider Dirichlet boundary condition and 
exclude the boundary edges/faces in the complexity estimates.

For the ND algorithm, we define
level 0 to be the same as level 1. Following the
analysis in \cite{eijkhout2014introduction}, we have $4^{N-k}$ crosses
(separator fronts) at level $k$ and each front is of size
$4\LRc{\frac{n}{2^{(N+1-k)}}}q_0$ and all matrices are  dense. The
factorization cost of the ND algorithm in 2D is then given by:

{\bf 2D ND algorithm}
\begin{align}
    Factor &= \sum_{k=1}^{N} 4^{(N-k)}\LRp{4\LRc{\frac{n}{2^{(N+1-k)}}}q_0}^3 \\\label{2D_ND_fact}
                &= \mc{O}\LRp{16 q_0^3 N_T^{3/2} \LRs{1-\frac{1}{\sqrt{N_T}}}}.
\end{align}

The memory requirement is given by

\begin{align}
    Memory &= \sum_{k=1}^{N} 4^{(N-k)}\LRp{4\LRc{\frac{n}{2^{(N+1-k)}}}q_0}^2 \\\label{2D_ND_mem}
                &= \mc{O}\LRp{8 q_0^2 N_T log_2\LRp{N_T}}.
\end{align}
As the Schur complement matrices are dense, the cost for the back
solve is same as that for memory. Similarly the estimates in 3D are
given as follows:

{\bf 3D ND algorithm}

\begin{align}
    Factor &= \sum_{k=1}^{N} 8^{(N-k)}\LRp{12\LRc{\frac{n}{2^{(N+1-k)}}}^2q_0}^3 \\\label{3D_ND_fact}
                &= \mc{O}\LRp{ 31 q_0^3 N_T^{2} \LRs{1-\frac{1}{N_T}}}.
\end{align}

\begin{align}
    Memory &= \sum_{k=1}^{N} 8^{(N-k)}\LRp{12\LRc{\frac{n}{2^{(N+1-k)}}}^2q_0}^2 \\\label{3D_ND_mem}
                &= \mc{O}\LRp{18 q_0^2 N_T^{4/3} \LRs{1-\frac{1}{N_T^{1/3}}}}.
\end{align}

Unlike \cite{poulson2012fast}, here we have not included the
factorization and memory costs for the matrix multiplication
$A_{k,BI}A_{k,II}^{-1}A_{k,IB}$. The reason is that the asymptotic
complexity for ND, ML, and EML is unaffected by this additional cost.
For EML in particular, the inclusion of this cost makes the analysis much more complicated
because of the different solution orders involved at different
levels. As shall be shown, our numerical results in section
\secref{Poisson_smooth} indicate that the asymptotic estimates derived
in this section are in good agreement with the numerical results.

As ML is a special case of EML with zero enrichment, it is sufficient to show the estimates
for EML. In this case 
 we still have $4^{N-k}$ fronts at level $k$ and each front is of the size $4q_k$.
The factorization and memory costs in 2D are then given by:

{\bf 2D EML algorithm}
\begin{align}
    Factor &= \sum_{k=1}^{N} 4^{(N-k)}\LRp{4 q_{k}}^3 \\
           &= 64 q_0^3 \sum_{k=1}^{N} 4^{(N-k)}\alpha_k^3 \\\label{2D_EML_fact}
           &= \mc{O}\LRp{64 q_0^3 \LRc{\frac{1}{4}\LRp{1+\frac{\alpha_N^3}{3}}N_T - \frac{\alpha_N^3}{3}}}.
\end{align}

\begin{align}
    Memory &= \sum_{k=1}^{N} 4^{(N-k)}\LRp{4 q_{k}}^2 \\
           &= 16 q_0^2 \sum_{k=1}^{N} 4^{(N-k)}\alpha_k^2 \\\label{2D_EML_mem}
           &= \mc{O}\LRp{16 q_0^2 \LRc{\frac{1}{4}\LRp{1+\frac{\alpha_N^2}{3}}N_T - \frac{\alpha_N^2}{3}}}.
\end{align}

Similarly in 3D we have:

{\bf 3D EML algorithm}

\begin{align}
    Factor &= \sum_{k=1}^{N} 8^{(N-k)}\LRp{12 q_{k}}^3 \\
           &= 1728 q_0^3 \sum_{k=1}^{N} 8^{(N-k)}\alpha_k^3 \\\label{3D_EML_fact}
           &= \mc{O}\LRp{1728 q_0^3 \LRc{\frac{1}{8}\LRp{1+\frac{\alpha_N^3}{7}}N_T - \frac{\alpha_N^3}{7}}},
\end{align}

\begin{align}
    Memory &= \sum_{k=1}^{N} 8^{(N-k)}\LRp{12 q_{k}}^2 \\
           &= 144 q_0^2 \sum_{k=1}^{N} 8^{(N-k)}\alpha_k^2 \\\label{3D_EML_mem}
           &= \mc{O}\LRp{144 q_0^2 \LRc{\frac{1}{8}\LRp{1+\frac{\alpha_N^2}{7}}N_T - \frac{\alpha_N^2}{7}}}.
\end{align}

Here, $\alpha_k=\frac{q_k}{q_0}$. To enable a direct
comparison with ND,  we have taken $\alpha_k=\alpha_N, \quad \forall
k=1,2,\cdots,N$.  As a result, the actual cost is less than the
estimated ones as $\alpha_k < \alpha_N, \quad \forall k<N$. 
Note that either ML or EML iterative algorithm \ref{al:multilevel_iterative} requires additional cost of $\mc{O}\LRp{N_T q_0^3}$ for factorization and of
$\mc{O}\LRp{N_T q_0^2}$ for memory and back solves due to block-Jacobi smoothing.
Since these additional costs
are less than the costs for ML and EML coarse solvers, they increase the
overall complexity of the algorithm by at most a constant, and hence can be omitted.

\begin{remark}
    From the above complexity estimates we observe that the
    factorization cost of our multilevel algorithms scales like
    $\mc{O}\LRp{\alpha_N^3q_0^3N_T}$ in both 2D and 3D. Compared to the
    cost for ND algorithms which is $\mc{O}\LRp{q_0^3N_T^{3/2}}$ in 2D and
    $\mc{O}\LRp{q_0^3N_T^{2}}$ in 3D, a significant gain (independent of
    spatial dimensions) can be achieved using our methods. Similarly,
    the memory cost has reduced to $\mc{O}\LRp{\alpha_N^2q_0^2N_T}$ independent of dimensions as opposed to
    $\mc{O}\LRp{q_0^2N_Tlog_2(N_T)}$ in 2D and $\mc{O}\LRp{q_0^2N_T^{4/3}}$ in 3D
    for the ND algorithm. Here, $\alpha_N=1$ for ML whereas it is
    greater than one for  EML. On the other hand, the
    memory and computational costs required by
    multigrid is typically $\mc{O}\LRp{N_Tq_0}$. 
    Thus the proposed multilevel algorithms are
    $\mc{O}\LRp{\alpha_N^3q_0^2}$ times more expensive in
    computation cost and require $\mc{O}\LRp{\alpha_N^2q_0}$ more  memory
    compared to standard multigrid algorithms. The cost of the
    multilevel algorithms lie in between direct (ND) solvers and
    multigrid solvers.
\end{remark}

\section{Numerical results}
\seclab{numerical} In this section we test the multilevel algorithm
\ref{al:multilevel_iterative} on elliptic, transport, and
convection-diffusion equations. Except for the transport equation in
section \secref{transport}, the domain $\Omega$ is a standard unit square
$[0,1]^2$ discretized with structured quadrilateral
elements. Dirichlet boundary condition is enforced strongly by means
of the trace unknowns on the boundary.
For the transport equation, we take $\Omega=[0,2]^2$
and inflow boundary condition is enforced through trace
unknowns while  outflow boundary
condition is employed on the remaining boundary. The number of levels in the multilevel
hierarchy and the corresponding number of quadrilateral elements are
shown in Table \ref{tab:level_elements_square} and they are used in all
numerical experiments.

\begin{table}[h!b!t!]
\begin{center}
\begin{tabular}{ | r || c c  c  c  c  c c c| }
\hline
    Levels ($N$) &  2 & 3 & 4 & 5 & 6 & 7 & 8 & 9\\
\hline
    Elements & $4^2$ & $8^2$ & $16^2$ & $32^2$ & $64^2$ & $128^2$ & $256^2$ & $512^2$\\
\hline
\end{tabular}
\caption{\label{tab:level_elements_square} The multilevel hierarchy.} 
\end{center}
\end{table}

We note that even though the iterative multilevel algorithm
works as a solver for most of the PDEs considered in this paper,
it either converges slowly or diverges 
for some of the difficult 
cases. Hence
throughout the numerical section we report the iteration counts for
GMRES, preconditioned by one v-cycle of the multilevel algorithm
\ref{al:multilevel_iterative} with the number of block-Jacobi
smoothing steps taken as $m_1=m_2=2$.

The UMFPACK \cite{davis2004algorithm} library is used for the
factorization and all the experiments are carried out in MATLAB in
serial mode. The specifications of the machine used for the
experiments is as follows. The cluster has 24 cores (2 sockets, 12
cores/socket) and 2 threads per core. The cores are Intel Xeon E5-2690
v3 with frequency 2.6 GHz and the total RAM is 256 GB.

The stopping tolerance for the residual is set to be $10^{-9}$ in the
GMRES algorithm. The maximum number of iterations is limited to $200$.
In the tables of subsequent sections by ``*'' we mean that
the algorithm has reached the maximum number of iterations.

\subsection{Elliptic equation}
\subsubsection{Example I: Poisson equation smooth solution}
\seclab{Poisson_smooth} In this section we test the multilevel
algorithm on the Poisson equation with a smooth
exact solution given by
$$u^e=\frac{1}{\pi^2}\sin(\pi x)\cos(\pi y).$$ The forcing is chosen
such that it corresponds to the exact solution, and the exact solution
is used to enforce the boundary conditions. The other parameters in
equation \eqref{model_problem} are taken as ${\bf K}=\mc{I}$, where
$\mc{I}$ is the identity matrix, and $\betab=0$.

In Table \ref{tab:ML_EML_exp1} we show the number
of GMRES iterations preconditioned by one v-cycle of iterative ML and EML
algorithms \ref{al:multilevel_iterative}.  First, the number of iterations for EML is much less than
that for ML and for high-order ($\p > 3$) solutions and fine meshes
EML performs like a direct solver upto the specified tolerance. That
is, the number of preconditioned GMRES iterations is $0$.  This is
expected due to: 1) smooth exact solution, and 2) exponential
convergence of high-order solutions, and thus the initial guess
computed by EML is, within tolerance, same as the direct
solution. As expected for ML, increasing solution order decreases the number of
GMRES iterations; this is again due to the smooth exact solution and the high-order accuracy of HDG. 


In Figures \figref{ml_vs_nd} and \figref{eml_vs_nd} we compare the
time-to-solution of GMRES preconditioned by ML and EML algorithms to the one obtained from
direct solver with nested dissection. 
As can be seen, the ratios (ND/ML and ND/EML)
are greater than one for either large number of elements or high solution
orders. That is, both ML and EML are faster than ND when accurate solution is desirable. 
In Figure \figref{eml_vs_ml}, the EML and ML algorithms are compared
and, apart from $\p=1$, EML is faster than ML though the speedup is
not significant in this example. For high solution orders, ML
behaves more like EML as the mesh is refined.

We now compare the memory usage of ML and EML against ND. Figures
\figref{memory_eml_vs_nd} and \figref{memory_ml_vs_nd} show the
ratio of memory usages (costs) of EML and ML algorithms relatively to ND. We can see that regardless of meshsize and solution order, both ML and EML requires (much) less memory than ND does. In particular,
ML requires almost 8 times less memory than ND at the finest mesh
for any solution order. For EML, memory saving is more significant as the mesh and/or solution order are refined, and EML
is six times less memory demanding than ND with sixth-order solution at the finest mesh size.


Figure \figref{memory_eml_vs_ml} compares the memory usage between EML
and ML. As expected, EML always requires more memory than ML due to
the enrichment. However, the maximum ratio is around 2.1 and since we
limit the maximum enrichment order to $10$, memory requirements at
high orders for both methods are similar.  This is also the reason
that all the curves in Figure \figref{memory_eml_vs_ml} converge to a constant value as
the mesh is refined. 
As the maximum enrichment order can be tuned, depending on the memory
availability of computing infrastructure, one can have the flexibility
to adjust the EML algorithm to adapt to memory and computation demands
of the problem under consideration. For example,  we can perform more
iterations for less memory  to achieve a given accuracy.

Next we verify the complexity estimates derived in section
\secref{complexity}.  Figures \figref{theory_expt_fact},
\figref{theory_expt_solve} and \figref{theory_expt_memory} show that
the numerical results agree well with our estimates except for the
factorization cost of ND in Figure \figref{ND_theory_expt_fact}, which seems to indicate that the asymptotic complexity of
$\mc{O}(N_T^{3/2})$
has not yet been reached. Since the results in Figures
\figref{theory_expt_fact}, \figref{theory_expt_solve} and
\figref{theory_expt_memory} are independent of the PDE under
consideration, in the subsequent sections we study only the iteration
counts. As long as the iterations do not increase much when the
mesh and the solution order are refined, we can obtain a scalable algorithm whose cost can be estimated based on the factorization and memory
costs derived in section \secref{complexity}.

\begin{table}[h!b!t!]
\centering
\begin{tabular}{|c|c|c|c|c|c|c||c|c|c|c|c|c|}
\hline
& \multicolumn{6}{c||}{ML with GMRES} & \multicolumn{6}{c|}{EML
with GMRES}\\
\cline{2-13}
\!\!\! $N$ \!\!\!\! &  \multicolumn{6}{c||}{\!\!\scriptsize  $\p$\!\!} &  \multicolumn{6}{c|}{\!\!\scriptsize $\p$\!\!}\\
\hline
& \scriptsize 1 & \scriptsize 2 & \scriptsize 3 & \scriptsize 4 & \scriptsize 5 & \scriptsize 6 & \scriptsize 1 & \scriptsize 2 & \scriptsize 3 & \scriptsize 4 & \scriptsize 5 & \scriptsize 6 \\
\cline{2-13}
2 & 3  & 3  & 3  & 3  & 2  & 1 & 2  & 3 & 3 & 2 & 2 & 1 \\
3 & 6  & 6  & 5  & 4  & 4  & 2 & 4  & 4 & 4 & 3 & 2 & 1 \\
4 & 10 & 8  & 8  & 6  & 5  & 4 & 7  & 5 & 3 & 2 & 0 & 0 \\
5 & 14 & 11 & 11 & 8  & 7  & 4 & 9  & 5 & 3 & 0 & 0 & 0 \\
6 & 21 & 16 & 16 & 12 & 10 & 6 & 12 & 6 & 2 & 0 & 0 & 0 \\
7 & 31 & 24 & 22 & 17 & 13 & 8 & 16 & 5 & 0 & 0 & 0 & 0 \\
8 & 44 & 34 & 32 & 23 & 19 & 11 & 19 & 3 & 0 & 0 & 0 & 0 \\
9 & 63 & 49 & 45 & 33 & 26 & 16 & 22 & 0 & 0 & 0 & 0 & 0 \\
\hline
\end{tabular}
\caption{\label{tab:ML_EML_exp1} Example I: number of ML- and EML-preconditioned GMRES iterations as the mesh is refined (increasing $N$) and the solution order $p$  increases.}
\end{table}



\begin{figure}[h!b!t!]
    \subfigure[EML versus ND]{
    \includegraphics[trim=0cm 3.1cm 0.25cm 3.5cm,clip=true,width=0.307\textwidth]{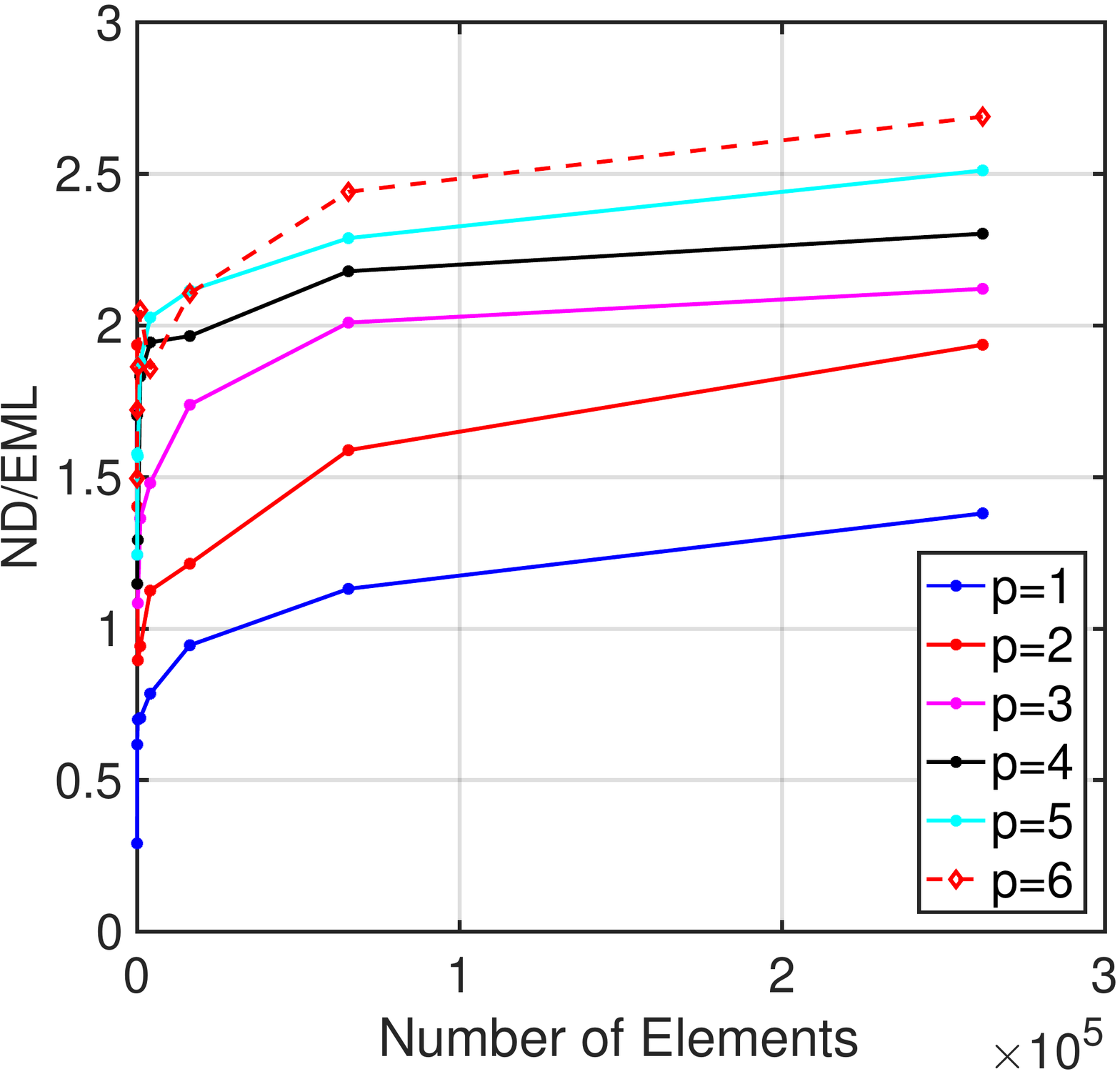}
    \figlab{eml_vs_nd}
  }
    \subfigure[ML versus ND]{
    \includegraphics[trim=0cm 3.7cm 0.25cm 5.3cm,clip=true,width=0.307\textwidth]{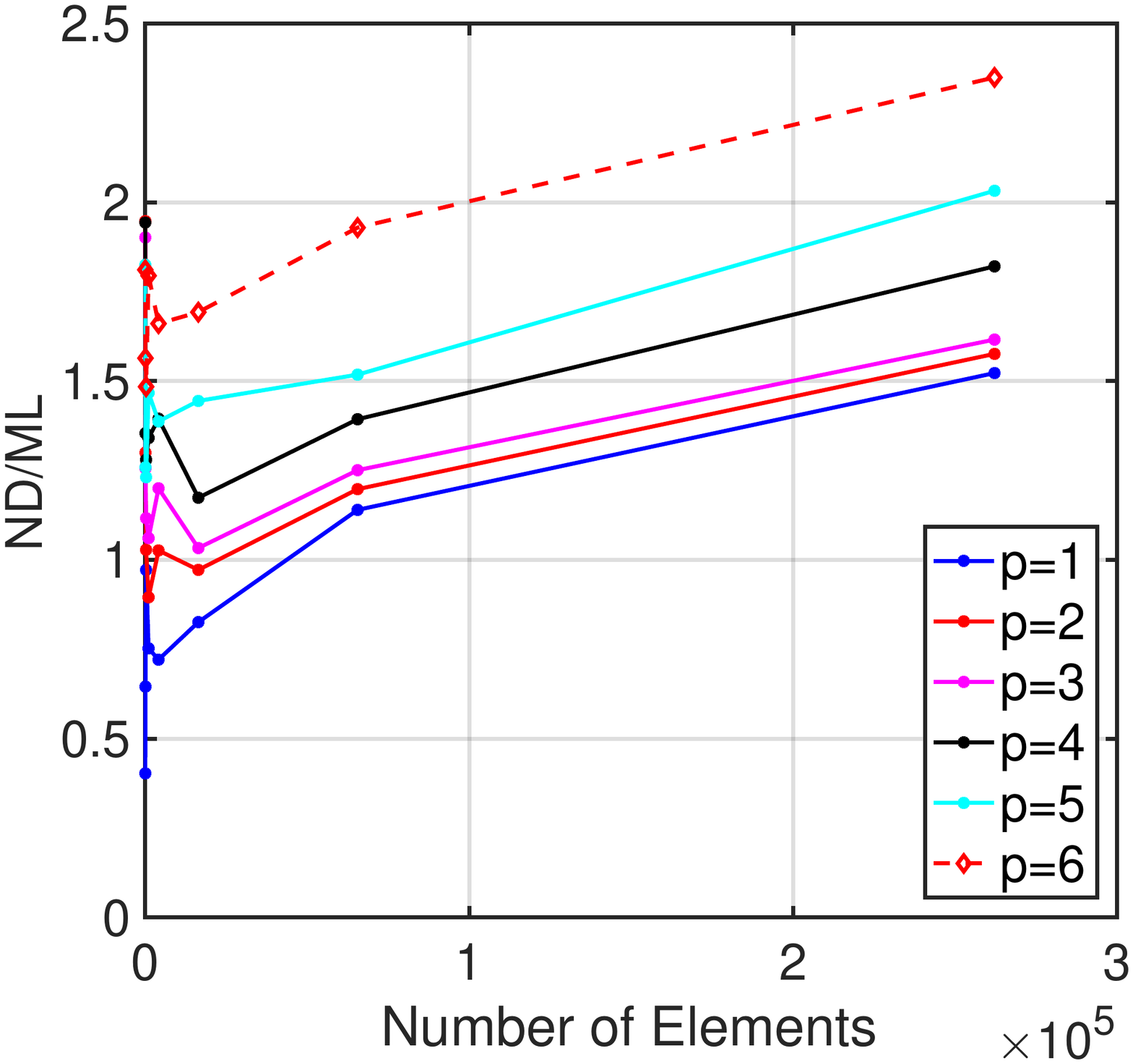}
    \figlab{ml_vs_nd}
  }
    \subfigure[EML versus ML]{
    \includegraphics[trim=0cm 3cm 0.25cm 4cm,clip=true,width=0.307\textwidth]{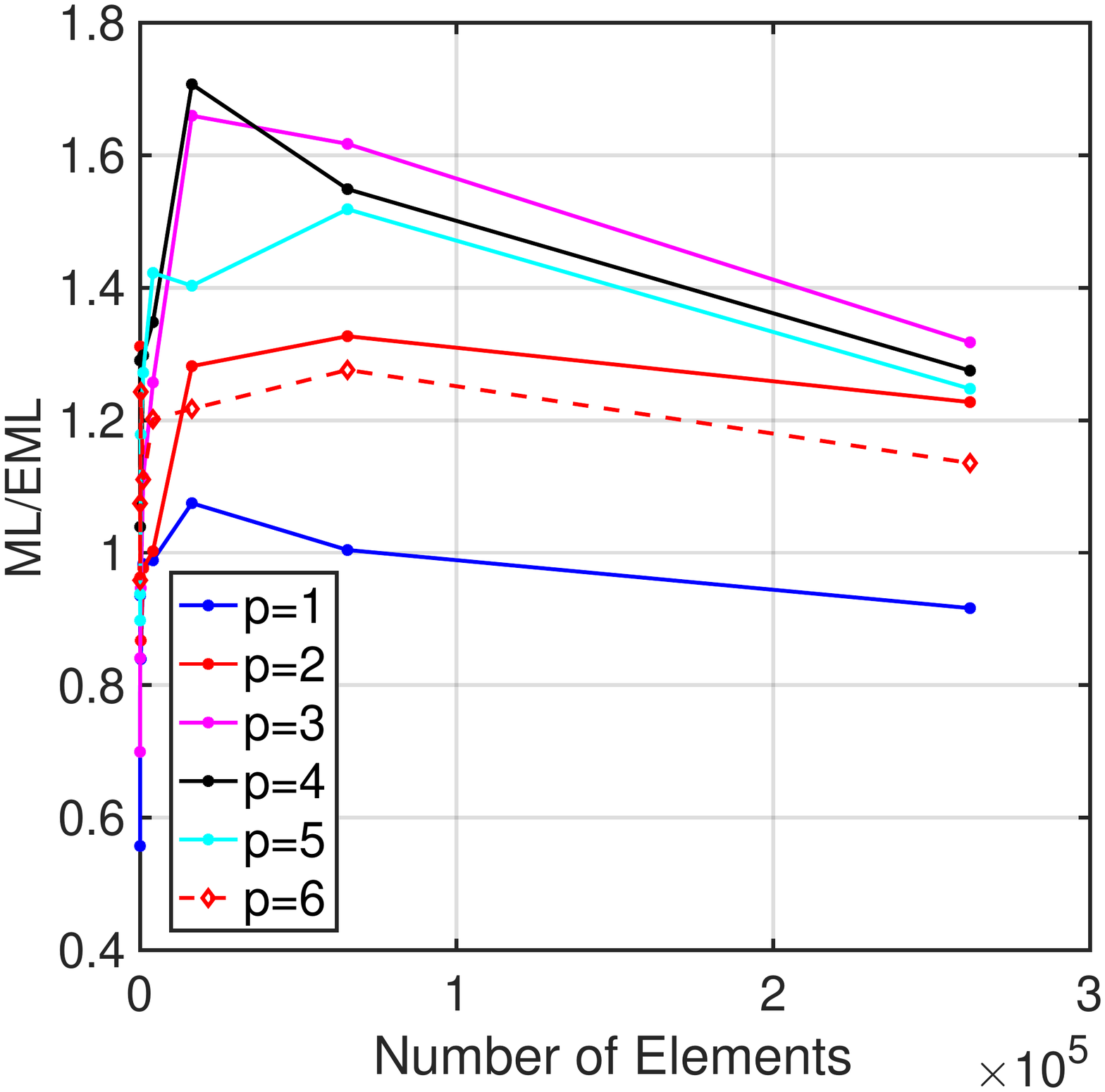}
    \figlab{eml_vs_ml}
  }
    \caption{A comparison of time-to-solution for  EML, ML, and ND algorithms.
    }
  \figlab{Speedup_plots}
\end{figure}

\begin{figure}[h!b!t!]
    \subfigure[EML versus ND]{
    \includegraphics[trim=0cm 3.1cm 0.25cm 3.5cm,clip=true,width=0.307\textwidth]{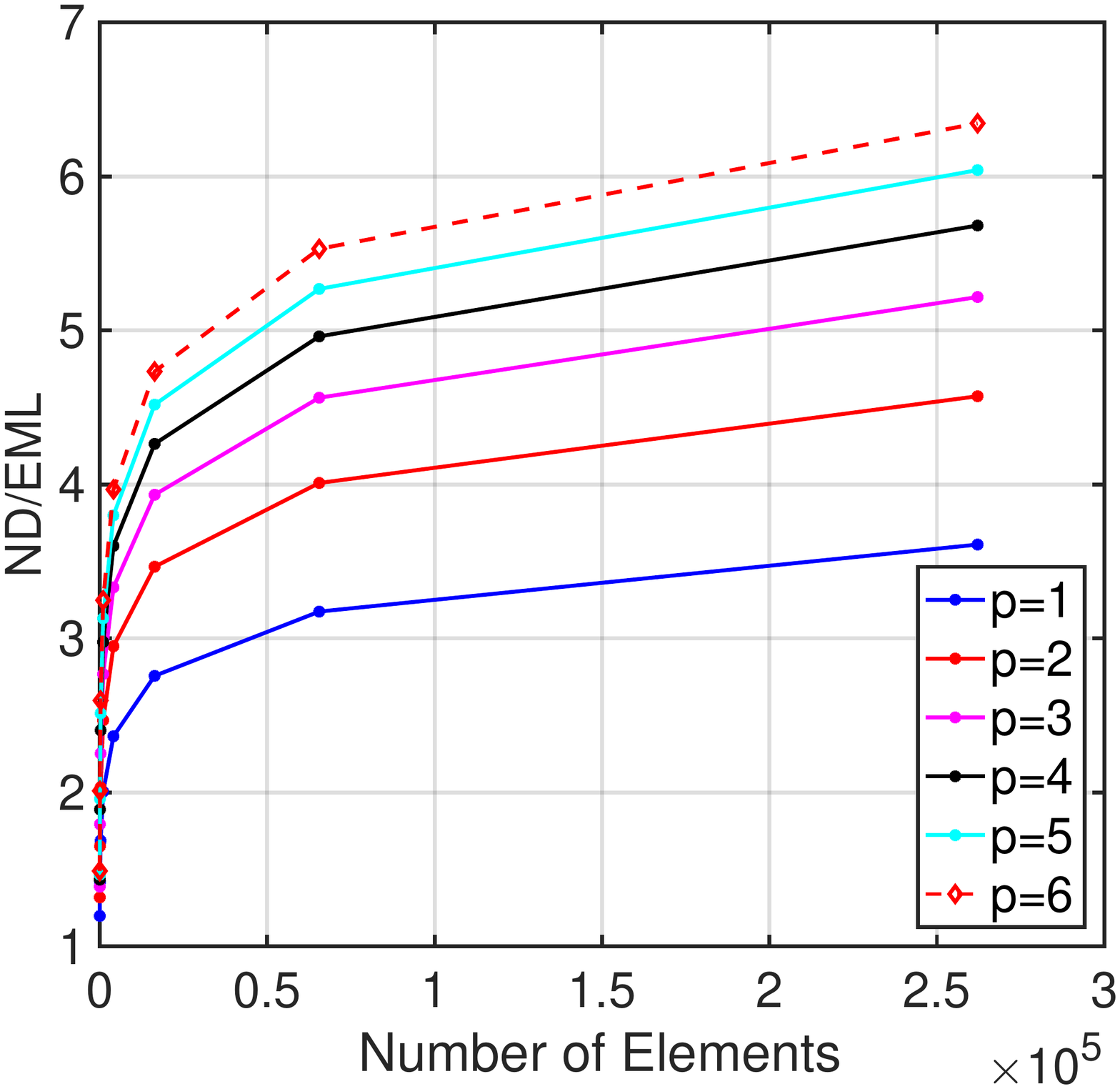}
    \figlab{memory_eml_vs_nd}
  }
    \subfigure[ML versus ND]{
    \includegraphics[trim=0cm 3.1cm 0.25cm 3.5cm,clip=true,width=0.307\textwidth]{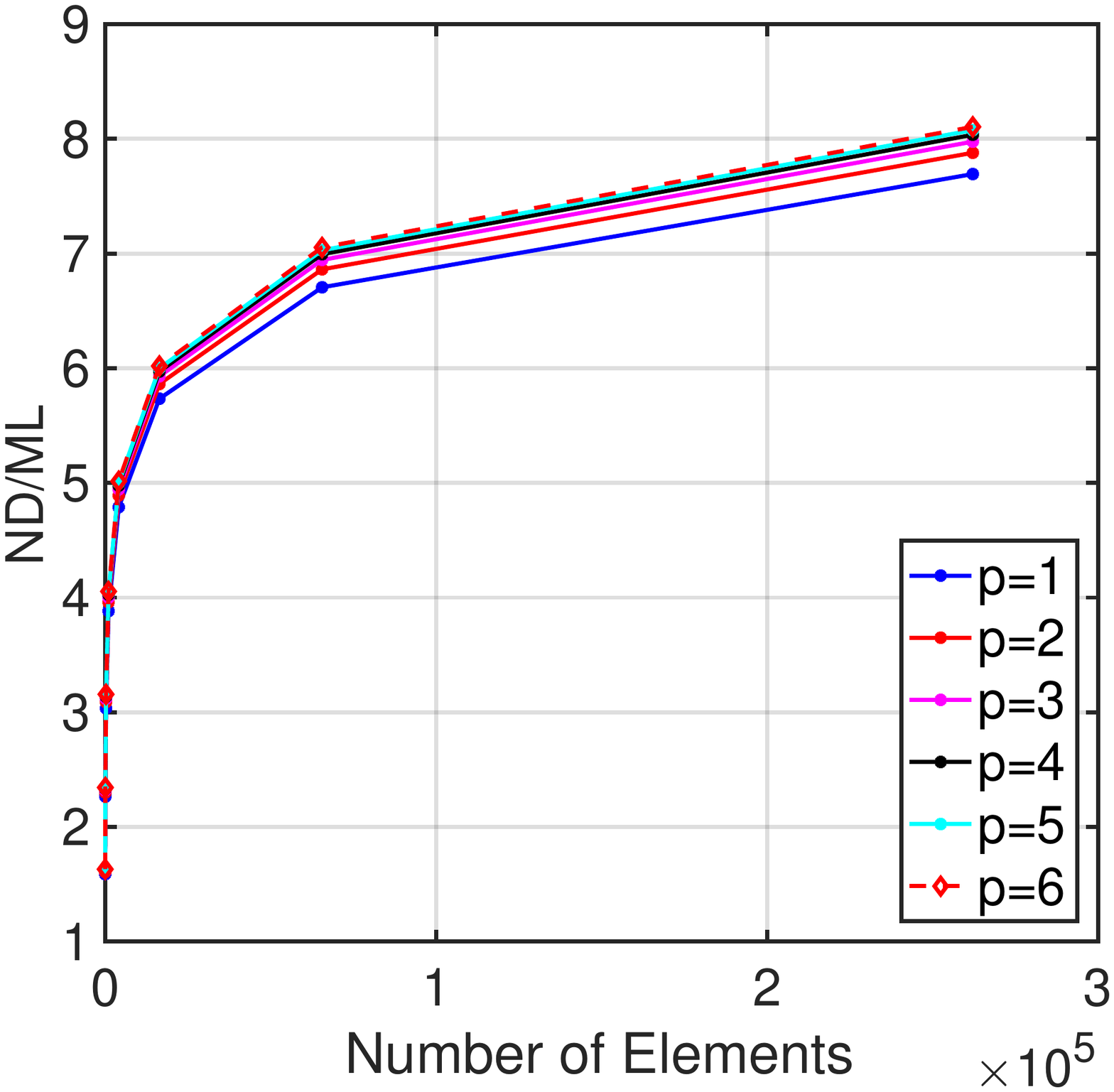}
    \figlab{memory_ml_vs_nd}
  }
    \subfigure[EML versus ML]{
    \includegraphics[trim=0cm 3.1cm 0.25cm 3.5cm,clip=true,width=0.307\textwidth]{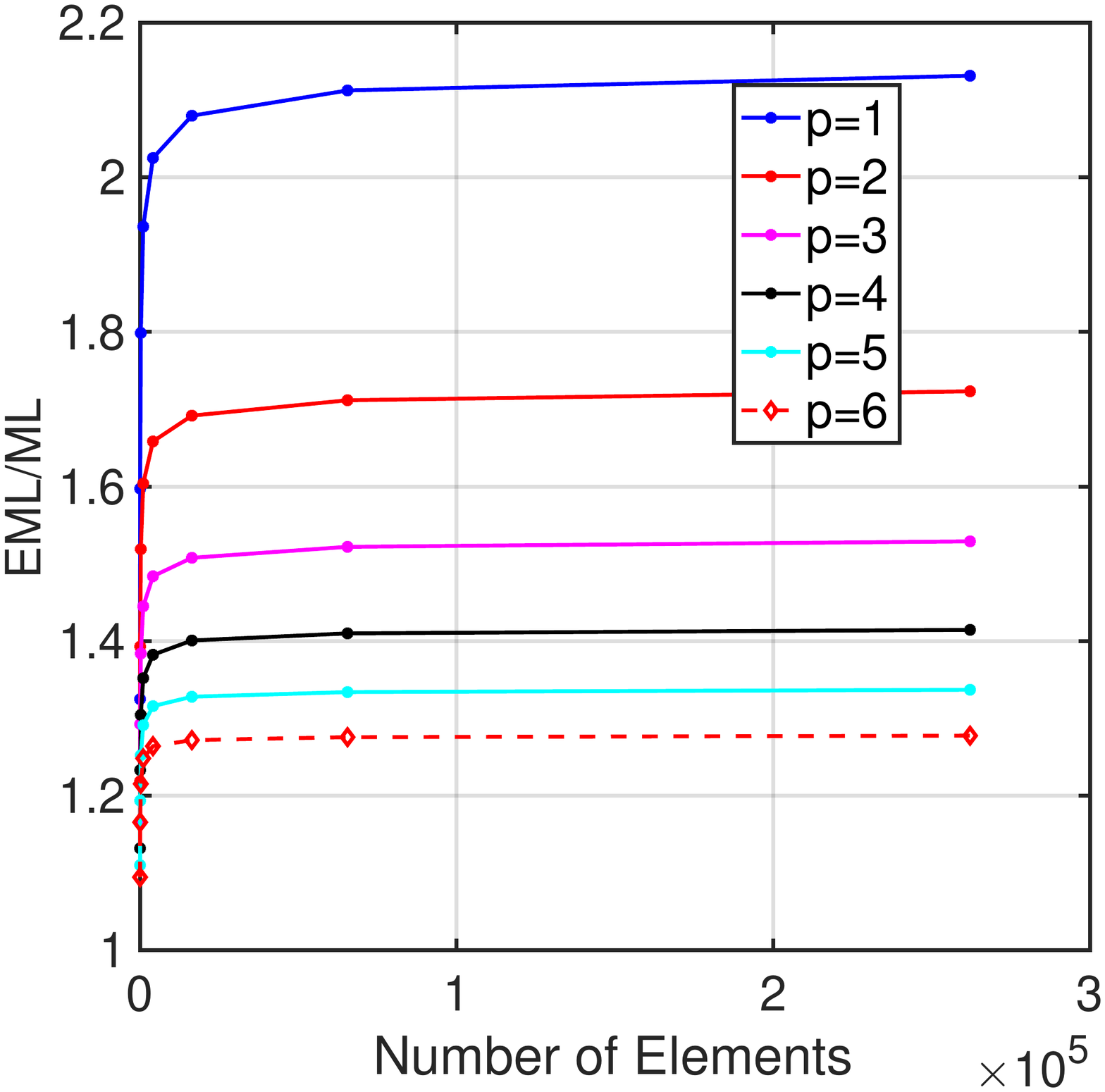}
    \figlab{memory_eml_vs_ml}
  }
    \caption{A comparison of memory requirement for EML, ML, and ND algorithms.}
  \figlab{Memory_plots}
\end{figure}

\begin{figure}[h!b!t!]
    \subfigure[ND]{
    \includegraphics[trim=0.15cm 3.8cm 1.6cm 5.3cm,clip=true,width=0.307\textwidth]{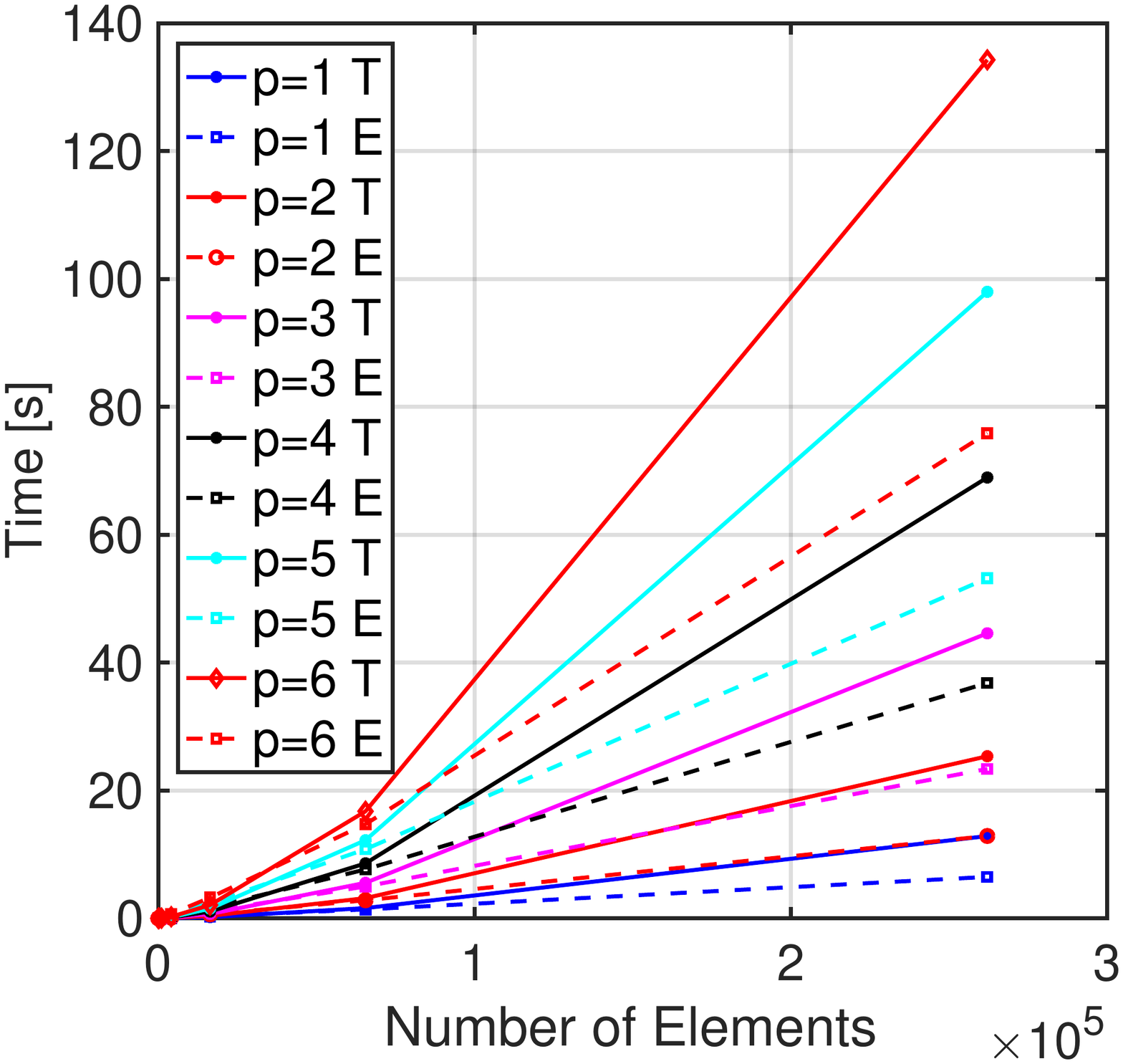}
    \figlab{ND_theory_expt_fact}
  }
    \subfigure[ML]{
    \includegraphics[trim=0.1cm 3.1cm 0.25cm 4cm,clip=true,width=0.307\textwidth]{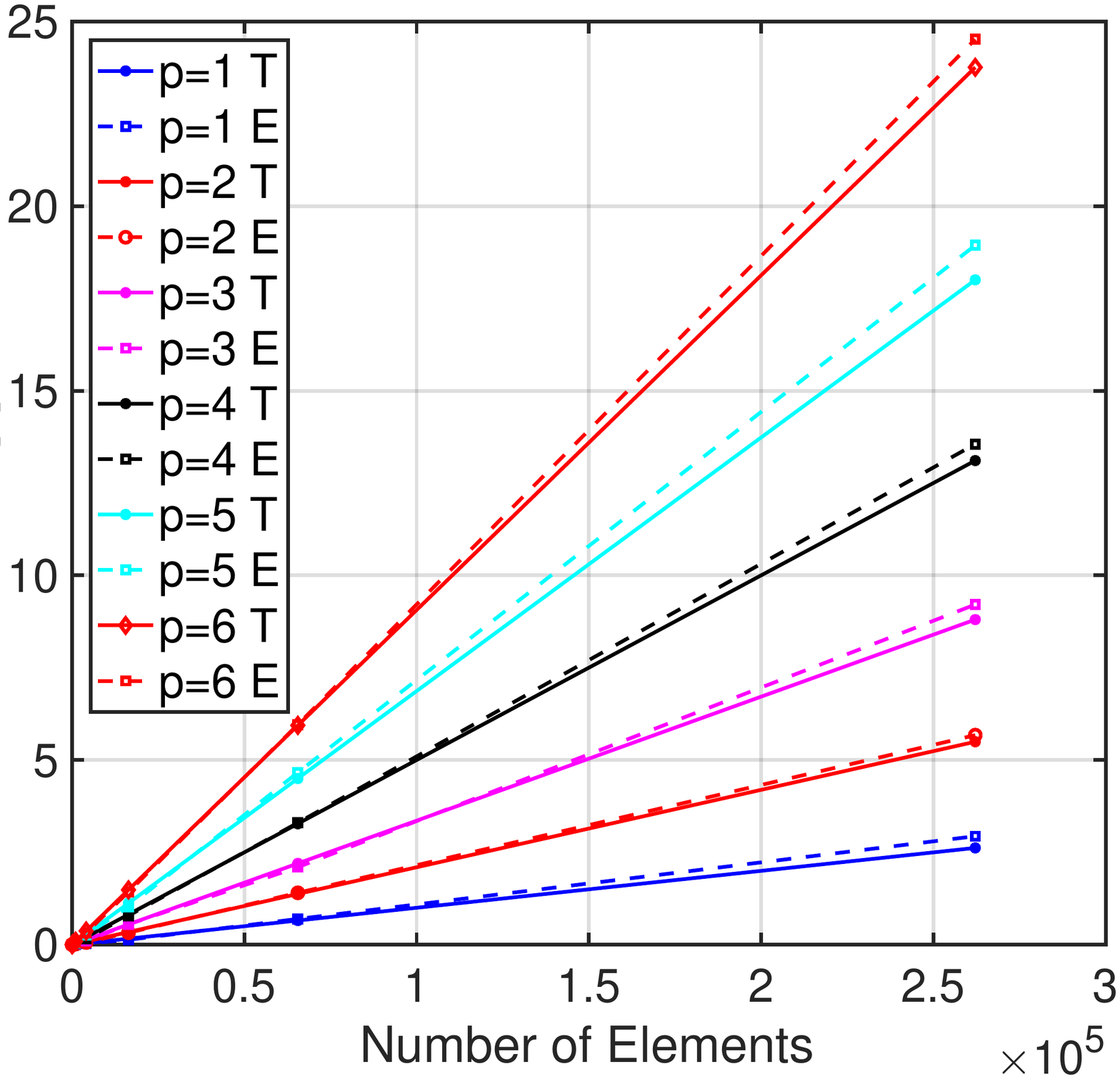}
    \figlab{ML_theory_expt_fact}
  }
    \subfigure[EML]{
    \includegraphics[trim=0.5cm 3.25cm 0.5cm 4cm,clip=true,width=0.307\textwidth]{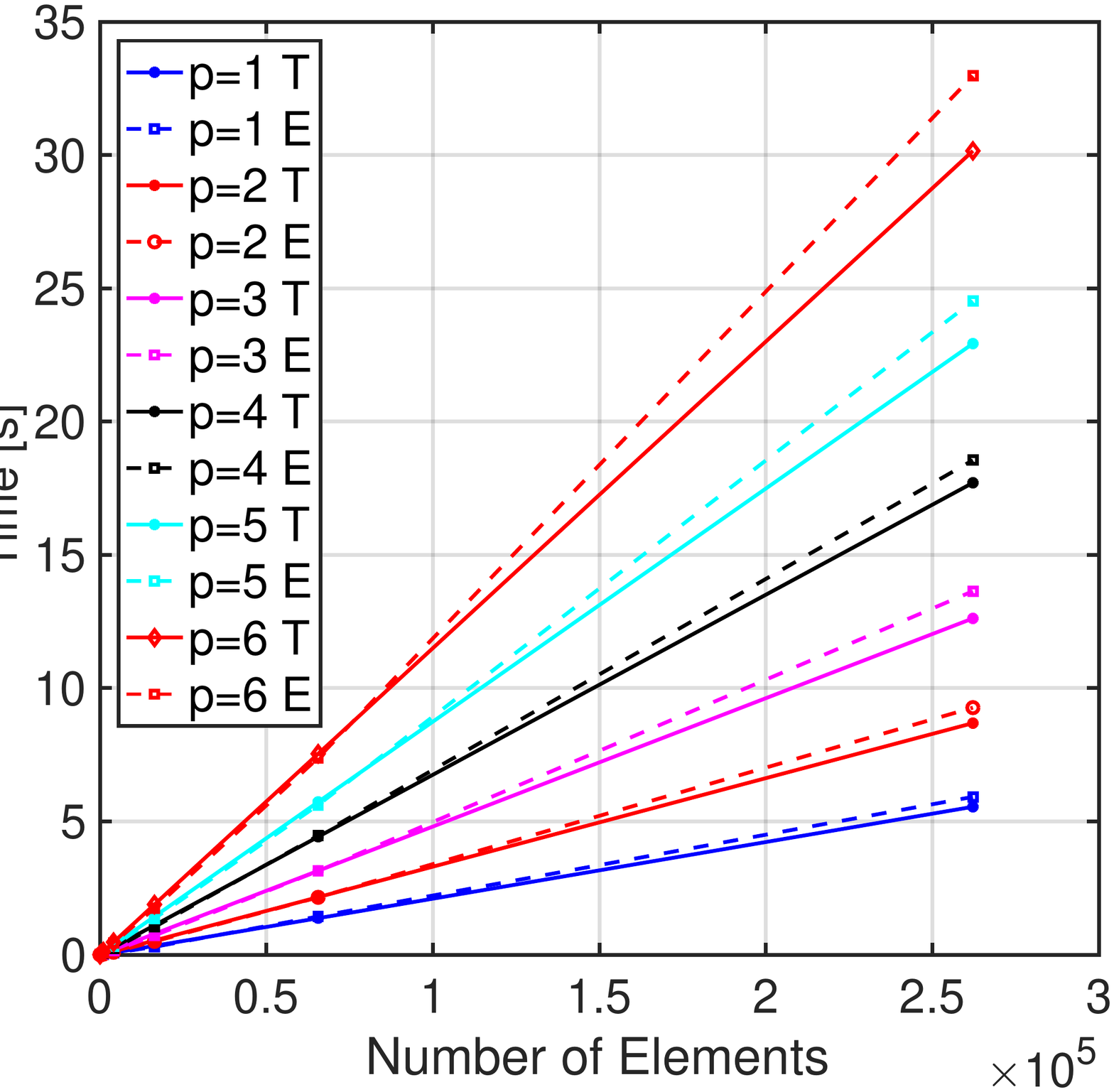}
    \figlab{EML_theory_expt_fact}
  }
    \caption{Asymptotic and numerical estimates of factorization time complexity for EML, ML, and ND. Here, T stands for the 
    theoretically estimated complexity derived in section \secref{complexity} and E for numerical experiment.}
  \figlab{theory_expt_fact}
\end{figure}

\begin{figure}[h!b!t!]
    \subfigure[ND]{
    \includegraphics[trim=0.15cm 3.8cm 1.4cm 5.1cm,clip=true,width=0.307\textwidth]{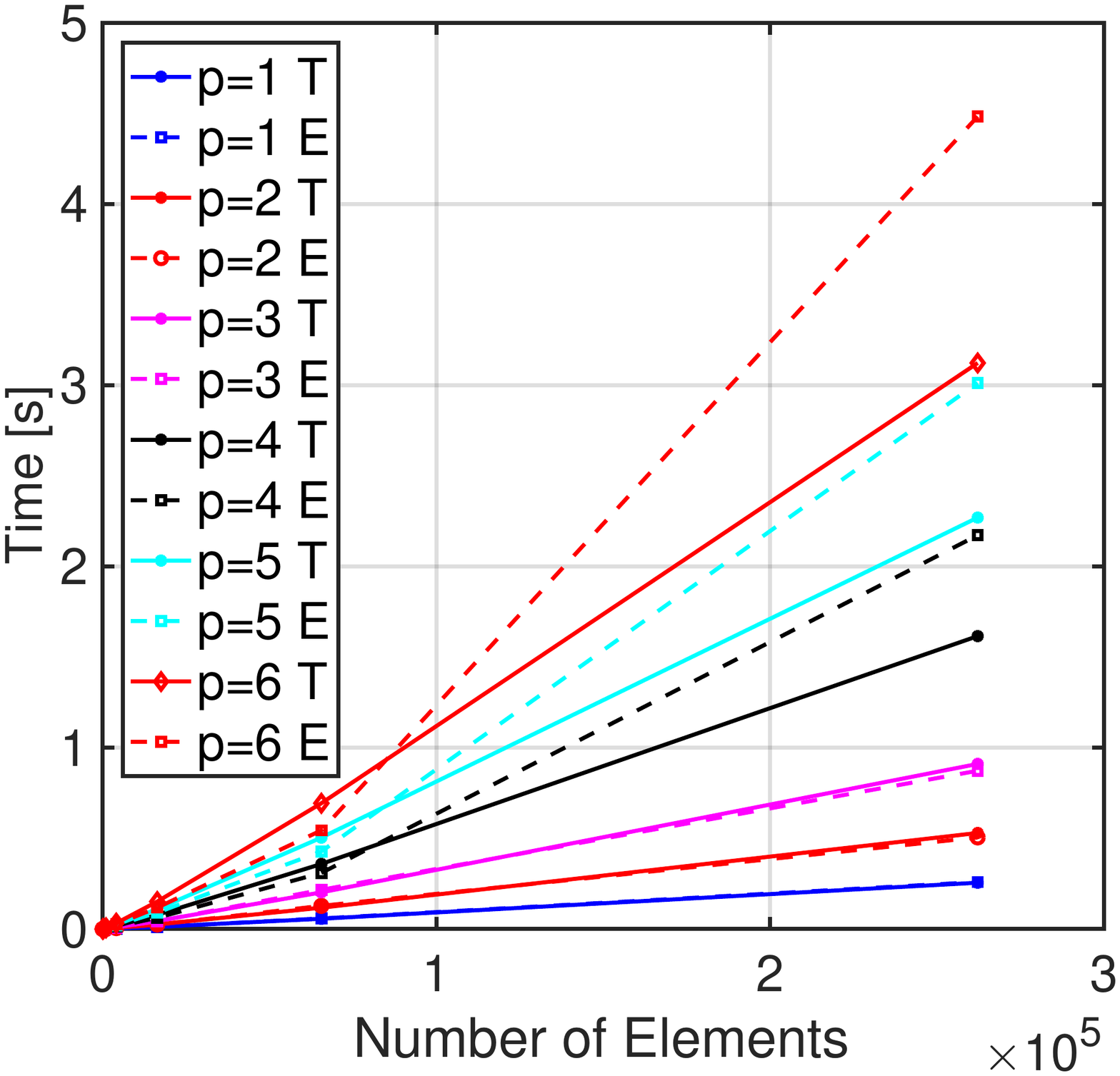}
  }
    \subfigure[ML]{
    \includegraphics[trim=0cm 3.1cm 0.25cm 4cm,clip=true,width=0.307\textwidth]{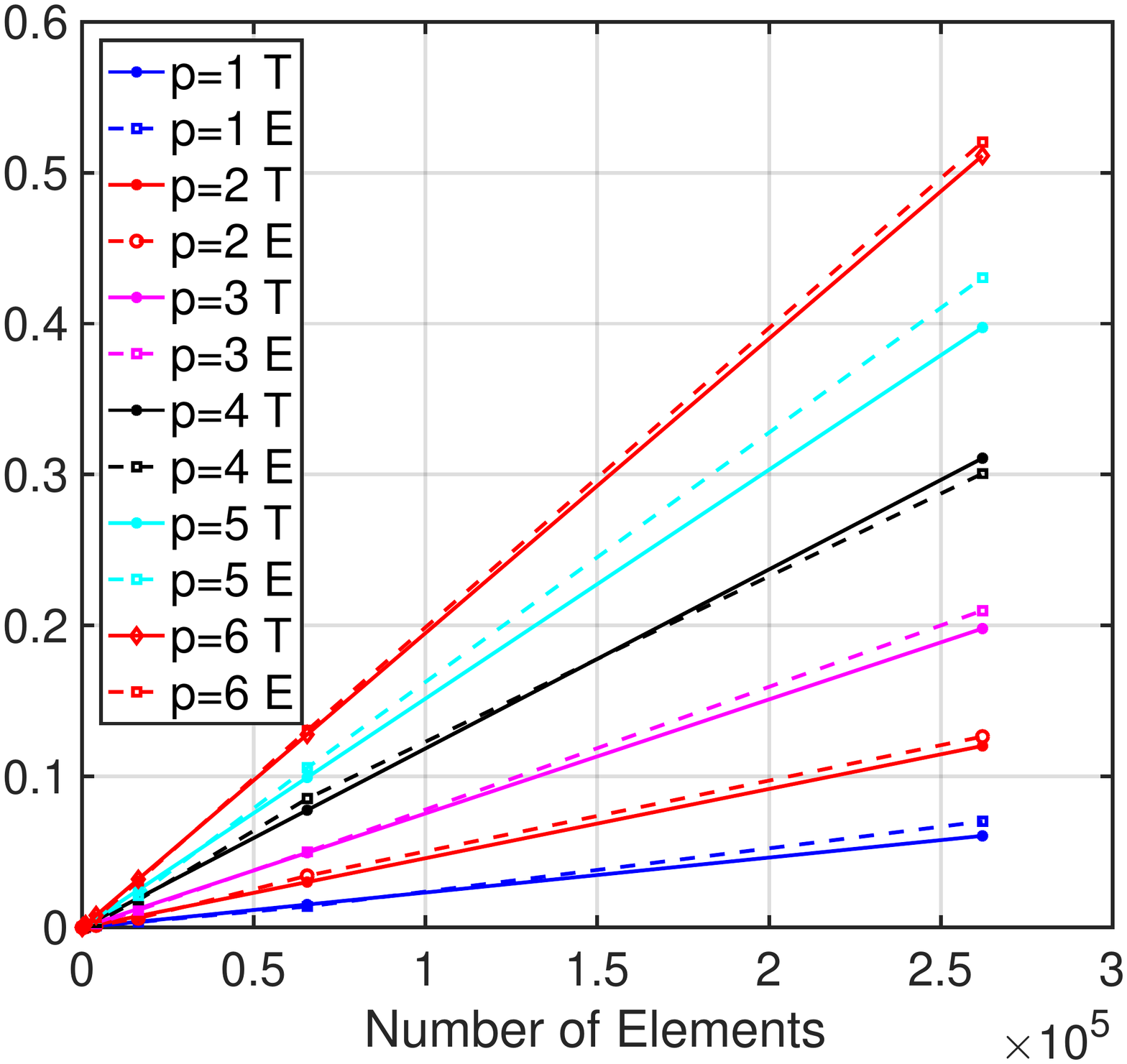}
  }
    \subfigure[EML]{
    \includegraphics[trim=0.4cm 3.1cm 0.25cm 4cm,clip=true,width=0.307\textwidth]{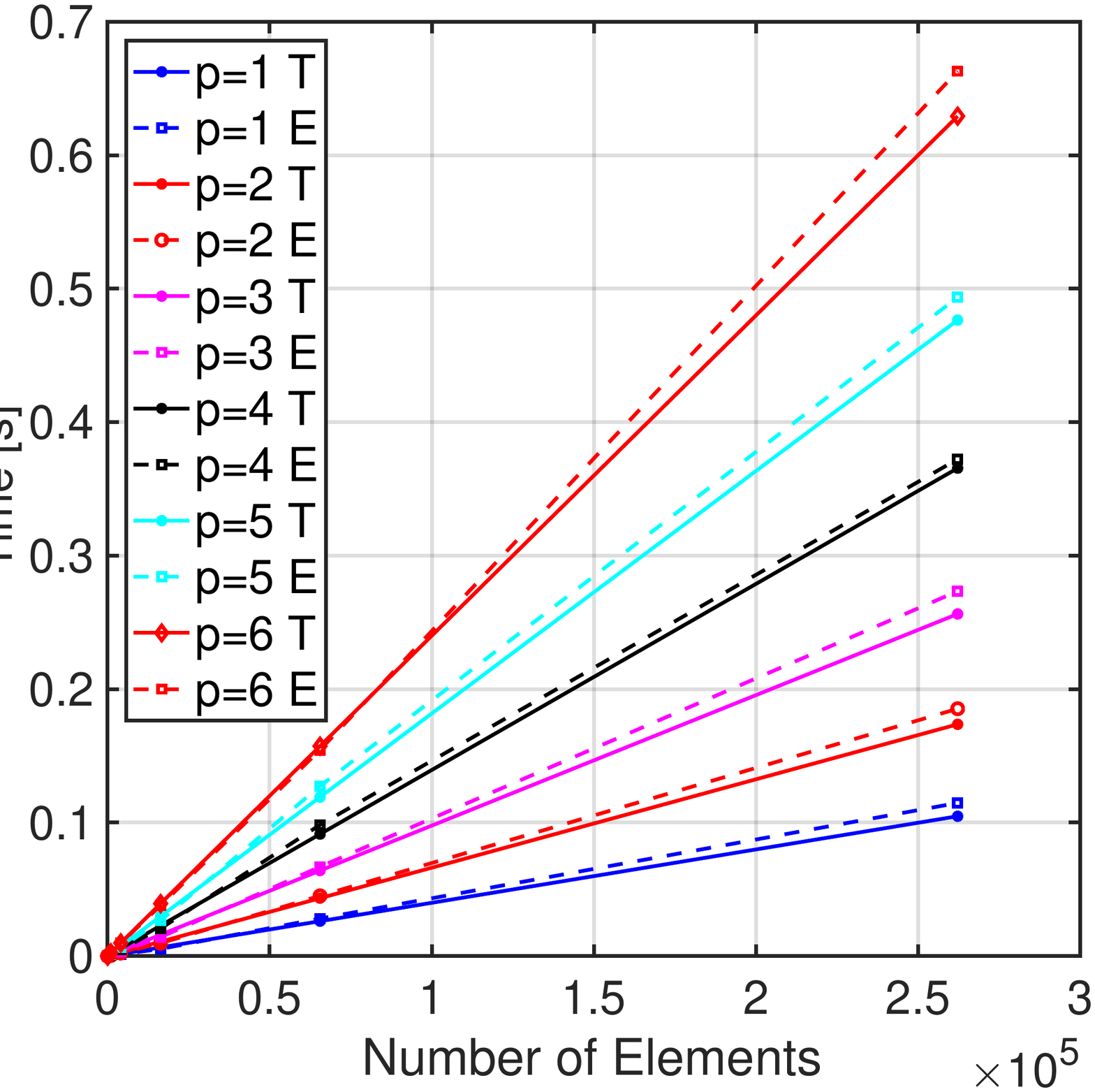}
  }
    \caption{Asymptotic and numerical estimates of back solve time complexity for EML,
      ML, and ND. Here, T stands for the theoretically
      estimated complexity derived in section \secref{complexity} and E for
      numerical experiment.}
  \figlab{theory_expt_solve}
\end{figure}

\begin{figure}[h!b!t!]
    \subfigure[ND]{
    \includegraphics[trim=0cm 3.1cm 0.25cm 3.5cm,clip=true,width=0.307\textwidth]{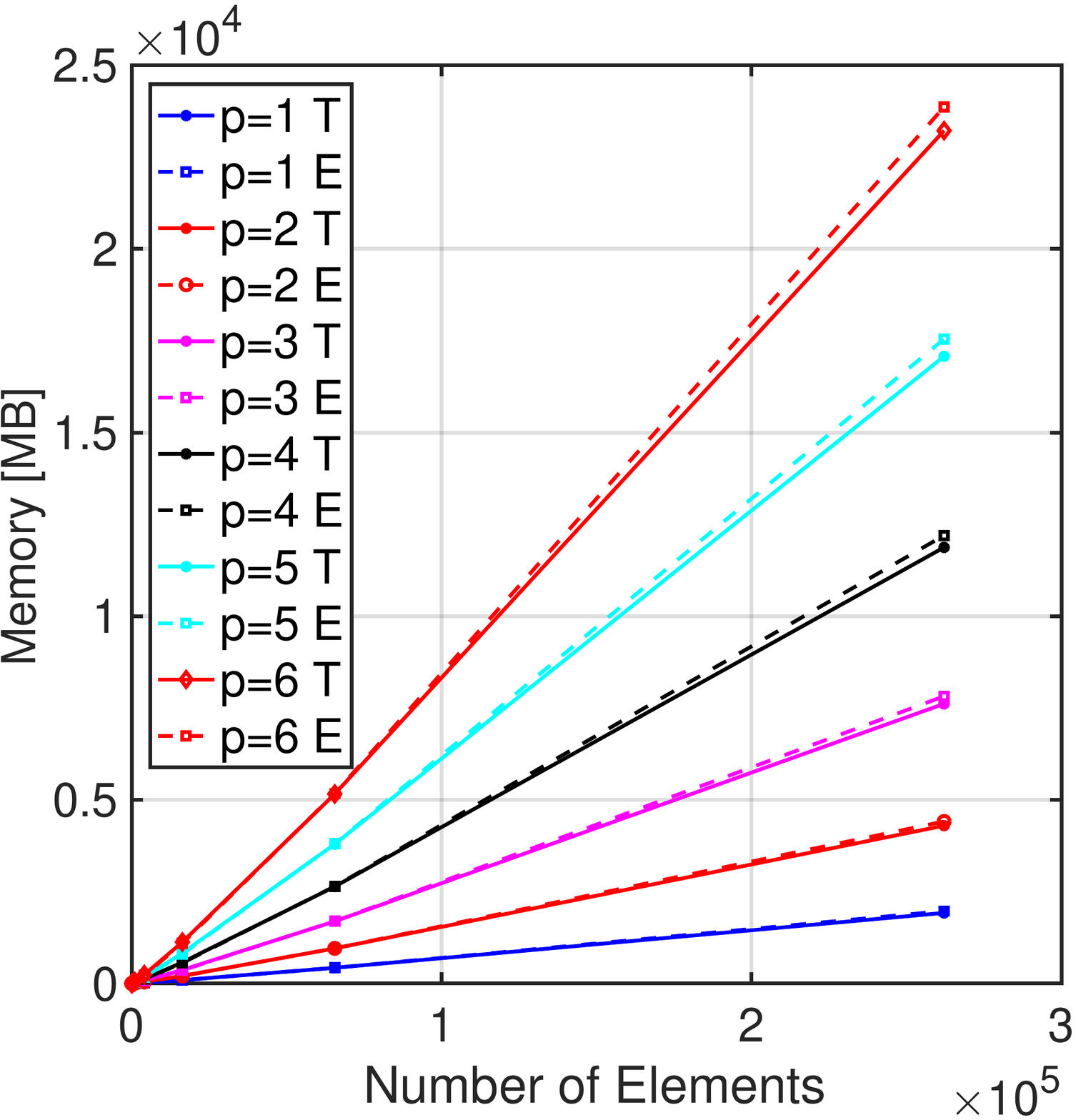}
  }
    \subfigure[ML]{
    \includegraphics[trim=0cm 2.9cm 0.25cm 3.5cm,clip=true,width=0.307\textwidth]{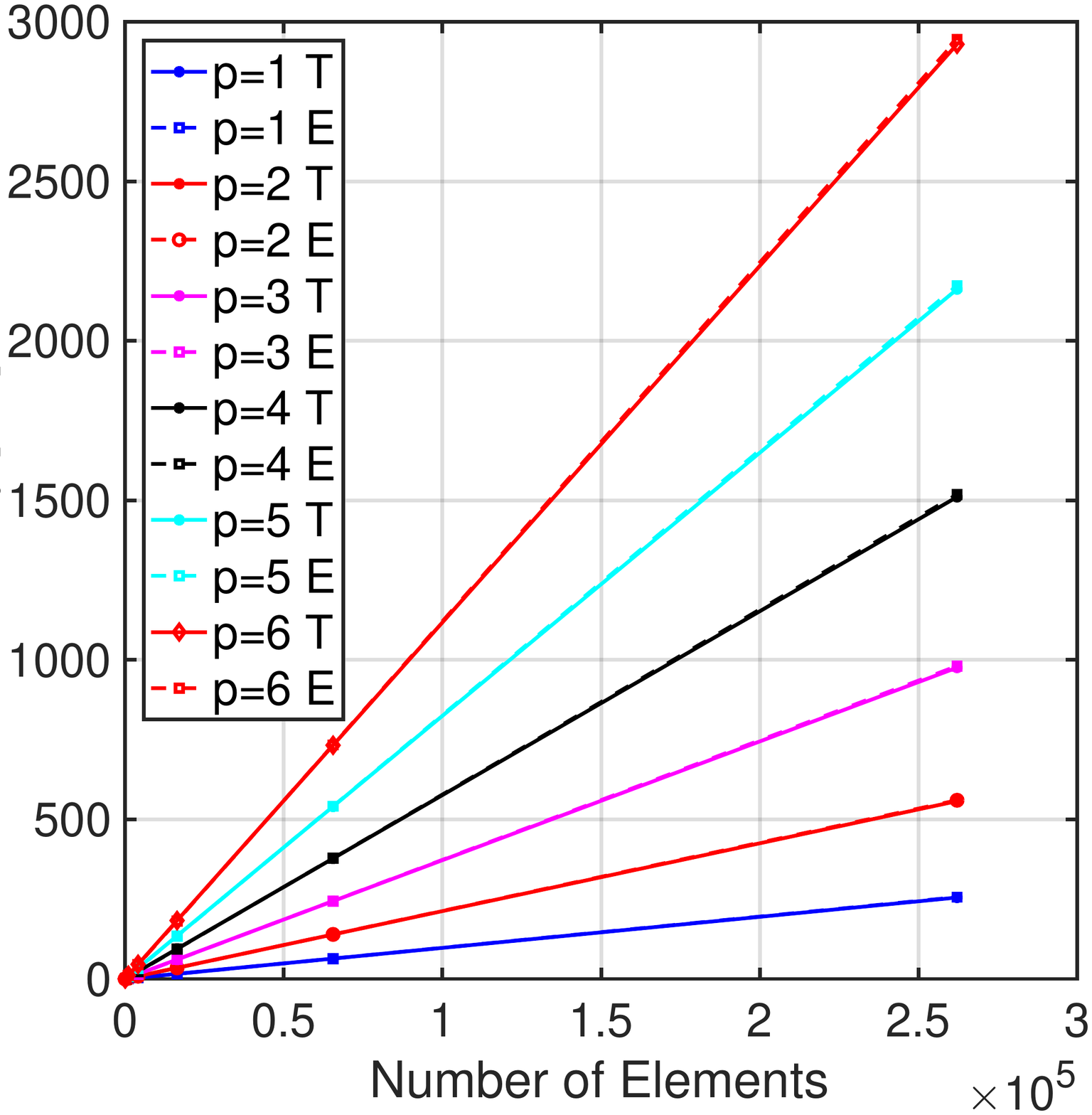}
  }
    \subfigure[EML]{
    \includegraphics[trim=0cm 3.1cm 0.25cm 3.5cm,clip=true,width=0.307\textwidth]{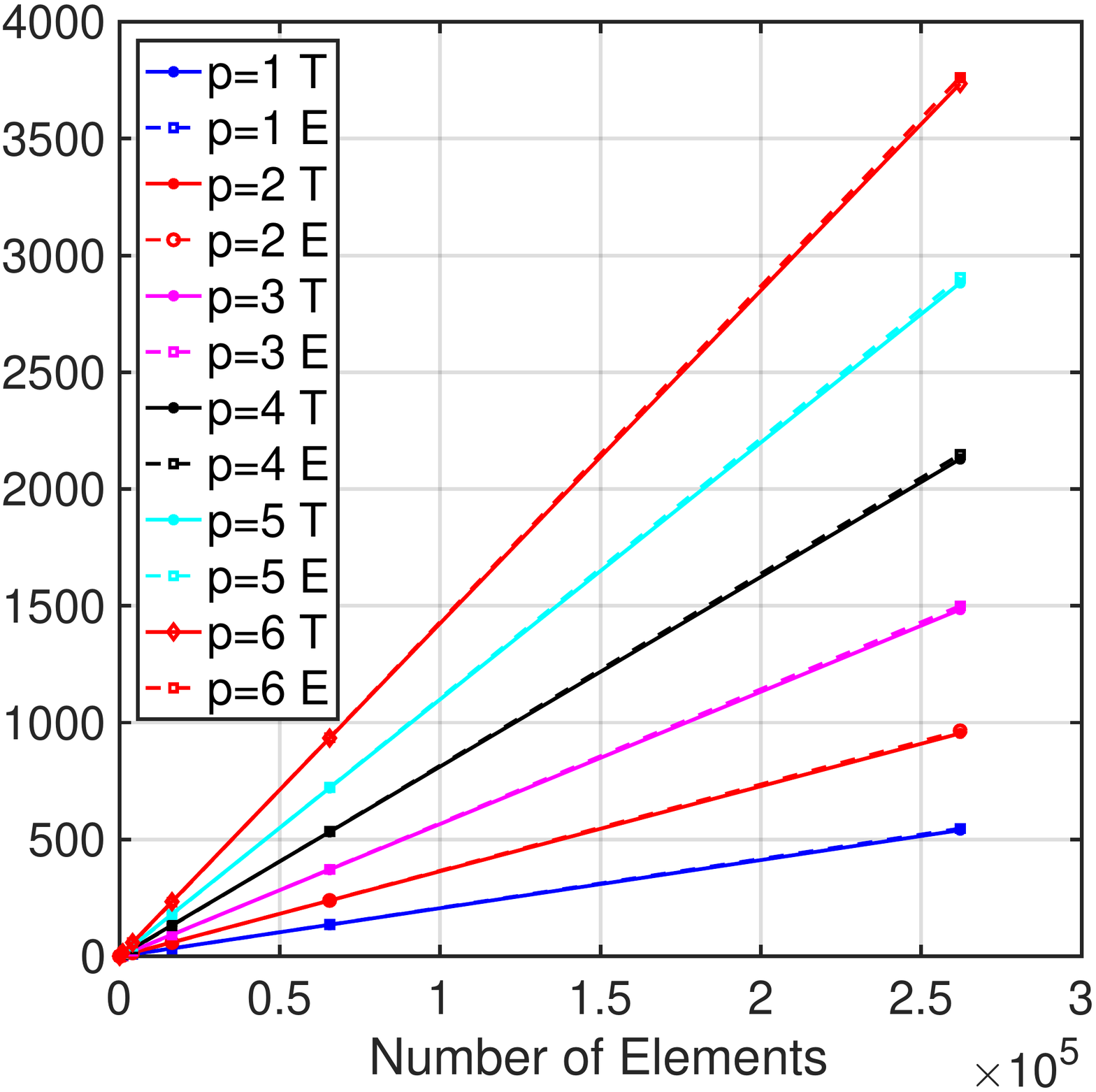}
  }
    \caption{Asymptotic and numerical estimates of memory complexity for EML,
      ML, and ND. Here, T stands for the 
    theoretically estimated complexity derived in section \secref{complexity} and E for numerical experiment.}
  \figlab{theory_expt_memory}
\end{figure}

\subsubsection{Example II: Discontinuous highly heterogeneous permeability}

In this section we test the robustness of the algorithm for elliptic PDE with a highly
discontinuous and heterogeneous permeability field. To that end, we
take $\betab=0$ and ${\bf K} = \kappa \mc{I}$ in
\eqref{model_problem}, where $\kappa$ is chosen according to example 2
in \cite{ho2016hierarchical} and is shown in Figure \figref{perm} for
three different meshes. The forcing and boundary condition in
\eqref{model_problem} are chosen as $f=1$ and $g_D=0$. This is a
difficult test case as the permeability varies by four orders of
magnitude and is also highly heterogeneous as seen in Figure
\figref{perm}.

\begin{figure}[h!b!t!]
    \subfigure[$N$ = 6]{
    \includegraphics[trim=3cm 7cm 2.5cm 7cm,clip=true,width=0.3\textwidth]{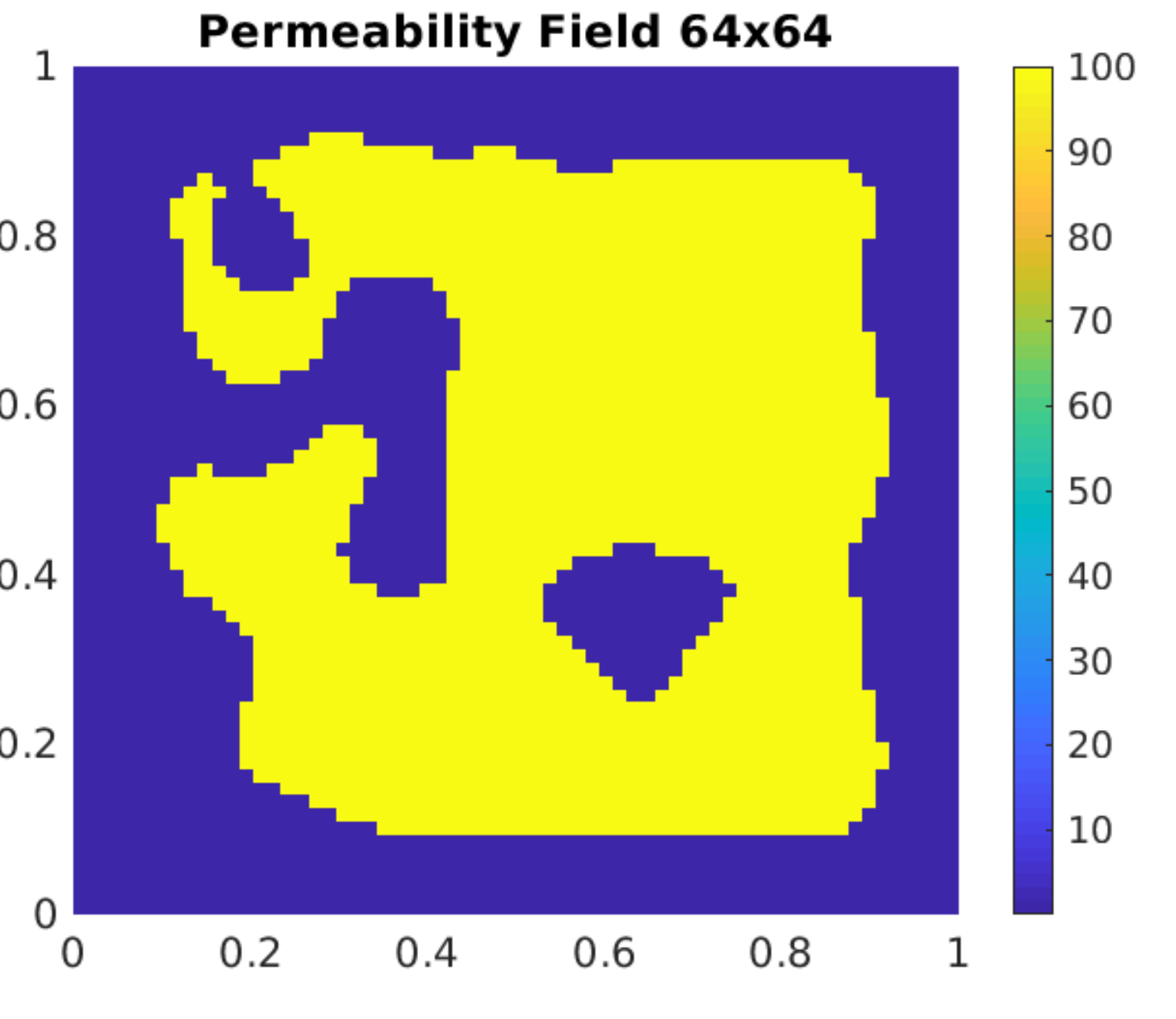}
    \figlab{perm1}
  }
    \subfigure[$N$ = 7]{
    \includegraphics[trim=3cm 7cm 2.5cm 7cm,clip=true,width=0.3\textwidth]{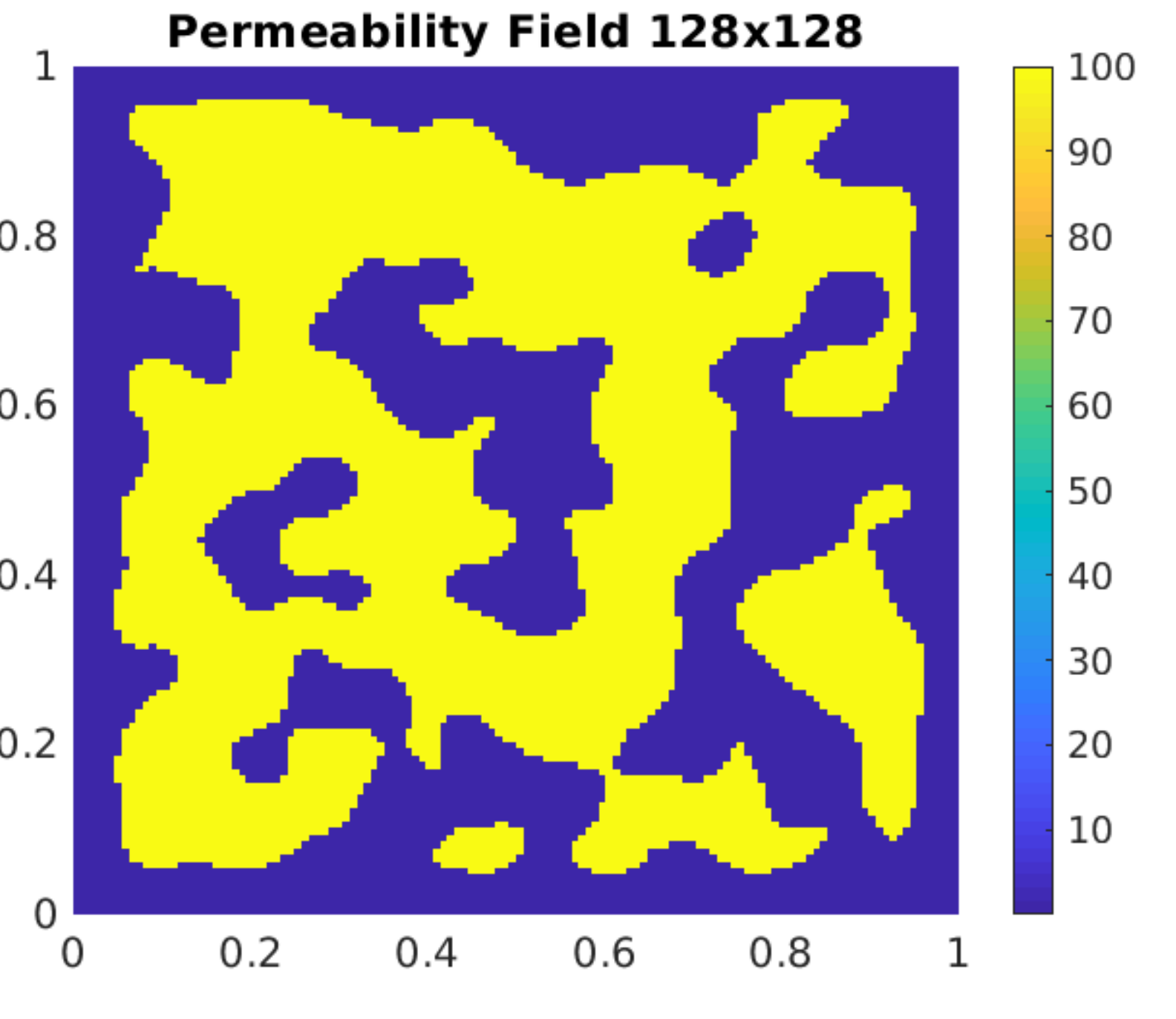}
    \figlab{perm2}
  }
    \subfigure[$N$ = 8]{
    \includegraphics[trim=3cm 7cm 2.5cm 7cm,clip=true,width=0.3\textwidth]{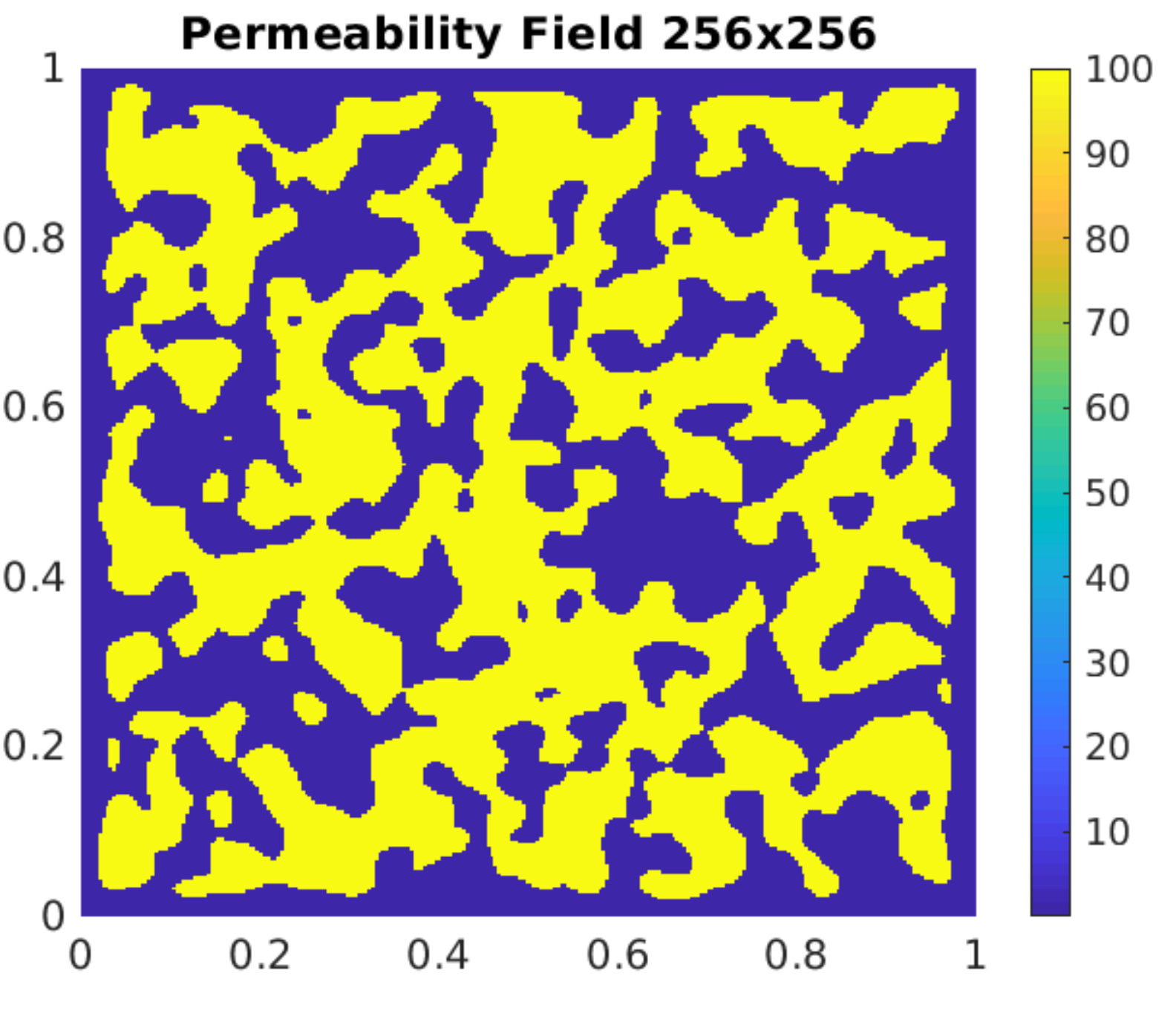}
    \figlab{perm3}
  }
    \caption{Discontinuous and heterogeneous permeability field
      \cite{ho2016hierarchical} on $64^2$, $128^2$ and
      $256^2$ meshes.}
  \figlab{perm}
\end{figure}

Tables \ref{tab:ML_exp2} and \ref{tab:EML_exp2} show the number of
iterations taken by ML- and EML-preconditioned GMRES, and in Table
\ref{tab:GMG_exp2} we compare them with those taken by GMRES preconditioned by one v-cycle of geometric
multigrid proposed in \cite{wildey2018unified}.
In Tables \ref{tab:ML_exp2}, \ref{tab:EML_exp2} and \ref{tab:GMG_exp2} the error, $\max{|\lambda_{direct}-\lambda|}$,
after 200 iterations is shown in the
parentheses. Here, $\lambda_{direct}$ and $\lambda$ denote trace
solution vectors obtained by the direct solver and the corresponding
iterative solver (ML, EML or geometric multigrid), respectively.
The results show that the geometric multigrid in
\cite{wildey2018unified} yields the least number of iterations or
more accurate approximation when convergence is not attained with
200 iterations.  This is expected since ML and EML algorithms have
smoothing only on the fine level while smoothing is performed on all
levels (in addition to the local smoothing) for the geometric multigrid
algorithm, and for elliptic-type PDEs performing smoothing on all
  levels typically provides better performance \cite{trottenberg2000multigrid}. However, the proposed algorithms in this paper
targets beyond elliptic PDEs, and for that reason it is not clear if
smoothing on coarser levels helps reduce the iteration counts
\cite{reusken1994multigrid}. We would also like to point out that the least number of iterations in geometric multigrid does
not translate directly to least overall time to solution which in turn depends on time per iteration and set-up cost for the three methods. In future work we will compare the overall time to solution for ML, EML and geometric multigrid. This challenging example clearly shows the benefits
of high-order nature of HDG, that is, for all three methods as the solution order increases
the solution  is not only more accurate but also
obtained with less number of GMRES iterations.


Between ML and EML, we can see that EML requires less number of
iterations and attains more accuracy at higher levels (see, e.g., the results with $N=8$ in Tables
\ref{tab:ML_exp2} and \ref{tab:EML_exp2}). The benefit of enrichment is clearly observed
for $\p=\LRc{3,4,5}$, in which EML is almost four orders of magnitude more
accurate than ML (see last row and columns $4-6$ of Tables \ref{tab:ML_exp2} and \ref{tab:EML_exp2}). 
For coarser meshes, the iteration counts of ML and
EML are similar. Finally, it is interesting to notice that 
columns $4-7$ in Table \ref{tab:EML_exp2}, corresponding to
$\p=\LRc{3,4,5,6}$, for EML (highlighted in blue) have similar
iteration counts and accuracy when compared to columns $3-6$ in Table
\ref{tab:GMG_exp2}, corresponding to $\p=\LRc{2,3,4,5}$, for the geometric
multigrid method.

\begin{table}[h!b!t!]
\centering
\begin{tabular}{|c|c|c|c|c|c|c|}
\hline
$N$    & \scriptsize $\p=1$ & \scriptsize $\p=2$ & \scriptsize $\p=3$ & \scriptsize $\p=4$ & \scriptsize $\p=5$ & \scriptsize $\p=6$ \\
\cline{1-7}

    6 & 178  & 137  & 107  & 92 & 79 & 73 \\
    7 & * ($10^{-5}$) & * ($10^{-8}$)& 167 & 138 & 113 & 99 \\
    8 & * ($10^{-2}$)& * ($10^{-2}$)& * ($10^{-2}$)& * ($10^{-4}$)& * ($10^{-5}$)& * ($10^{-7}$)\\
\hline
\end{tabular}
    \caption{\label{tab:ML_exp2} Example II: number of ML-preconditioned GMRES iterations as the mesh and solution order are refined. 
     }
\end{table}

\begin{table}[h!b!t!]
\centering
\begin{tabular}{|c|c|c|c|c|c|c|}
\hline
$N$    & \scriptsize $\p=1$ & \scriptsize $\p=2$ & \scriptsize $\p=3$ & \scriptsize $\p=4$ & \scriptsize $\p=5$ & \scriptsize $\p=6$ \\
\cline{1-7}

    6 & 180  & 132  & \myblue{115 } & \myblue{98 } & \myblue{81 } & \myblue{72} \\
    7 & * ($10^{-7}$) & 178 & \myblue{155}& \myblue{132}& \myblue{112}& \myblue{93}\\
    8 & * ($10^{-4}$)& * ($10^{-5}$)& \myblue{* ($10^{-6}$)}& \myblue{* ($10^{-8}$)}& \myblue{195}&\myblue{176}\\

\hline
\end{tabular}
\caption{\label{tab:EML_exp2} Example II: number of EML-preconditioned GMRES iterations as the mesh and solution order are refined. 
    }
\end{table}

\begin{table}[h!b!t!]
\centering
\begin{tabular}{|c|c|c|c|c|c|c|}
\hline
$N$    & \scriptsize $\p=1$ & \scriptsize $\p=2$ & \scriptsize $\p=3$ & \scriptsize $\p=4$ & \scriptsize $\p=5$ & \scriptsize $\p=6$ \\
\cline{1-7}

    6 & 157 & \myred{114} & \myred{97} & \myred{85} & \myred{74} & 69 \\
    7 & * ($10^{-8}$) & \myred{157}& \myred{129}& \myred{112}& \myred{98}& 88\\
    8 & * ($10^{-6}$)& \myred{* ($10^{-7}$)}& \myred{* ($10^{-8}$)}& \myred{190}& \myred{176}& 170\\
\hline
\end{tabular}
\caption{\label{tab:GMG_exp2} Example II: number of geometric-multigrid-preconditioned GMRES iterations as the mesh and solution order are refined. 
    }
\end{table}

\subsection{Example III: Transport equation}
\seclab{transport} In this section we apply ML and EML to a pure
transport equation. To that end, we take ${\bf K} = {\bf 0}$ in
\eqref{model_problem}.  Similar to \cite{MR2511736,iHDG}, we consider the velocity field
$\betab=(1+\sin(\pi y/2), 2)$, forcing $f=0$, and the inflow boundary
conditions
\[
g = 
\left\{
\begin{array}{ll}
1 & x = 0, 0 \le y \le 2 \\
\sin^6\LRp{\pi x} & 0< x \le 1, y = 0 \\
0 & 1 \le x \le 2, y = 0
\end{array}
\right.
.
\]
 The solution is shown in Figure \figref{Shock} and the difficulty of
 this test case comes from the presence of a curved discontinuity
 (shock) emanating from the inflow to the outflow boundaries. In
 Table \ref{tab:exp_3_precond} we show the iteration counts for both ML-
 and EML-preconditioned GMRES for different solution orders and mesh
 levels. As can be seen,  while $h-$scalability is not attained with both ML
 and EML, $p-$scalability is observed for both methods, i.e., the number of GMRES iterations
 is almost constant for all solution orders.
 Again EML takes less iteration counts than ML for all cases.


Table \ref{tab:exp_3_bJ_precond} shows the iteration counts for
block-Jacobi preconditioned GMRES. Compared to ML and EML in Table
\ref{tab:exp_3_precond}, the iteration counts for block-Jacobi are
higher, and for levels 7 and 8 block-Jacobi does not converge within the
maximum number of iteration counts. This indicates though both ML and EML do not give $h-$scalable
results, they provide a global coupling for the two-level algorithm
and thus help in reducing the total number of
iterations. Moreover, both the ML and EML algorithms are robust (with respect to
convergence) even for solution with shock. It is important to point out 
that for pure transport problems it is in general not trivial to obtain
$h-$scalable results unless some special smoother, which follows the direction of
  convection, is used \cite{bank1981comparison,hackbusch1997downwind,kim2004uniformly,wang1999crosswind}. For this reason, the moderate growth in the iteration counts for both ML and EML algorithms is encouraging.

  Next, we test the algorithms on a smooth exact solution (see Figure \figref{Smooth}) given by
  \[
      u^e=\frac{1}{\pi}\sin(\pi x)\cos(\pi y).
  \]
All the other parameters are the same as those for the discontinuous
solution considered above.  Tables \ref{tab:exp_3_smooth_precond} and
\ref{tab:exp_3_smooth_bJ_precond} show the number of ML-, EML-, and
block-Jacobi-preconditioned GMRES iterations. Table
\ref{tab:exp_3_smooth_precond} shows that the performance of ML-
and EML-preconditioned GMRES is similar to the one observed for
the elliptic equation with smooth solution in Table
\ref{tab:ML_EML_exp1}. Block-Jacobi preconditioned GMRES, on the
other hand, is more or less independent of the smoothness of the
solution as Table \ref{tab:exp_3_smooth_bJ_precond} for
the smooth solution is very similar to Table \ref{tab:exp_3_bJ_precond} for the discontinuous
solution. Thus this example demonstrates
that, unlike many standard iterative algorithms which depend on
the nature of the PDE under consideration, the performance of ML and EML algorithms seems
to depend only on the smoothness of the solution and otherwise is 
independent of the underlying PDE. This behavior is expected, again, thanks to their root in direct solver
strategy.

\begin{figure}[h!b!t!]
    \subfigure[Discontinuous solution]{
\includegraphics[trim=3.5cm 8.5cm 4cm 9cm,clip=true,width=0.48\textwidth,height=0.38\textwidth]{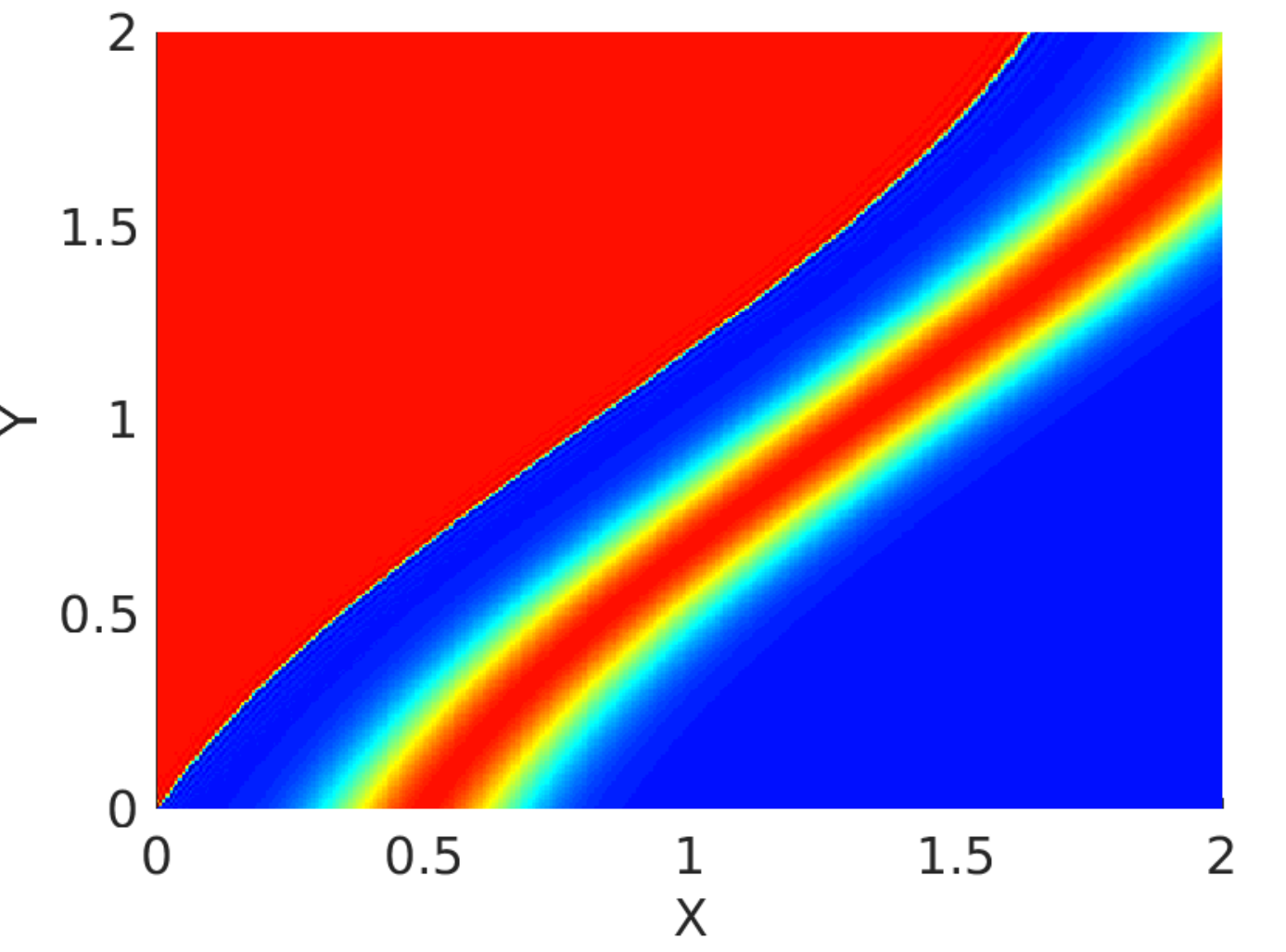}
    \figlab{Shock}
  }
    \subfigure[Smooth solution]{
\includegraphics[trim=1cm 6cm 2cm 6cm,clip=true,width=0.48\textwidth]{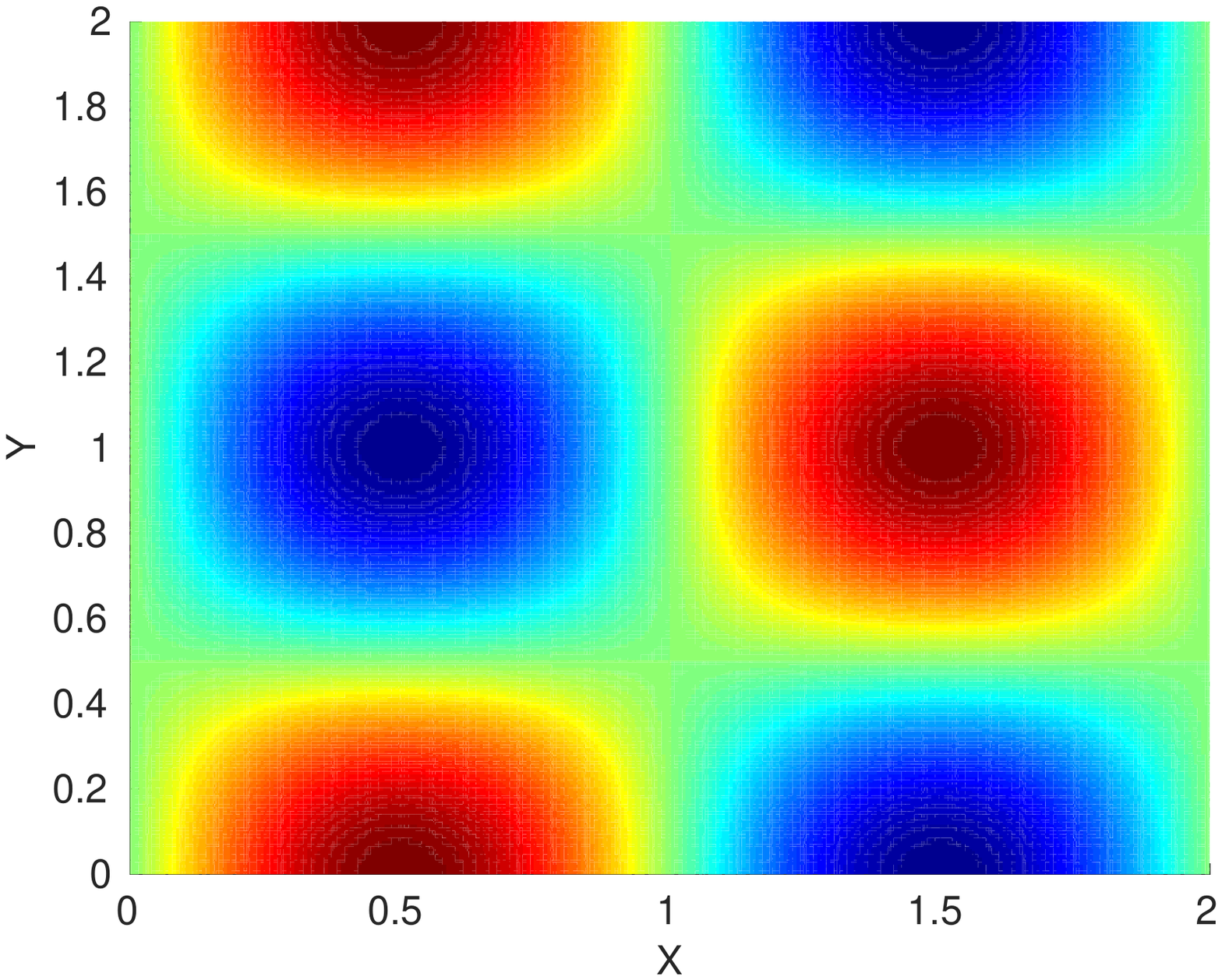}
    \figlab{Smooth}
  }
    \caption{Discontinuous and smooth solution for transport equation on a $64\times 64$ uniform mesh and $\p=6$ solution order.}
  \figlab{Shock_smooth}
\end{figure}


\begin{table}[h!b!t!]
\centering
\begin{tabular}{|c|c|c|c|c|c|c||c|c|c|c|c|c|}
\hline
& \multicolumn{6}{c||}{ML with GMRES} & \multicolumn{6}{c|}{EML
with GMRES}\\
\cline{2-13}
\!\!\! $N$ \!\!\!\! &  \multicolumn{6}{c||}{\!\!\scriptsize  $\p$\!\!} &  \multicolumn{6}{c|}{\!\!\scriptsize $\p$\!\!}\\
\hline
& \scriptsize 1 & \scriptsize 2 & \scriptsize 3 & \scriptsize 4 & \scriptsize 5 & \scriptsize 6 & \scriptsize 1 & \scriptsize 2 & \scriptsize 3 & \scriptsize 4 & \scriptsize 5 & \scriptsize 6 \\
\cline{2-13}
2 & 4  & 5  & 5  & 5  & 5  & 5  & 3  & 4  & 5  & 5  & 5  & 5\\
3 & 7  & 7  & 7  & 7  & 7  & 7  & 6  & 6  & 7  & 7  & 7  & 7\\
4 & 10 & 11 & 11 & 11 & 10 & 11 & 8  & 9  & 10 & 10 & 10 & 10\\
5 & 16 & 17 & 16 & 18 & 17 & 17 & 10 & 14 & 15 & 15 & 16 & 16\\
6 & 25 & 27 & 26 & 28 & 27 & 28 & 15 & 22 & 23 & 24 & 25 & 26\\
7 & 41 & 44 & 44 & 46 & 45 & 47 & 21 & 34 & 39 & 43 & 43 & 44\\
8 & 66 & 76 & 79 & 82 & 81 & 83 & 31 & 55 & 67 & 75 & 77 & 79\\
\hline
\end{tabular}
\caption{\label{tab:exp_3_precond} Example III. Discontinuous solution: number of  ML- and EML-preconditioned GMRES iterations.}
\end{table}

\begin{table}[h!b!t!]
\centering
\begin{tabular}{|c|c|c|c|c|c|c|}
\hline
$N$ & \scriptsize $\p=1$ & \scriptsize $\p=2$ & \scriptsize $\p=3$ & \scriptsize $\p=4$ & \scriptsize $\p=5$ & \scriptsize $\p=6$\\
\cline{1-7}
2 & 7   & 7  & 7  & 7   & 7   & 7  \\
3 & 9   & 9  & 10 & 10  & 10  & 9  \\
4 & 14  & 14 & 14 & 15  & 14  & 14 \\
5 & 24  & 24 & 24 & 25  & 25  & 25 \\
6 & 43  & 41 & 44 & 45  & 46  & 46 \\
7 & 78  & 76 & *  & *   & *   & *  \\
8 & 146 & *  & *  & *   & *   & *  \\
\hline
\end{tabular}
\caption{\label{tab:exp_3_bJ_precond} Example III. Discontinuous solution: number of block-Jacobi preconditioned GMRES iterations.}
\end{table}

\begin{table}[h!b!t!]
\centering
\begin{tabular}{|c|c|c|c|c|c|c||c|c|c|c|c|c|}
\hline
& \multicolumn{6}{c||}{ML with GMRES} & \multicolumn{6}{c|}{EML
with GMRES}\\
\cline{2-13}
\!\!\! $N$ \!\!\!\! &  \multicolumn{6}{c||}{\!\!\scriptsize  $\p$\!\!} &  \multicolumn{6}{c|}{\!\!\scriptsize $\p$\!\!}\\
\hline
& \scriptsize 1 & \scriptsize 2 & \scriptsize 3 & \scriptsize 4 & \scriptsize 5 & \scriptsize 6 & \scriptsize 1 & \scriptsize 2 & \scriptsize 3 & \scriptsize 4 & \scriptsize 5 & \scriptsize 6 \\
\cline{2-13}
2 & 5  & 4  & 5  & 4  & 4  & 3  & 3  & 4  & 4 & 4 & 3 & 3\\
3 & 6  & 6  & 6  & 5  & 5  & 3  & 5  & 5  & 5 & 4 & 3 & 2\\
4 & 10 & 9  & 9  & 8  & 7  & 5  & 7  & 7  & 6 & 4 & 2 & 1\\
5 & 15 & 14 & 13 & 13 & 10 & 8  & 8  & 9  & 7 & 2 & 1 & 0\\
6 & 23 & 22 & 20 & 20 & 17 & 13 & 10 & 11 & 2 & 1 & 0 & 0\\
7 & 36 & 37 & 37 & 34 & 30 & 17 & 12 & 7  & 1 & 0 & 0 & 0\\
8 & 58 & 63 & 65 & 62 & 48 & 21 & 12 & 2  & 0 & 0 & 0 & 0\\
\hline
\end{tabular}
\caption{\label{tab:exp_3_smooth_precond} Example III. Smooth solution: number of  ML- and EML-preconditioned GMRES iterations.}
\end{table}

\begin{table}[h!b!t!]
\centering
\begin{tabular}{|c|c|c|c|c|c|c|}
\hline
$N$ & \scriptsize $\p=1$ & \scriptsize $\p=2$ & \scriptsize $\p=3$ & \scriptsize $\p=4$ & \scriptsize $\p=5$ & \scriptsize $\p=6$\\
\cline{1-7}
2 & 6   & 7  & 7  & 7   & 7   & 7  \\
3 & 9   & 9  & 9  & 9   & 9   & 9  \\
4 & 14  & 13 & 14 & 14  & 14  & 14 \\
5 & 23  & 22 & 23 & 24  & 24  & 24 \\
6 & 41  & 39 & 43 & 44  & 45  & 68 \\
7 & 76  & 74 & *  & *   & *   & *  \\
8 & 142 & *  & *  & *   & *   & *  \\
\hline
\end{tabular}
\caption{\label{tab:exp_3_smooth_bJ_precond} Example III. Smooth solution: number of block-Jacobi preconditioned GMRES iterations.}
\end{table}
\subsection{Convection-diffusion equation}
In this section we test the proposed algorithms for the convection-diffusion
equation in both diffusion- and convection-dominated regimes. To that
end, we consider $f=0$ in \eqref{model_problem}. We shall take some standard
 problems that are often  used to test the robustness of  multigrid algorithms
for convection-diffusion equations.

\subsubsection{Example IV}
\begin{figure}[h!b!t!]
\subfigure[$\alpha = 10$]{
\includegraphics[trim=3.5cm 8cm 4cm 9.15cm,clip=true,width=0.48\textwidth]{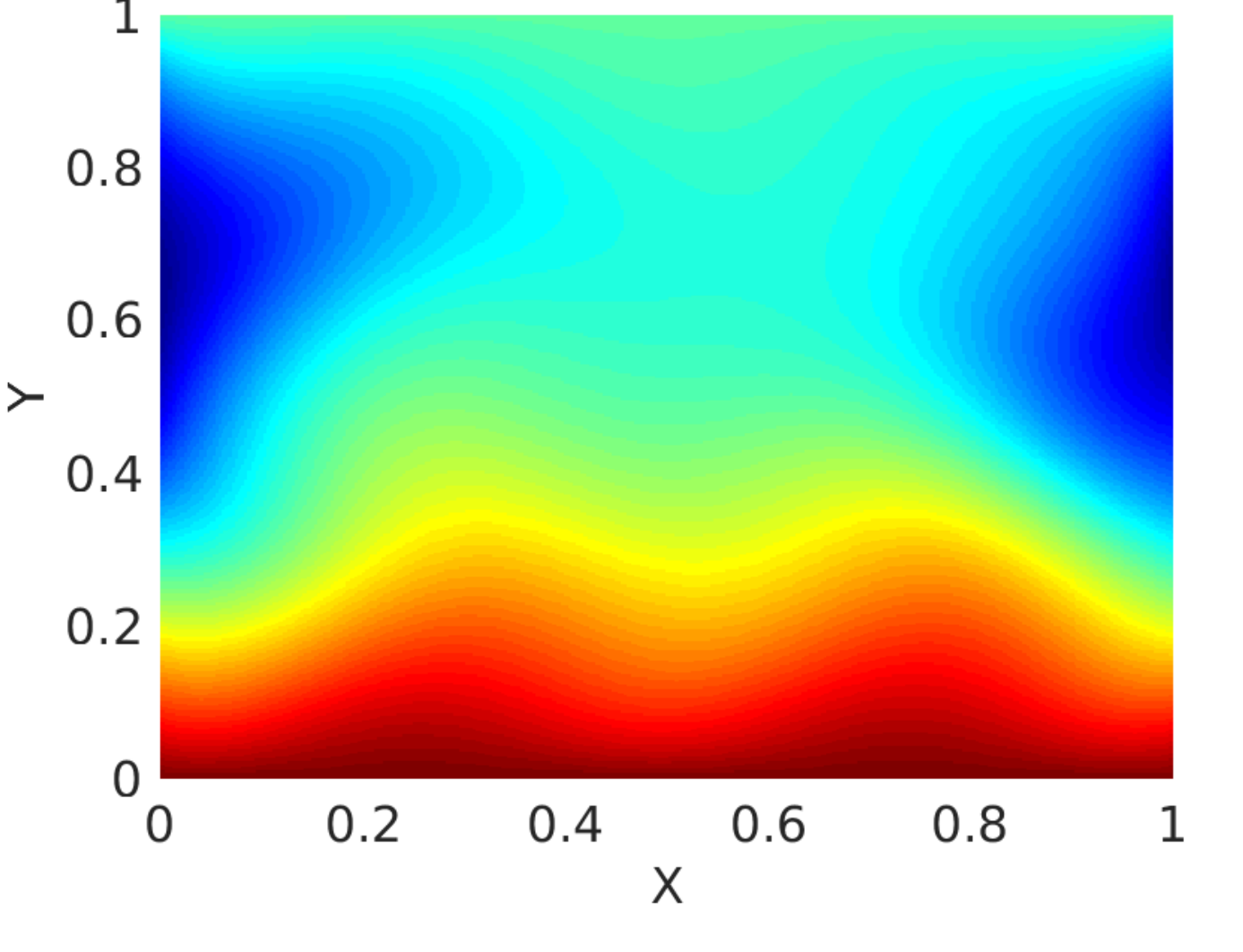}
    \label{a10}
}
\subfigure[$\alpha = 10^2$]{
\includegraphics[trim=3.5cm 8cm 4cm 9.15cm,clip=true,width=0.48\textwidth]{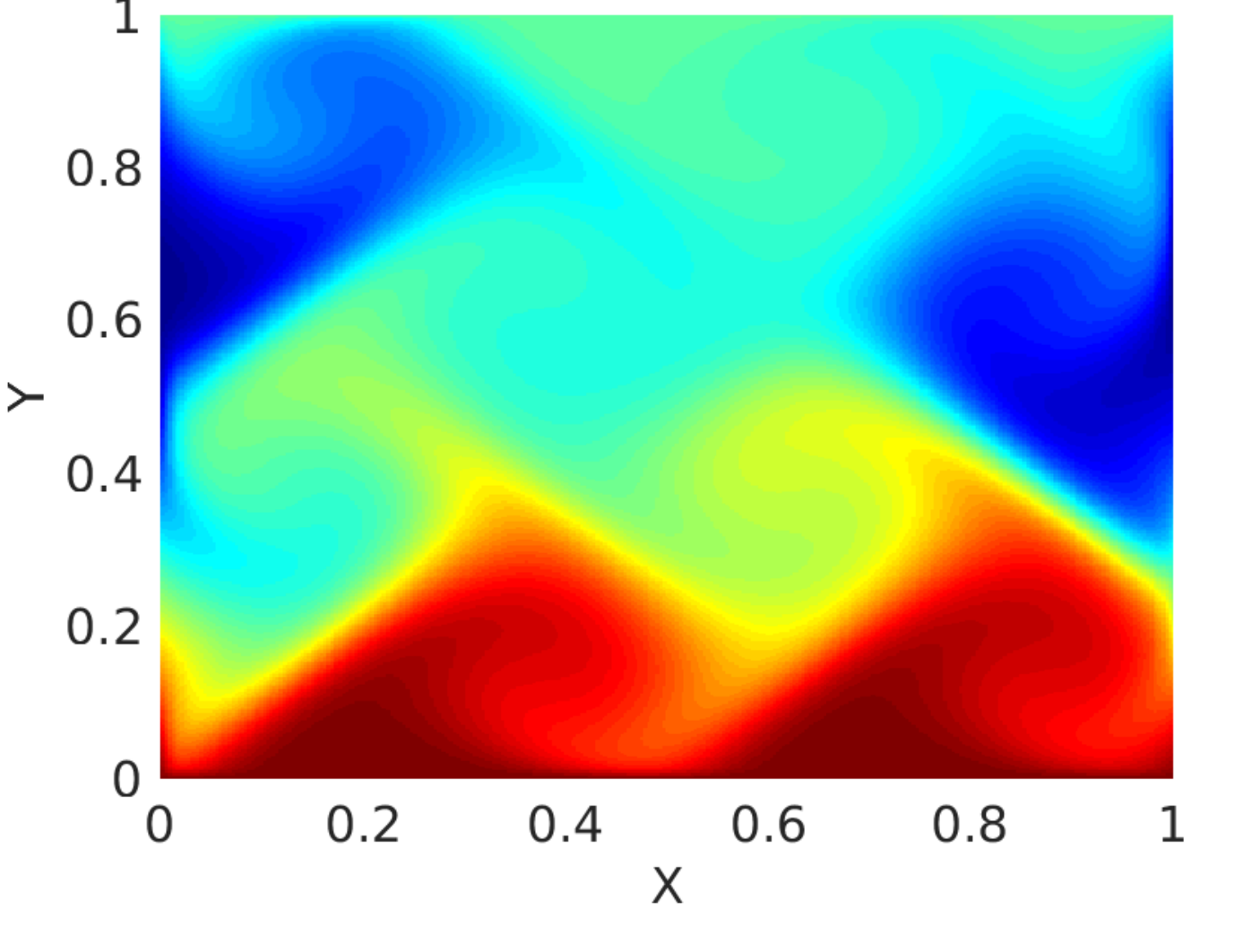}
    \label{a100}
}
\subfigure[$\alpha = 10^3$]{
\includegraphics[trim=3.5cm 8cm 4cm 9.15cm,clip=true,width=0.48\textwidth]{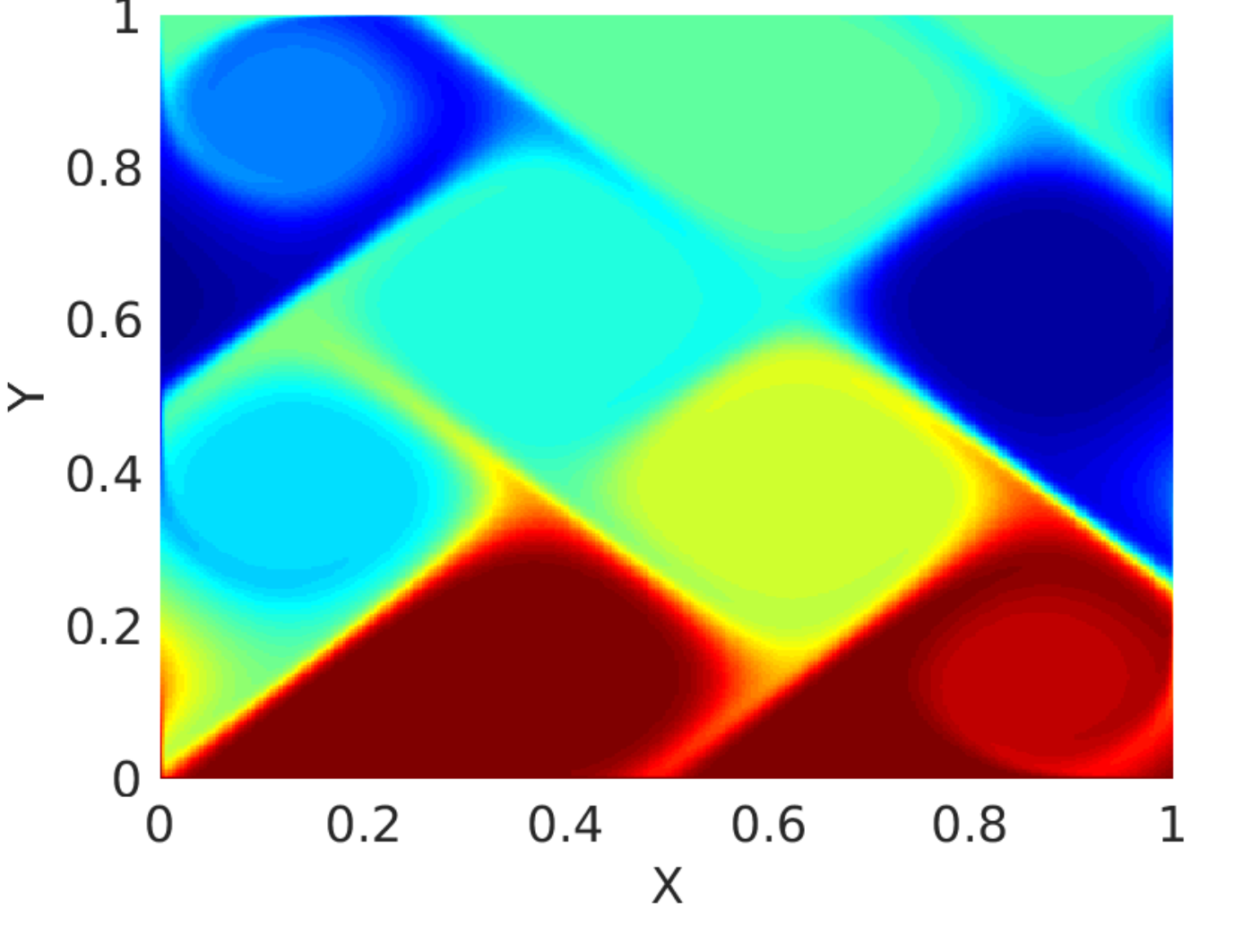}
    \label{a1000}
}
\subfigure[$\alpha = 10^4$]{
\includegraphics[trim=3.5cm 8cm 4cm 9.15cm,clip=true,width=0.48\textwidth]{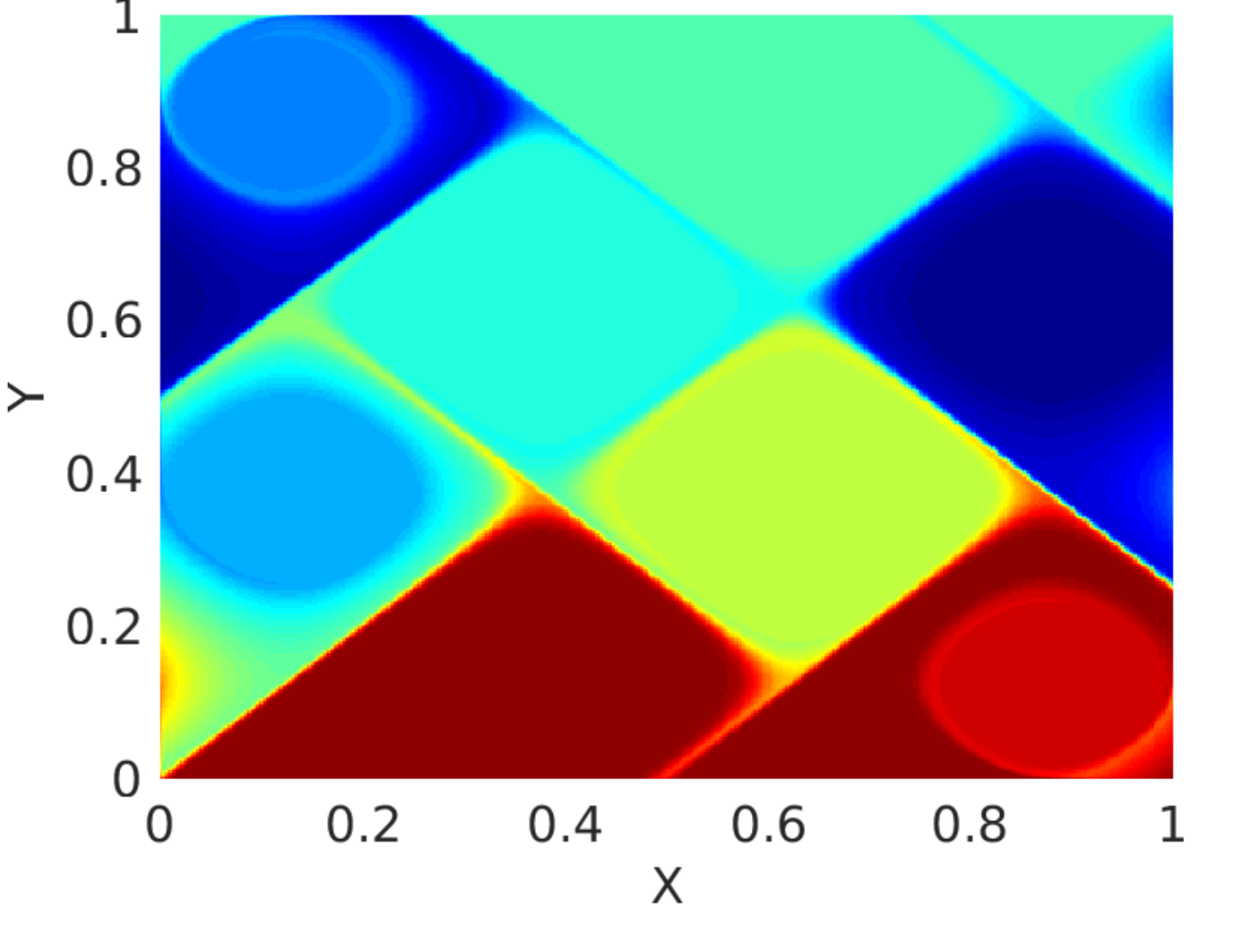}
    \label{a10000}
}

\caption{Example IV: solutions of the convection-diffusion equation for different values of $\alpha$ on a $64\times 64$ uniform mesh and $\p=6$ solution order.}
\label{exp_4}
\end{figure}

Here we consider an example similar to the one in
\cite{gillman2014direct}.  In particular, we take ${\bf K}=\mc{I}$,
$g_D=\cos(2y)(1-2y)$ and $\betab=(-\alpha\cos(4\pi y),-\alpha\cos(4\pi
x))$ in \eqref{model_problem}, where $\alpha$ is a parameter which
determines the magnitude of convection velocity. In Figure
\ref{exp_4}, solutions for different values of $\alpha$ ranging in
$\LRs{10,10^4}$ are shown. As $\alpha$ increases, the problem becomes
more convection-dominated and shock-like structures are formed.

In Tables \ref{tab:exp_4_a10}-\ref{tab:exp_4_a10000} are the
iteration counts for ML- and EML-preconditioned GMRES with various values of $\alpha$. We
observe the following. In all cases, as expected, the
iteration counts for  EML are less than for  ML. As the mesh is refined
we see growth in iterations for both ML and EML, though it is less for EML
than for  ML. With increase in solution order the iterations remain
(almost) constant, and in many cases decrease.
For mildly-to-moderately convection-dominated, i.e. $\alpha \in \LRs{10,10^3}$,
both ML and EML are robust in the sense that their
iteration counts negligibly vary with respect to
$\alpha$.
For $\alpha=10^4$, i.e. strongly convection-dominated,  we see an increase in iteration counts 
for both algorithms, though the growth is much less pronounced for EML
than for  ML (especially with low solution orders  $\p=\LRc{1,2,3,4}$). 

\begin{table}[h!b!t!]
\centering
\begin{tabular}{|c|c|c|c|c|c|c||c|c|c|c|c|c|}
\hline
& \multicolumn{6}{c||}{ML with GMRES} & \multicolumn{6}{c|}{EML
with GMRES}\\
\cline{2-13}
\!\!\! $N$ \!\!\!\! &  \multicolumn{6}{c||}{\!\!\scriptsize  $\p$\!\!} &  \multicolumn{6}{c|}{\!\!\scriptsize $\p$\!\!}\\
\hline
& \scriptsize 1 & \scriptsize 2 & \scriptsize 3 & \scriptsize 4 & \scriptsize 5 & \scriptsize 6 & \scriptsize 1 & \scriptsize 2 & \scriptsize 3 & \scriptsize 4 & \scriptsize 5 & \scriptsize 6 \\
\cline{2-13}
6 & 29 & 24 & 24 & 20 & 20 & 17 & 13 & 12 & 12 & 12 & 11 & 11\\
7 & 42 & 35 & 34 & 29 & 28 & 24 & 17 & 15 & 16 & 14 & 15 & 15\\
8 & 60 & 49 & 49 & 40 & 39 & 33 & 22 & 20 & 19 & 20 & 21 & 21\\
\hline
\end{tabular}
\caption{\label{tab:exp_4_a10} Example IV. $\alpha=10$: number of iterations for  ML- and EML-preconditioned GMRES.}
\end{table}

\begin{table}[h!b!t!]
\centering
\begin{tabular}{|c|c|c|c|c|c|c||c|c|c|c|c|c|}
\hline
& \multicolumn{6}{c||}{ML with GMRES} & \multicolumn{6}{c|}{EML
with GMRES}\\
\cline{2-13}
\!\!\! $N$ \!\!\!\! &  \multicolumn{6}{c||}{\!\!\scriptsize  $\p$\!\!} &  \multicolumn{6}{c|}{\!\!\scriptsize $\p$\!\!}\\
\hline
& \scriptsize 1 & \scriptsize 2 & \scriptsize 3 & \scriptsize 4 & \scriptsize 5 & \scriptsize 6 & \scriptsize 1 & \scriptsize 2 & \scriptsize 3 & \scriptsize 4 & \scriptsize 5 & \scriptsize 6 \\
\cline{2-13}
6 & 27 & 25 & 23 & 22 & 22 & 20 & 15 & 15 & 16 & 16 & 16 & 17\\
7 & 39 & 36 & 34 & 32 & 32 & 30 & 22 & 21 & 22 & 22 & 23 & 23\\
8 & 55 & 51 & 49 & 46 & 46 & 42 & 29 & 29 & 31 & 31 & 33 & 34\\
\hline
\end{tabular}
\caption{\label{tab:exp_4_a100} Example IV. $\alpha=10^2$: number of iterations for  ML- and EML-preconditioned GMRES.}
\end{table}

\begin{table}[h!b!t!]
\centering
\begin{tabular}{|c|c|c|c|c|c|c||c|c|c|c|c|c|}
\hline
& \multicolumn{6}{c||}{ML with GMRES} & \multicolumn{6}{c|}{EML
with GMRES}\\
\cline{2-13}
\!\!\! $N$ \!\!\!\! &  \multicolumn{6}{c||}{\!\!\scriptsize  $\p$\!\!} &  \multicolumn{6}{c|}{\!\!\scriptsize $\p$\!\!}\\
\hline
& \scriptsize 1 & \scriptsize 2 & \scriptsize 3 & \scriptsize 4 & \scriptsize 5 & \scriptsize 6 & \scriptsize 1 & \scriptsize 2 & \scriptsize 3 & \scriptsize 4 & \scriptsize 5 & \scriptsize 6 \\
\cline{2-13}
6 & 39 & 29 & 24 & 18 & 16 & 15 & 14 & 13 & 12 & 12 & 11 & 12\\
7 & 49 & 34 & 26 & 21 & 21 & 20 & 15 & 15 & 15 & 15 & 15 & 16\\
8 & 62 & 41 & 35 & 28 & 29 & 27 & 22 & 22 & 22 & 22 & 23 & 23\\
\hline
\end{tabular}
\caption{\label{tab:exp_4_a1000} Example IV. $\alpha=10^3$: number of iterations for  ML- and EML-preconditioned GMRES.}
\end{table}

\begin{table}[h!b!t!]
\centering
\begin{tabular}{|c|c|c|c|c|c|c||c|c|c|c|c|c|}
\hline
& \multicolumn{6}{c||}{ML with GMRES} & \multicolumn{6}{c|}{EML
with GMRES}\\
\cline{2-13}
\!\!\! $N$ \!\!\!\! &  \multicolumn{6}{c||}{\!\!\scriptsize  $\p$\!\!} &  \multicolumn{6}{c|}{\!\!\scriptsize $\p$\!\!}\\
\hline
& \scriptsize 1 & \scriptsize 2 & \scriptsize 3 & \scriptsize 4 & \scriptsize 5 & \scriptsize 6 & \scriptsize 1 & \scriptsize 2 & \scriptsize 3 & \scriptsize 4 & \scriptsize 5 & \scriptsize 6 \\
\cline{2-13}
6 & 89 & 101 & 88 & 68 & 59 & 51 & 29 & 43 & 44 & 41 & 37 & 33\\
7 & 136 & 133 & 98 & 75 & 65 & 55 & 40 & 48 & 44 & 38 & 35 & 34\\
8 & 192 & 141 & 101 & 76 & 67 & 56 & 46 & 45 & 39 & 35 & 34 & 34\\
\hline
\end{tabular}
\caption{\label{tab:exp_4_a10000} Example IV. $\alpha=10^4$: number of iterations for  ML- and EML-preconditioned GMRES.}
\end{table}

\subsubsection{Example V}
\begin{figure}[h!b!t!]
    \subfigure[$\kappa = 10^{-1}$]{
\includegraphics[trim=3.5cm 8cm 4cm 9.15cm,clip=true,width=0.48\textwidth]{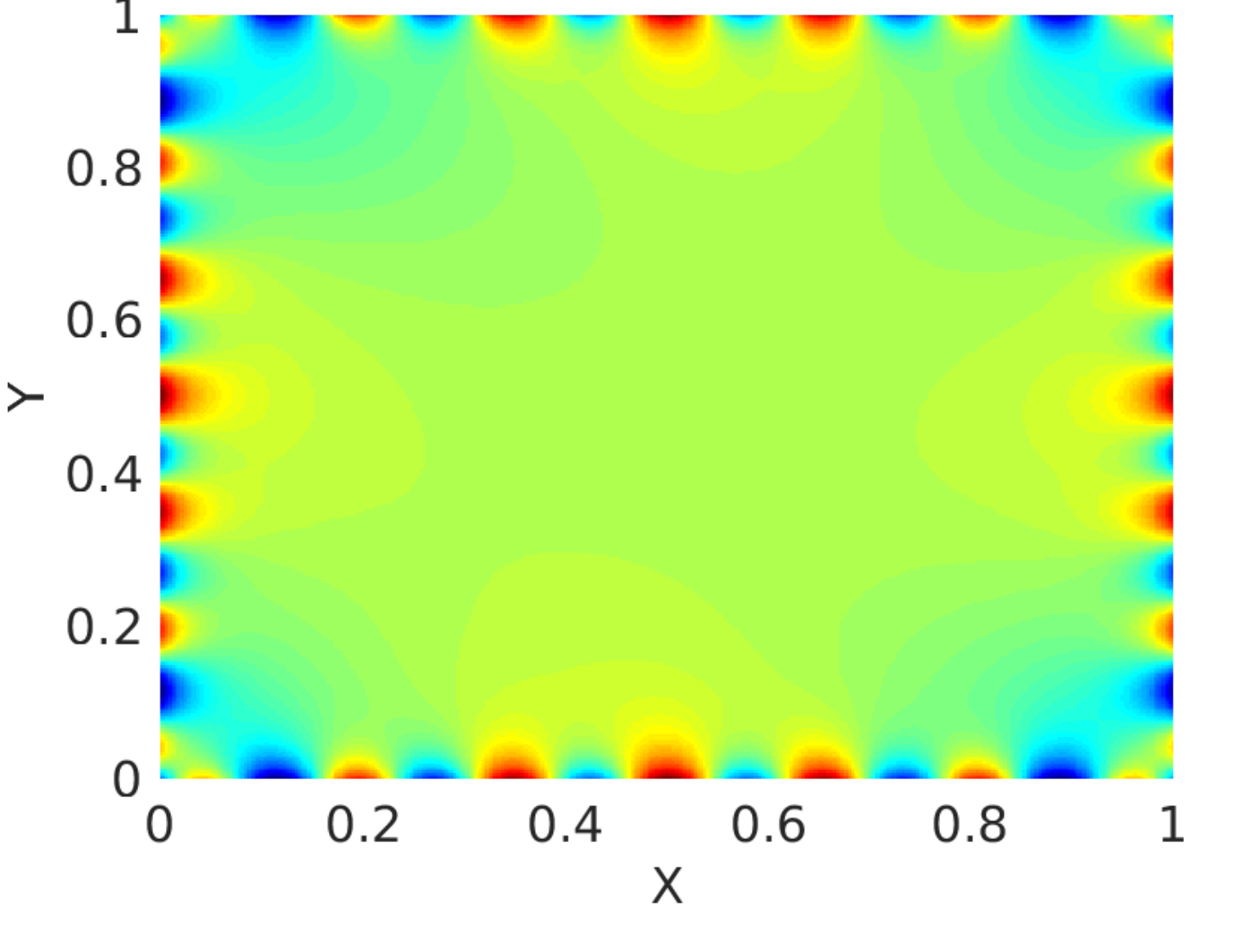}
    \label{k10}
}
    \subfigure[$\kappa = 10^{-2}$]{
\includegraphics[trim=3.5cm 8cm 4cm 9.23cm,clip=true,width=0.48\textwidth]{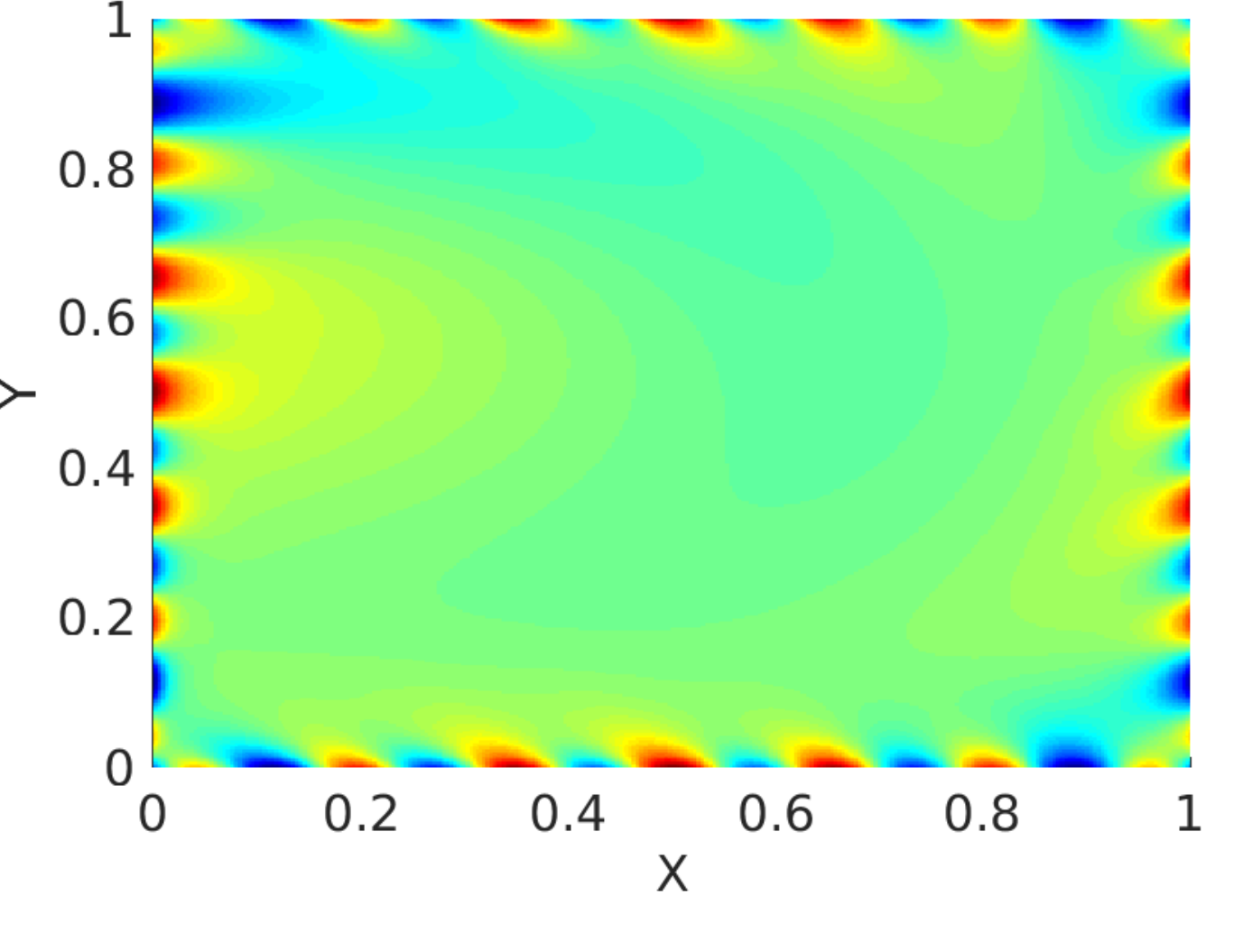}
    \label{k100}
}
    \subfigure[$\kappa = 10^{-3}$]{
\includegraphics[trim=3.5cm 8cm 4cm 9.15cm,clip=true,width=0.48\textwidth]{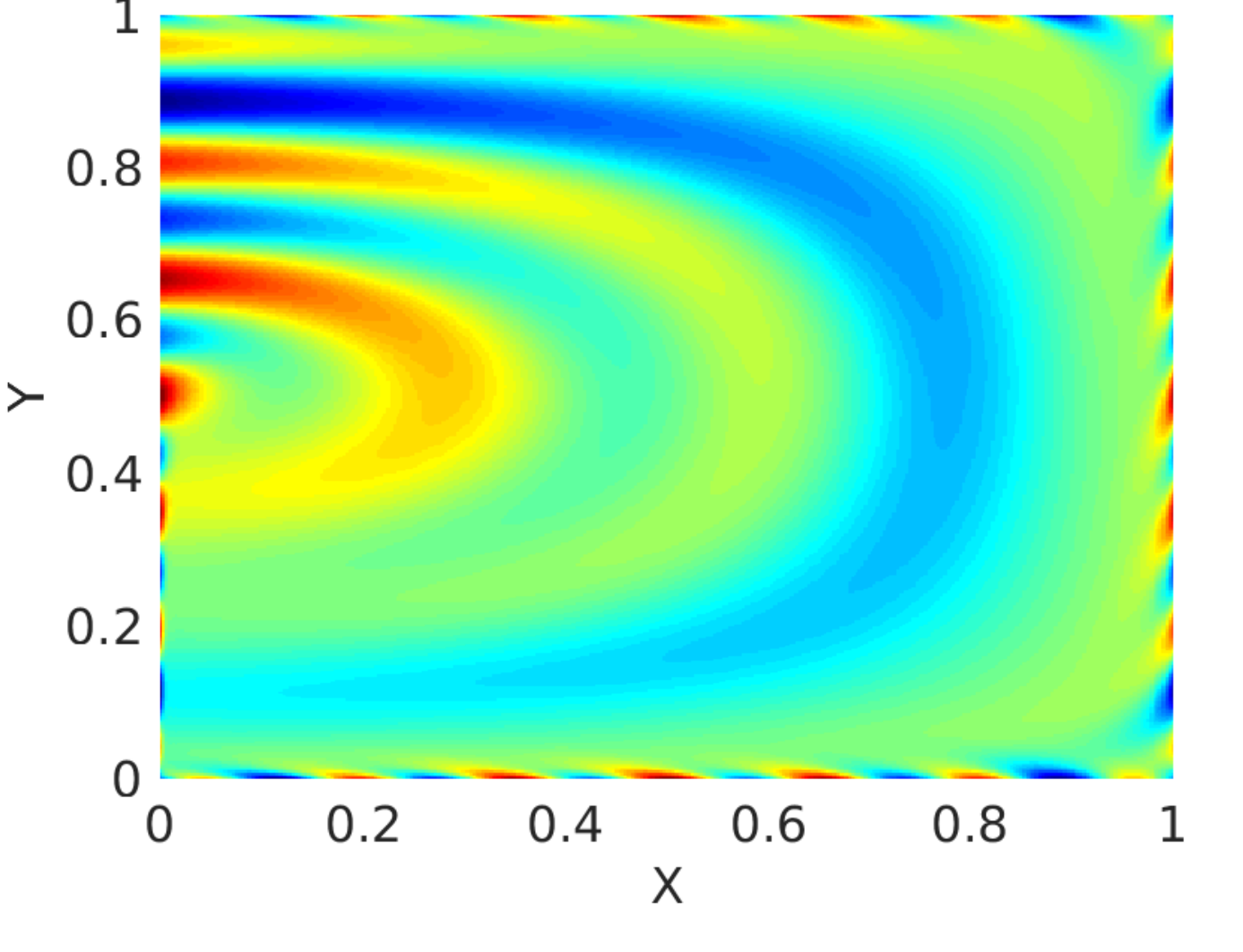}
    \label{k1000}
}
    \subfigure[$\kappa = 10^{-4}$]{
\includegraphics[trim=3.5cm 8cm 4cm 9.15cm,clip=true,width=0.48\textwidth]{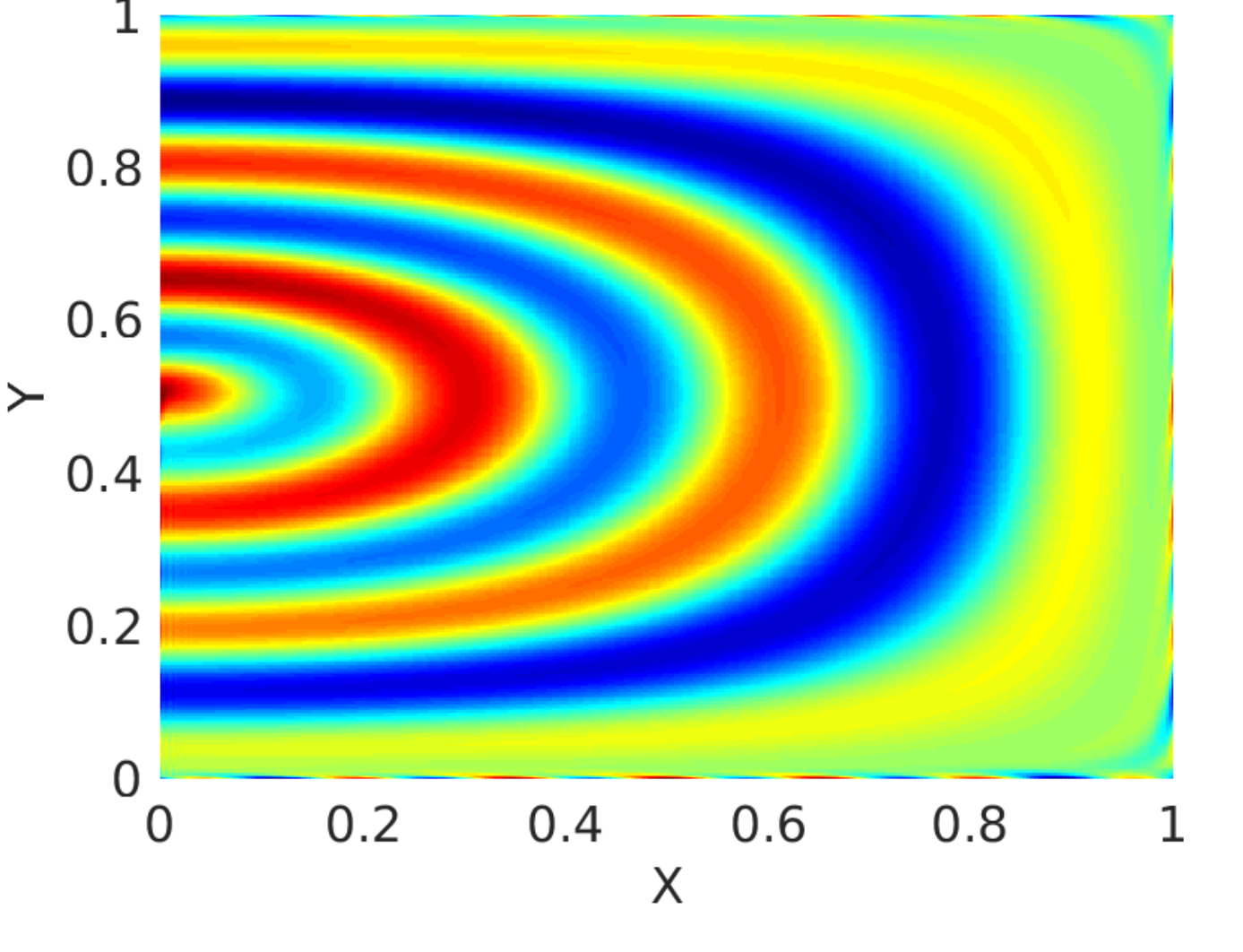}
    \label{k10000}
}

    \caption{Example V: solutions of the convection-diffusion equation for different values of $\kappa$ on a $64\times 64$ uniform mesh and $\p=6$ solution order.}
\label{exp_5}
\end{figure}

In this section a test case for multigrid method in
\cite{reusken1994multigrid} is considered. The parameters for this
example are $\betab=((2y-1)(1-x^2), 2xy(y-1))$, $g_D=\sin(\pi
x)+\sin(13\pi x)+\sin(\pi y)+\sin(13\pi y)$ and ${\bf K}=\kappa
\mc{I}$. The solution fields for various values of $\kappa$ in the
range $\LRs{10^{-1},10^{-4}}$ are shown in Figure \ref{exp_5}. As can
be observed, the problem becomes more convection dominated when the
diffusion coefficient $\kappa$ decreases. Tables
\ref{tab:exp_5_k10}-\ref{tab:exp_5_k10000} present the iteration
counts of ML- and EML-preconditioned GMRES for different values of $\kappa$.  Again, the
iteration counts for EML are less than for ML. As the mesh is refined
we see growth in iterations for both ML and EML, though it is less for
EML than for ML. An outlier is the case of $\kappa = 10^{-4}$ in
Table \ref{tab:exp_5_k10000}, where the number of EML-preconditioned
GMRES iterations reduces as the mesh is refined and ML-preconditioned
GMRES does not converge for
$p=2$ on mesh levels $6$ and $7$. By the time of writing, we have not yet found the reason for this behavior. 

\begin{table}[h!b!t!]
\centering
\begin{tabular}{|c|c|c|c|c|c|c||c|c|c|c|c|c|}
\hline
& \multicolumn{6}{c||}{ML with GMRES} & \multicolumn{6}{c|}{EML
with GMRES}\\
\cline{2-13}
\!\!\! $N$ \!\!\!\! &  \multicolumn{6}{c||}{\!\!\scriptsize  $\p$\!\!} &  \multicolumn{6}{c|}{\!\!\scriptsize $\p$\!\!}\\
\hline
& \scriptsize 1 & \scriptsize 2 & \scriptsize 3 & \scriptsize 4 & \scriptsize 5 & \scriptsize 6 & \scriptsize 1 & \scriptsize 2 & \scriptsize 3 & \scriptsize 4 & \scriptsize 5 & \scriptsize 6 \\
\cline{2-13}
6 & 28 & 23 & 23 & 21 & 21 & 19 & 13 & 13 & 13 & 13 & 13 & 13\\
7 & 40 & 34 & 34 & 30 & 29 & 26 & 18 & 17 & 18 & 16 & 18 & 18\\
8 & 58 & 47 & 48 & 43 & 41 & 37 & 23 & 22 & 23 & 23 & 25 & 25\\
\hline
\end{tabular}
\caption{\label{tab:exp_5_k10} Example V. $\kappa=10^{-1}$:
number of iterations for  ML- and EML-preconditioned GMRES.}
\end{table}

\begin{table}[h!b!t!]
\centering
\begin{tabular}{|c|c|c|c|c|c|c||c|c|c|c|c|c|}
\hline
& \multicolumn{6}{c||}{ML with GMRES} & \multicolumn{6}{c|}{EML
with GMRES}\\
\cline{2-13}
\!\!\! $N$ \!\!\!\! &  \multicolumn{6}{c||}{\!\!\scriptsize  $\p$\!\!} &  \multicolumn{6}{c|}{\!\!\scriptsize $\p$\!\!}\\
\hline
& \scriptsize 1 & \scriptsize 2 & \scriptsize 3 & \scriptsize 4 & \scriptsize 5 & \scriptsize 6 & \scriptsize 1 & \scriptsize 2 & \scriptsize 3 & \scriptsize 4 & \scriptsize 5 & \scriptsize 6 \\
\cline{2-13}
6 & 30 & 24 & 23 & 21 & 20 & 18 & 15 & 14 & 14 & 13 & 12 & 12\\
7 & 41 & 34 & 33 & 30 & 29 & 26 & 19 & 17 & 17 & 15 & 16 & 16\\
8 & 56 & 46 & 47 & 42 & 40 & 36 & 23 & 21 & 21 & 21 & 22 & 23\\
\hline
\end{tabular}
\caption{\label{tab:exp_5_k100} Example V. $\kappa=10^{-2}$:
number of iterations for  ML- and EML-preconditioned GMRES.
}
\end{table}

\begin{table}[h!b!t!]
\centering
\begin{tabular}{|c|c|c|c|c|c|c||c|c|c|c|c|c|}
\hline
& \multicolumn{6}{c||}{ML with GMRES} & \multicolumn{6}{c|}{EML
with GMRES}\\
\cline{2-13}
\!\!\! $N$ \!\!\!\! &  \multicolumn{6}{c||}{\!\!\scriptsize  $\p$\!\!} &  \multicolumn{6}{c|}{\!\!\scriptsize $\p$\!\!}\\
\hline
& \scriptsize 1 & \scriptsize 2 & \scriptsize 3 & \scriptsize 4 & \scriptsize 5 & \scriptsize 6 & \scriptsize 1 & \scriptsize 2 & \scriptsize 3 & \scriptsize 4 & \scriptsize 5 & \scriptsize 6 \\
\cline{2-13}
6 & 46 & 21 & 19 & 16 & 16 & 15 & 23 & 13 & 13 & 12 & 11 & 12\\
7 & 55 & 29 & 27 & 24 & 23 & 22 & 26 & 17 & 17 & 16 & 16 & 16\\
8 & 61 & 43 & 41 & 36 & 35 & 32 & 31 & 24 & 24 & 23 & 24 & 24\\
\hline
\end{tabular}
\caption{\label{tab:exp_5_k1000} Example V. $\kappa=10^{-3}$:
number of iterations for  ML- and EML-preconditioned GMRES.}
\end{table}

\begin{table}[h!b!t!]
\centering
\begin{tabular}{|c|c|c|c|c|c|c||c|c|c|c|c|c|}
\hline
& \multicolumn{6}{c||}{ML with GMRES} & \multicolumn{6}{c|}{EML
with GMRES}\\
\cline{2-13}
\!\!\! $N$ \!\!\!\! &  \multicolumn{6}{c||}{\!\!\scriptsize  $\p$\!\!} &  \multicolumn{6}{c|}{\!\!\scriptsize $\p$\!\!}\\
\hline
& \scriptsize 1 & \scriptsize 2 & \scriptsize 3 & \scriptsize 4 & \scriptsize 5 & \scriptsize 6 & \scriptsize 1 & \scriptsize 2 & \scriptsize 3 & \scriptsize 4 & \scriptsize 5 & \scriptsize 6 \\
\cline{2-13}
6 & 139 & *  & 159 & 28  & 25  & 21 & 68 & 161 & 107 & 20 & 17 & 15\\
7 & 175 & *  & 33  & 27  & 23  & 21 & 58 & 24 & 21 & 17 & 15 & 14\\
8 & 184 & 43 & 34  & 30  & 26  & 24 & 51 & 20 & 18 & 16 & 16 & 16\\
\hline
\end{tabular}
\caption{\label{tab:exp_5_k10000} Example V. $\kappa=10^{-4}$:
number of iterations for  ML- and EML-preconditioned GMRES.}
\end{table}

\subsubsection{Example VI}
\begin{figure}[h!b!t!]
\subfigure[$\alpha = 10$]{
\includegraphics[trim=3.5cm 8cm 4cm 9.15cm,clip=true,width=0.48\textwidth]{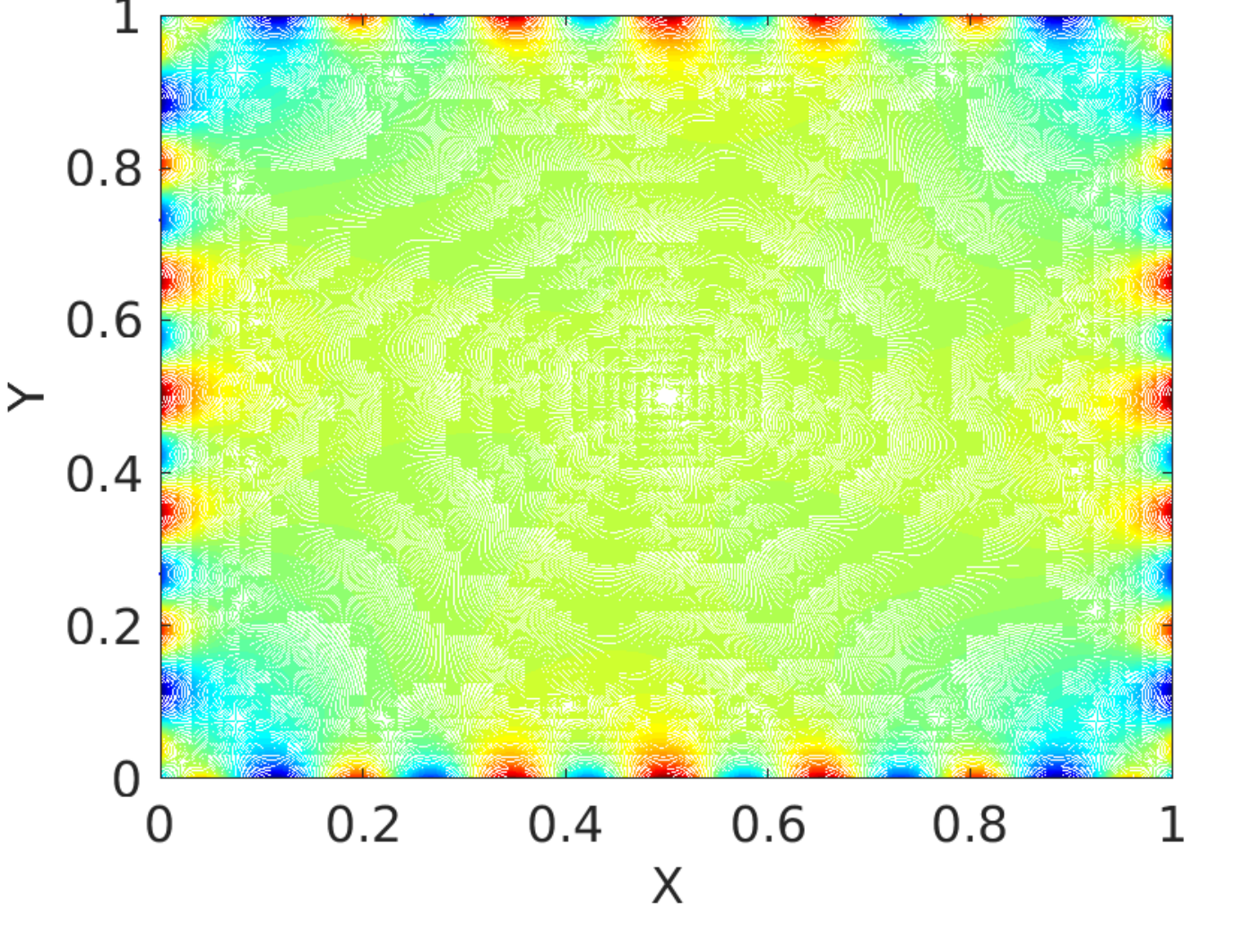}
    \label{exp6_a10}
}
    \subfigure[$\alpha = 10^{2}$]{
\includegraphics[trim=3.5cm 8cm 4cm 9.23cm,clip=true,width=0.48\textwidth]{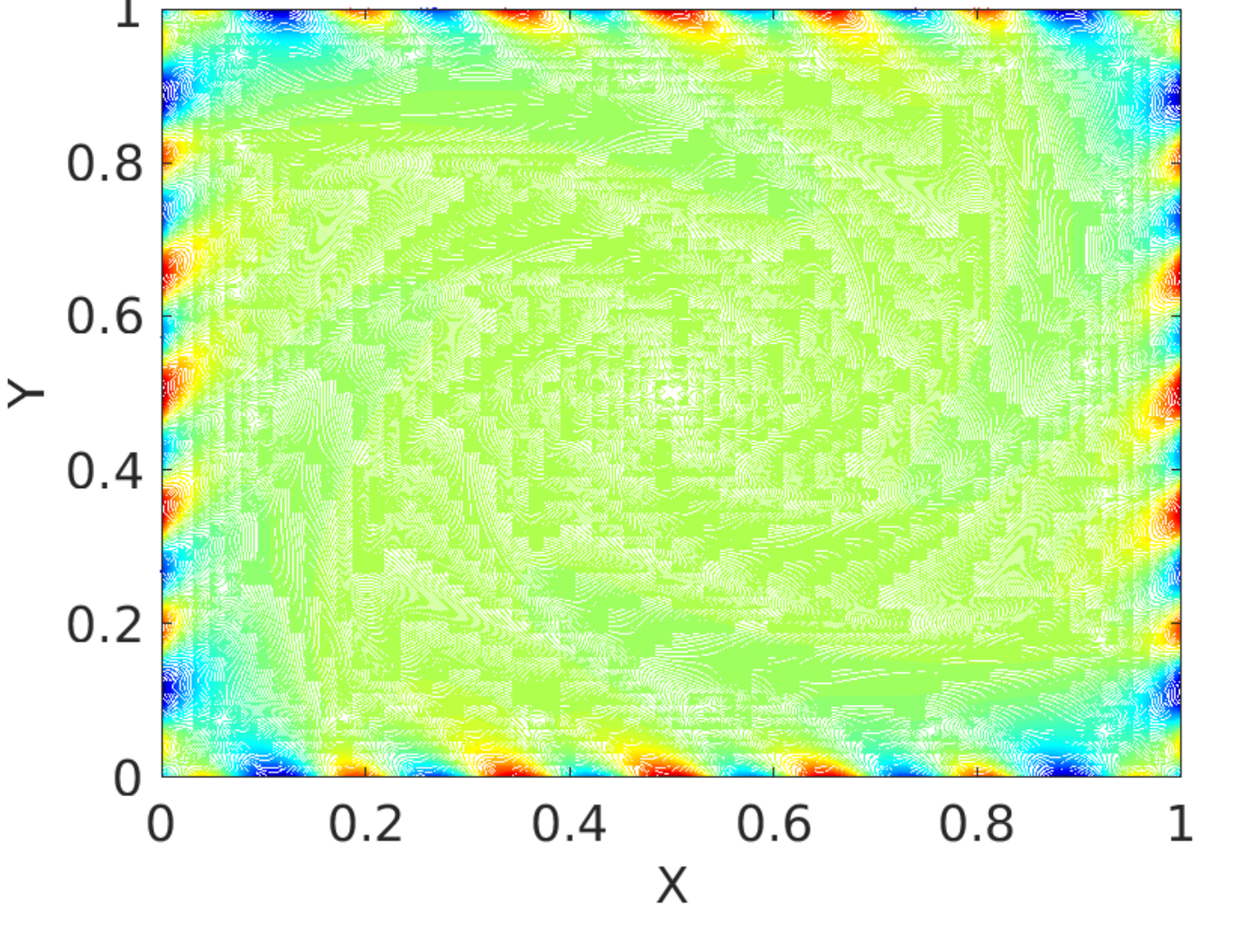}
    \label{exp6_a100}
}
    \subfigure[$\alpha = 10^{3}$]{
\includegraphics[trim=3.5cm 8cm 4cm 9.15cm,clip=true,width=0.48\textwidth]{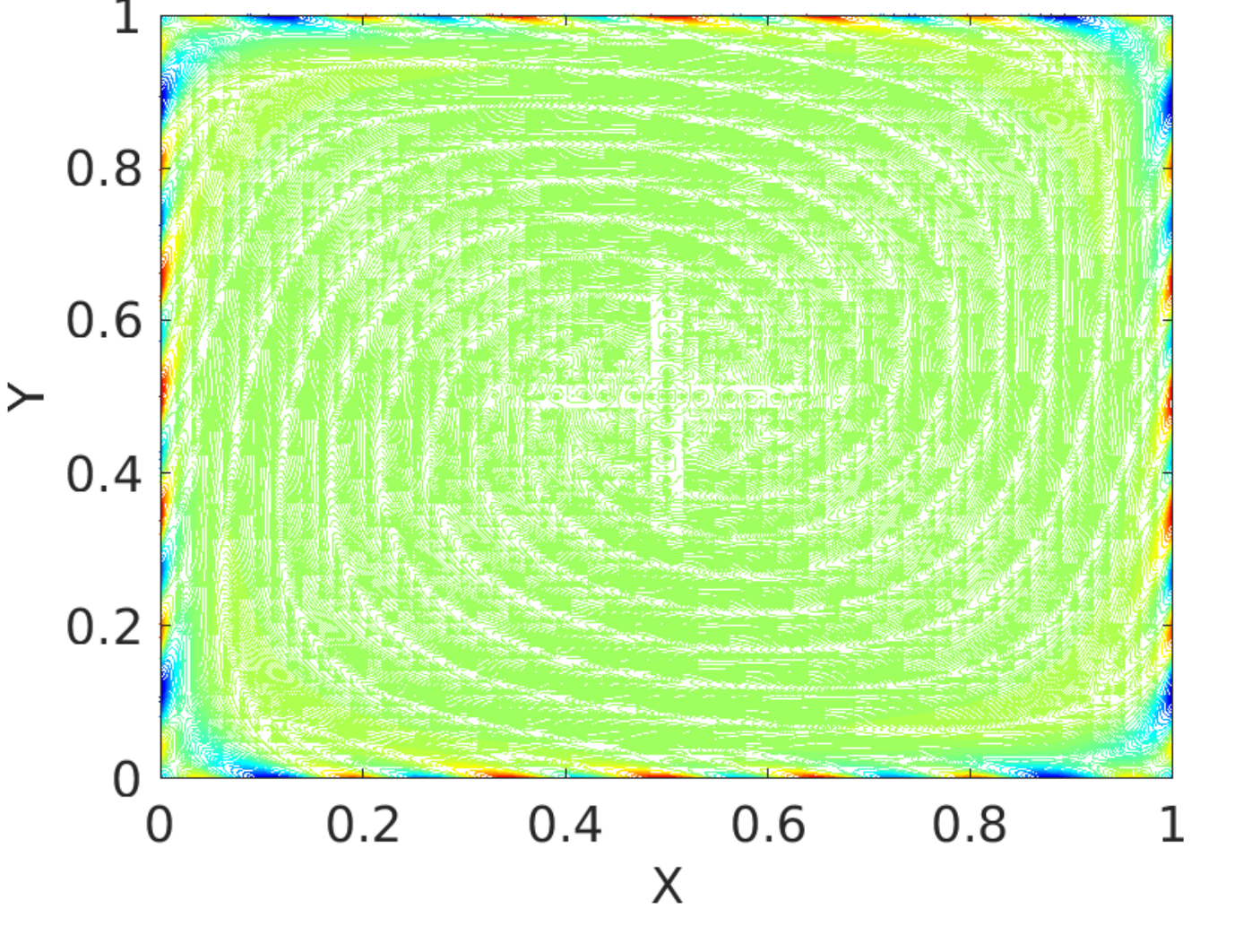}
    \label{exp6_a1000}
}
    \subfigure[$\alpha = 10^{4}$]{
\includegraphics[trim=3.5cm 8cm 4cm 9.15cm,clip=true,width=0.48\textwidth]{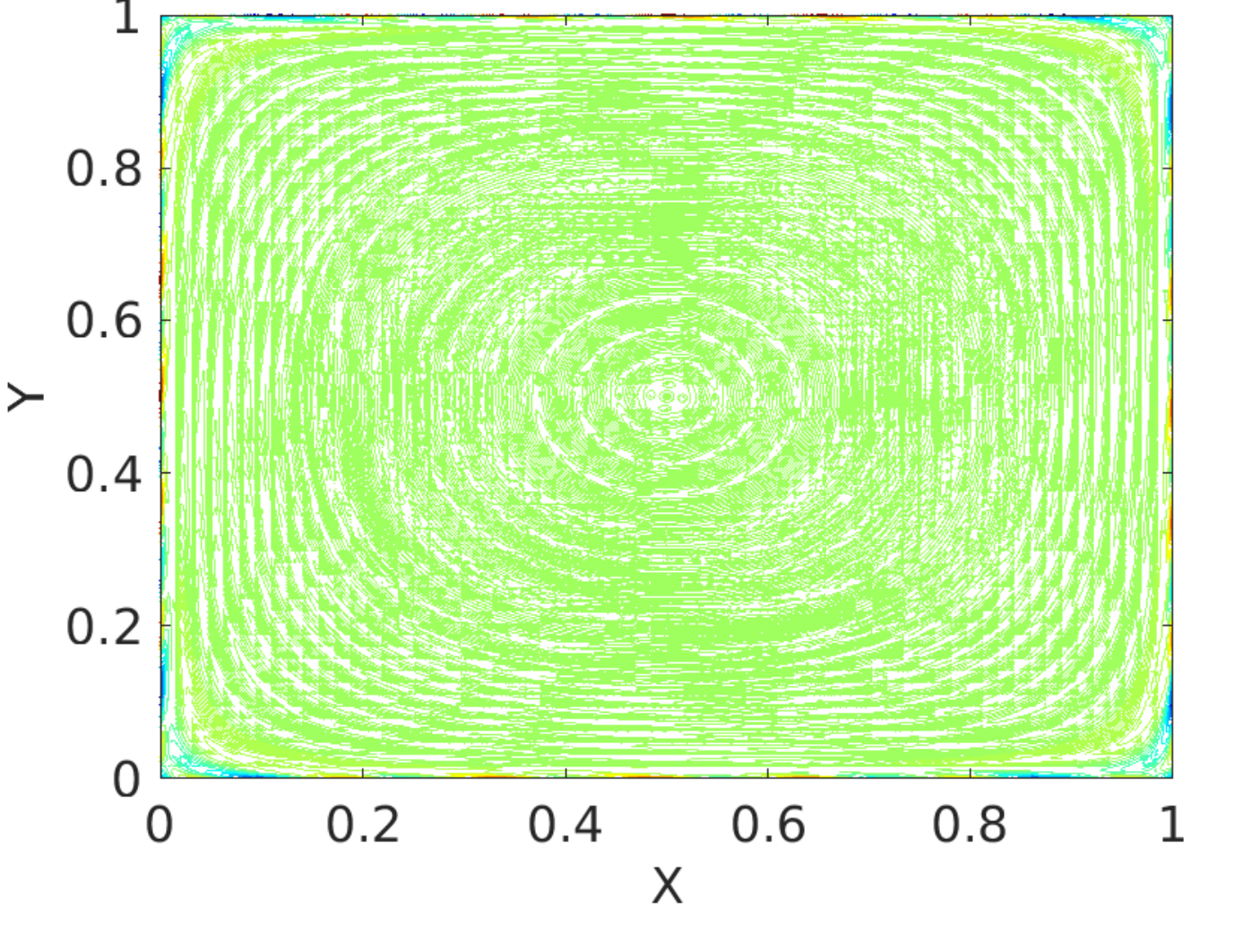}
    \label{exp6_a10000}
}

\caption{Example VI: solutions of the convection-diffusion equation for different $\alpha$ on a $64\times 64$ uniform mesh and $\p=6$ solution order.} 
\label{exp_6}
\end{figure}

Next we consider a test case from \cite{bader2003robust}. The
parameters are $\betab=(4\alpha x(x-1)(1-2y),-4\alpha y(y-1)(1-2x))$,
$g_D=\sin(\pi x)+\sin(13\pi x)+\sin(\pi y)+\sin(13\pi y)$ and ${\bf
  K}=\mc{I}$. The solution fields for four values of $\alpha$ in
$\LRs{10,10^{4}}$ are shown in Figure \ref{exp_6}. Similar to example
IV, the problem becomes more convection-dominated as $\alpha$
increases. However, one difference 
is that the
streamlines of the convection field in this case are circular
\cite{bader2003robust}. This is challenging for geometric
multigrid methods with Gauss-Seidel type smoothers if unknowns are not
ordered in the flow direction. Since the block-Jacobi method, which is
insensitive to direction, is used in ML and EML algorithms, we do not
encounter the same challenge here. In Tables
\ref{tab:exp_6_a10}-\ref{tab:exp_6_a10000} are the iteration counts of
ML- and EML-preconditioned GMRES for different values of $\alpha$. As expected, all
observations/conclusions made for example IV hold true for this
example as well. The
results for $\alpha \ge 10^3$ show that this case is, however, more
challenging. Indeed, this example requires more iterations for both ML
and EML, and in some cases convergence is not obtained within the 200-iteration constraint.

\begin{table}[h!b!t!]
\centering
\begin{tabular}{|c|c|c|c|c|c|c||c|c|c|c|c|c|}
\hline
& \multicolumn{6}{c||}{ML with GMRES} & \multicolumn{6}{c|}{EML
with GMRES}\\
\cline{2-13}
\!\!\! $N$ \!\!\!\! &  \multicolumn{6}{c||}{\!\!\scriptsize  $\p$\!\!} &  \multicolumn{6}{c|}{\!\!\scriptsize $\p$\!\!}\\
\hline
& \scriptsize 1 & \scriptsize 2 & \scriptsize 3 & \scriptsize 4 & \scriptsize 5 & \scriptsize 6 & \scriptsize 1 & \scriptsize 2 & \scriptsize 3 & \scriptsize 4 & \scriptsize 5 & \scriptsize 6 \\
\cline{2-13}
6 & 28 & 23 & 23 & 22 & 22 & 21 & 23 & 16 & 15 & 14 & 15 & 15\\
7 & 38 & 33 & 34 & 31 & 31 & 29 & 31 & 20 & 20 & 19 & 20 & 20\\
8 & 56 & 47 & 48 & 45 & 44 & 41 & 40 & 26 & 26 & 27 & 29 & 29\\
\hline
\end{tabular}
\caption{\label{tab:exp_6_a10} Example VI. $\alpha=10$: number of iterations for  ML- and EML-preconditioned GMRES.}
\end{table}

\begin{table}[h!b!t!]
\centering
\begin{tabular}{|c|c|c|c|c|c|c||c|c|c|c|c|c|}
\hline
& \multicolumn{6}{c||}{ML with GMRES} & \multicolumn{6}{c|}{EML
with GMRES}\\
\cline{2-13}
\!\!\! $N$ \!\!\!\! &  \multicolumn{6}{c||}{\!\!\scriptsize  $\p$\!\!} &  \multicolumn{6}{c|}{\!\!\scriptsize $\p$\!\!}\\
\hline
& \scriptsize 1 & \scriptsize 2 & \scriptsize 3 & \scriptsize 4 & \scriptsize 5 & \scriptsize 6 & \scriptsize 1 & \scriptsize 2 & \scriptsize 3 & \scriptsize 4 & \scriptsize 5 & \scriptsize 6 \\
\cline{2-13}
6 & 35 & 28 & 28 & 25 & 24 & 22 & 26 & 22 & 21 & 19 & 18 & 17\\
7 & 52 & 42 & 41 & 36 & 34 & 31 & 42 & 33 & 28 & 24 & 23 & 21\\
8 & 77 & 59 & 54 & 50 & 46 & 42 & 66 & 43 & 32 & 28 & 30 & 29\\
\hline
\end{tabular}
\caption{\label{tab:exp_6_a100} Example VI. $\alpha=10^2$: number of iterations for  ML- and EML-preconditioned GMRES.}
\end{table}

\begin{table}[h!b!t!]
\centering
\begin{tabular}{|c|c|c|c|c|c|c||c|c|c|c|c|c|}
\hline
& \multicolumn{6}{c||}{ML with GMRES} & \multicolumn{6}{c|}{EML
with GMRES}\\
\cline{2-13}
\!\!\! $N$ \!\!\!\! &  \multicolumn{6}{c||}{\!\!\scriptsize  $\p$\!\!} &  \multicolumn{6}{c|}{\!\!\scriptsize $\p$\!\!}\\
\hline
& \scriptsize 1 & \scriptsize 2 & \scriptsize 3 & \scriptsize 4 & \scriptsize 5 & \scriptsize 6 & \scriptsize 1 & \scriptsize 2 & \scriptsize 3 & \scriptsize 4 & \scriptsize 5 & \scriptsize 6 \\
\cline{2-13}
6 & 106 & 34 & 31 & 25 & 23 & 21 & 43  & 25 & 26 & 20 & 20 & 19\\
7 & *   & 52 & 47 & 41 & 37 & 35 & 152 & 40 & 37 & 32 & 31 & 30\\
8 & 180 & 92 & 77 & 70 & 64 & 62 & 114 & 66 & 61 & 55 & 56 & 53\\
\hline
\end{tabular}
\caption{\label{tab:exp_6_a1000} Example VI. $\alpha=10^3$: number of iterations for  ML- and EML-preconditioned GMRES.}
\end{table}

\begin{table}[h!b!t!]
\centering
\begin{tabular}{|c|c|c|c|c|c|c||c|c|c|c|c|c|}
\hline
& \multicolumn{6}{c||}{ML with GMRES} & \multicolumn{6}{c|}{EML
with GMRES}\\
\cline{2-13}
\!\!\! $N$ \!\!\!\! &  \multicolumn{6}{c||}{\!\!\scriptsize  $\p$\!\!} &  \multicolumn{6}{c|}{\!\!\scriptsize $\p$\!\!}\\
\hline
& \scriptsize 1 & \scriptsize 2 & \scriptsize 3 & \scriptsize 4 & \scriptsize 5 & \scriptsize 6 & \scriptsize 1 & \scriptsize 2 & \scriptsize 3 & \scriptsize 4 & \scriptsize 5 & \scriptsize 6 \\
\cline{2-13}
6 & 145 & 111 & 89  & 69 & 59 & 49 & 62  & 54 & 53 & 42 & 39 & 34\\
7 & 200 & 126 & 117 & 75 & 69 & 53 & 91  & 63 & 78 & 47 & 53 & 40\\
8 & * & 151 & * & 93 & 99 & 71 & * & 90 & * & 67 & 92 & 60\\
\hline
\end{tabular}
\caption{\label{tab:exp_6_a10000} Example VI. $\alpha=10^4$: number of iterations for  ML- and EML-preconditioned GMRES.}
\end{table}

Examples IV, V and VI show that both ML and EML preconditioners behave especially 
well in both diffusion-dominated and moderately convection-dominated
regimes. EML is more beneficial than ML in terms of robustness and
iteration counts, especially for low orders $p \leq 4$.

\section{Conclusion}
\seclab{conclusion} We have proposed a multilevel framework for HDG
discretizations exploiting the concepts of nested dissection, domain
decomposition, and high-order and variational structure of HDG
methods. The chief idea is to create coarse solvers for domain
decomposition methods by controlling the front growth of nested
dissection. This is achieved by  projecting the skeletal data
at different levels to either same or high-order polynomial on a set
of increasingly $h-$coarser edges/faces. When the same polynomial
order is used for the projection we name the method {\em multilevel (ML) algorithm}
and  {\em enriched multilevel (EML) algorithm} for higher
polynomial orders.  The coarse solver is combined with a block-Jacobi
fine scale solver to construct a two-level solver in the context of
domain decomposition methods. We show that the two-level approach can
also be interpreted as a multigrid algorithm with specific intergrid
transfer and smoothing operators on each level. Our complexity
estimates show that the cost of the multilevel algorithms is somewhat
in between the cost of nested dissection and standard multigrid
solvers.

We have conducted several numerical experiments with Poisson equation,
transport equation, and convection-diffusion equation in both
diffusion- and convection-dominated regimes. The numerical experiments
show that our algorithms are robust even for transport equation with
discontinuous solution and elliptic equation with highly
heterogeneous and discontinuous permeability. For convection-diffusion
equations the multilevel algorithms are scalable and reasonably robust
(with respect to changes in parameters of the underlying PDE) from diffusion-dominated to moderately
convection-dominated regimes.
EML is more beneficial than ML in terms of robustness and
iteration counts, especially for low orders $p \leq 4$.


We have demonstrated the applicability of our algorithms both as
iterative solvers and as preconditioners for various
prototypical PDEs in this work. One of the advantages of the 
algorithms is that they are designed not to depend on the nature of the PDE
being solved, but only on the smoothness  of the solution. Ongoing work is to study the performance of these algorithms for
wide variety of system of PDEs including, but not limited to, Stokes,
Navier-Stokes, and magnetohydrodynamics equations. Part of future research focuses on improving the algorithms for strongly convected and hyperbolic systems with discontinuous solutions to improve the iteration count and obtain $h-$optimal scaling.  

\section*{Acknowledgements}
 We are indebted to Professor Hari Sundar for sharing his high-order
finite element library {\bf homg}, on which we have implemented the multilevel algorithms and produced numerical results. We thank Dr. Tim Wildey for various fruitful discussions on this topic. The first author would like to thank Stephen Shannon for his generous help regarding the derivation of
projections used in this paper and also Shinhoo Kang for various fruitful discussions on this topic. 

\section*{References}

\bibliographystyle{elsarticle-num}

\bibliography{references,ceo,DtNbib}

\end{document}